\documentclass[11pt]{amsart}
\usepackage[letterpaper,margin=1.08in]{geometry}

\usepackage{amsmath,amssymb,mathtools,amsthm}
\usepackage{booktabs,array,tabularx,longtable}
\usepackage{microtype}
\usepackage[svgnames]{xcolor}
\usepackage{hyperref}
\usepackage{enumitem}
\usepackage{float}
\usepackage{xurl}

\hypersetup{
  colorlinks=true, allcolors=Navy, breaklinks=true,
  pdftitle={Exact Differential Spectra for Five Families of Power Functions},
  pdfauthor={Xiang Fan},
  pdfsubject={Power-function differential spectra and differential cryptanalysis},
  pdfkeywords={character sums, cryptography, differential spectra, finite fields, power functions},
  hypertexnames=false
}

\allowdisplaybreaks

\newtheorem{theorem}{Theorem}[section]
\newtheorem{proposition}[theorem]{Proposition}
\newtheorem{lemma}[theorem]{Lemma}
\newtheorem{corollary}[theorem]{Corollary}
\theoremstyle{definition}
\newtheorem{definition}[theorem]{Definition}
\theoremstyle{remark}
\newtheorem{remark}[theorem]{Remark}

\title{Exact Differential Spectra for Five Families of Power Functions}
\author{Xiang Fan}
\address{School of Mathematics, Sun Yat-sen University, Guangzhou 510275, China}
\email{fanx8@mail.sysu.edu.cn}

\begin{document}
\begin{abstract}
For a power function over a finite field, differential uniformity records
only the largest entry of its differential distribution table, whereas
the differential spectrum records every multiplicity and its frequency.
We determine exact differential spectra for five families of
low-uniformity power functions in odd characteristic, including families
previously studied primarily through uniformity bounds.  Our methods
combine spectrum-preserving exponent changes, factorial moments,
intersections of injective images, character sums, and point counts on
curves over finite fields.  For the ternary family of Yuan, Du, Qiao,
and Zhou, we give a fixed-cover formula and an explicit Artin
factorization for one pair-intersection cover; its odd-degree counts
reduce to traces of elliptic and genus-two curves.  Inverse-cubic
powers yield a conic evaluation in characteristic five and low-genus
trace formulas in larger characteristics.  Four ternary branches
associated with Dobbertin, Mills, M\"uller, Pott, and Willems share one
complete spectrum governed by an elliptic trace.  For the Yan--Han
family, an Artin--Prym decomposition gives all multiplicity frequencies
and exact differential uniformity, including a planar endpoint.
For Felke's family, we express extension-field spectra as
one-dimensional character sums while identifying the known
prime-field cases.  A standalone program independently reproduces the
finite-field calculations used in the proofs.
\end{abstract}
\maketitle

\section{Introduction}\label{sec:intro}

For a function $F:\mathbb F_q\to\mathbb F_q$, define
\[
 \delta_F(a,b):=\#\{x\in\mathbb F_q:F(x+a)-F(x)=b\}.
\]
The \emph{differential distribution table} (DDT) consists of these
counts.  For a nonzero input difference $a$ and a uniformly distributed
input, the differential probability of output difference $b$ is
$\delta_F(a,b)/q$.  The largest entry in a nonzero row defines the
\emph{differential uniformity}
$\Delta(F):=\max_{a\ne0,b}\delta_F(a,b)$.
For a power function $F(x)=x^d$, scaling the input identifies every
nonzero row with the row $a=1$, up to a permutation of target values.
Its \emph{differential spectrum} is therefore the sequence
\[
 \omega_j:=\#\{b\in\mathbb F_q:\delta_F(1,b)=j\}.
\]
Thus $\omega_j$ counts output differences occurring with probability
$j/q$.  The maximum alone does not determine the distribution of
these probabilities.  Complete spectra refine the differential
information used in cryptanalysis and also arise in the study of
finite-field codes and related structures
\cite{BihamShamir1991,Nyberg1994,Carlet2021,BlondeauCanteautCharpin2010,TangDingXiong2020}.

Several power families of interest have been characterized by their
uniformity without a complete account of the differential frequencies.
Yuan, Du, Qiao, and Zhou proved that a ternary permutation family has
uniformity three at the smallest degree and at most six at all larger
odd degrees; for inverse-cubic powers, they obtained
$5\leqslant\Delta\leqslant6$ in characteristic five and
$\Delta\leqslant8$ in characteristics greater than five
\cite[Section~4 and Theorems~3.2, 3.4, and~4.6]{YuanEtAl2026}.
Dobbertin, Mills, M\"uller, Pott, and Willems introduced two ternary
exponent families with uniformity at most two \cite{DMMP2003};
Leducq established equality outside the planar endpoints
\cite{Leducq2012}.  Yan and Han proved a bound of three for a further
ternary family \cite[Theorem~1]{YanHan2019}.  Felke proved uniformity
three for $p\equiv1\pmod4$ and odd extension degree, except for the
prime-field case $(p,n)=(17,1)$, where the value is two
\cite[Theorem~2.1]{Felke2020}.  The prime-field spectrum is already
known through the work of Jiang, Li, Li, and Qu
\cite[Theorem~2]{JiangLiLiQu2022}.  We retain these precise
attributions and determine the additional spectral data indicated
below.

The common objective is to recover every frequency $\omega_j$, not
merely the maximum nonzero multiplicity.  We use two complementary
counting principles.  The factorial moments
$\sum_b(\delta_F(1,b))_r$, where
$(z)_r:=\prod_{i=0}^{r-1}(z-i)$, enumerate ordered tuples of
pairwise distinct inputs belonging to a common fiber.
Alternatively, when a derivative decomposes into injective blocks,
its spectrum is determined by intersections of their target images.
In the families considered here, these counts become finite-field
character sums or point counts on curves.  Such geometric reductions
are used only with explicit exceptional-fiber and boundary
corrections.

Table~\ref{tab:application-roles} summarizes the original bounds,
the exact outputs of this paper, and the principal reductions.
\begin{table}[H]
\caption{Previously known uniformity results, exact spectral results, and main reductions.}
\label{tab:application-roles}
\centering
\footnotesize
\begin{tabularx}{\columnwidth}{@{}p{0.115\columnwidth}p{0.19\columnwidth}p{0.29\columnwidth}X@{}}
\toprule
Family & Previous result & Spectrum obtained here & Main reduction\\
\midrule
YDQZ & $\Delta\leqslant6$ & complete spectrum for odd extension degrees, via counts on curves fixed over $\mathbb F_3$ & inverse power map and three image blocks\\
inverse cubic & $5\leqslant\Delta\leqslant6$ for $p=5$; $\Delta\leqslant8$ for $p>5$ & complete spectra for $p=5$ and $p>5$ & exceptional fibers; moments of orders two or three\\
DMMP & $\Delta=2$ outside planar endpoints & common complete spectrum for all four branches & Kummer reduction; one elliptic trace\\
Yan--Han & $\Delta\leqslant3$ & complete spectrum and exact uniformity, including the endpoint & image intersections and Artin decomposition\\
Felke & $\Delta=3$ except $(17,1)$; prime-field spectrum known & spectrum as character sums for $p\equiv1\pmod4$, odd $n$ & image projectors; known prime-field spectrum\\
\bottomrule
\end{tabularx}
\end{table}

The results differ in how explicitly the remaining arithmetic is
resolved.  For the Yuan--Du--Qiao--Zhou family, the exact spectrum is
given by fixed-cover counts, and a pair-intersection count is further
reduced to an elliptic trace and a genus-two trace.  For the inverse-
cubic and Dobbertin--Mills--M\"uller--Pott--Willems families, the
necessary moments or block intersections admit explicit low-genus
trace expressions.  For Yan--Han, the configuration counts are
reconstructed through an Artin decomposition, with a planar
small-degree case treated separately.  For Felke, the full
extension-field answer is expressed in one-dimensional character
sums; the known prime-field formula is not claimed as new.  These
forms are exact, but an expression in fixed-cover point counts or
character sums is not, by itself, a closed list of numerical values
uniform in the field size.

The methods rest on two distinctions.  Replacing a power permutation
by an inverse exponent preserves its differential spectrum; sharing
a \emph{Kummer cover}, obtained by normalizing equations
$u^h=R(x)$ with $h$ prime to the characteristic, need not.
Likewise, formal moment identities determine a spectrum only after
the relevant support or rank conditions have been established.
The algebraic framework is developed in
Sections~\ref{sec:framework} and~\ref{sec:screen}; the five families
are treated in Sections~\ref{sec:ydqz}, \ref{sec:cover},
\ref{sec:dmmp}, \ref{sec:yanhan-screen}, and~\ref{sec:felke}.
Section~\ref{sec:elliptic} discusses a known-spectrum comparison
and the limitations of transferring these reductions.

These contributions build on established work on finite-field power
functions, character sums, monodromy, and curves
\cite{YanEtAl2022,AubryHerbaut2018,AubryHerbautVoloch2019,AubryHerbautIssa2023}.
The novelty claimed here is in the stated family-specific spectral
computations, common trace descriptions, and geometric factorizations,
not in the underlying classical point-counting identities.

\section{Recovering Spectra from Moments and Image Intersections}\label{sec:framework}

\subsection{From fiber counts to factorial-moment inversion}

Let $k=\mathbb F_q$, $k_m=\mathbb F_{q^m}$, and let $\pi:X\to Y$ be a finite separable morphism of degree $n$ between smooth projective geometrically connected curves.  Fiber points are counted without ramification multiplicity.  Write $\nu_m^{\rm pr}(y):=\#\pi^{-1}(y)(k_m)$.  For Frobenius-stable finite boundary sets $D_X,D_Y$, let $X^\circ:=\pi^{-1}(Y\setminus D_Y)\setminus D_X$ and denote the corresponding affine-core count by $\nu_m^\circ(y)$.  A finite-field function may further redefine finitely many pole or boundary inputs; its row count is denoted $N_m(y)$.

\begin{lemma}[Finite exceptional-input correction]\label{lem:finite-correction}
If the finite row differs from the affine core only on a finite set $E_m$, and $r_m(y)$ (resp. $s_m(y)$) counts exceptional inputs inserted at $y$ (resp. removed from the core over $y$), then
\[
 N_m(y)=\nu_m^\circ(y)+r_m(y)-s_m(y).
\]
Hence only targets represented among those finitely many removed or inserted values can change multiplicity.
\end{lemma}
\begin{proof}
The nonexceptional inputs agree; the last two terms remove the old exceptional values and insert the finite-function values.
\end{proof}

For $r\geqslant0$, let $\operatorname{Conf}_r(X/Y)$ be the ordered $r$-fold fiber product with all pairwise diagonals removed, with $\operatorname{Conf}_0(X/Y):=Y$.  Using the falling factorials defined above, put
\begin{align*}
 M_r(m)&:=\sum_{y\in Y(k_m)}(\nu_m^{\rm pr}(y))_r,\\
 \omega_j(m)&:=\#\{y:\nu_m^{\rm pr}(y)=j\}.
\end{align*}

\begin{theorem}[Exact configuration-moment identity]\label{thm:configuration-moment}
For every $m\geqslant1$ and $r\geqslant0$,
\[
 M_r(m)=\#\operatorname{Conf}_r(X/Y)(k_m).
\]
This includes branched fibers and requires no ramification correction.
\end{theorem}
\begin{proof}
Over each $y$, the configuration fiber is exactly the set of ordered $r$-tuples of distinct rational points in $\pi^{-1}(y)$, of cardinality $(\nu_m^{\rm pr}(y))_r$.
\end{proof}

\begin{theorem}[Exact spectrum inversion]\label{thm:exact-inversion}
For $0\leqslant j\leqslant n$,
\begin{align*}
 \omega_j(m)
 &=\frac1{j!}\sum_{r=j}^n\frac{(-1)^{r-j}}{(r-j)!}M_r(m)\\
 &=\frac1{j!}\sum_{r=j}^n\frac{(-1)^{r-j}}{(r-j)!}
   \#\operatorname{Conf}_r(X/Y)(k_m).
\end{align*}
\end{theorem}
\begin{proof}
For $\mathcal E_m(Z)=\sum_j\omega_j(m)Z^j$, one has $\mathcal E_m^{(r)}(1)=M_r(m)$; Taylor expansion at $Z=1$ and coefficient extraction give the formula.
\end{proof}

For a finite-type $k$-scheme $W$, put $Z(W,T):=\exp(\sum_{m\geqslant1}\#W(k_m)T^m/m)$.  The zeta functions of $\operatorname{Conf}_0,\ldots,\operatorname{Conf}_n$ together determine the spectra over all extensions, by the preceding two theorems.
For a smooth projective curve $C/k$, its \emph{Weil numerator}
$L_C(T)$ is determined by
\[
 Z(C,T)=\frac{L_C(T)}{(1-T)(1-qT)}.
\]
Here $L_C(T)=\det(1-T\operatorname{Frob}_q^*\mid
H^1(C_{\overline{k}},\mathbb Q_\ell))$, where $\ell\nmid q$
and $\operatorname{Frob}_q^*$ is pullback by the $q$-power
Frobenius endomorphism; its $m$th power has trace
$q^m+1-\#C(k_m)$.

\begin{theorem}[Spectral rationality and Frobenius recurrence]\label{thm:spectral-recurrence}
For a fixed finite separable cover, each sequence $(\omega_j(m))_{m\geqslant1}$ has a rational ordinary generating function and hence satisfies a finite linear recurrence over $\mathbb Q$.  More precisely,
\begin{align*}
 \sum_{m\geqslant1}\omega_j(m)T^m
 ={}&\frac1{j!}\sum_{r=j}^n\frac{(-1)^{r-j}}{(r-j)!}\\
 &\times T\frac{d}{dT}\log Z(\operatorname{Conf}_r(X/Y),T).
\end{align*}
If the relevant zeta functions factor into Weil or Artin factors, the reduced recurrence roots are among their Frobenius eigenvalues, up to cancellation.  The same conclusion holds after any fixed finite Frobenius-stable boundary and exceptional-input datum.
\end{theorem}
\begin{proof}
Dwork rationality gives rational $Z(W,T)$ for finite-type $W$ \cite{Dwork1960}, and
\[
 T\frac{d}{dT}\log Z(W,T)=\sum_{m\geqslant1}\#W(k_m)T^m.
\]
Combine this with Theorems~\ref{thm:configuration-moment} and \ref{thm:exact-inversion}. Fixed finite boundary data contribute only finite Frobenius-orbit counting sequences, whose generating functions are rational.
\end{proof}

\subsection{Recovering multiplicities from image intersections}

Moment inversion is not the only exact reconstruction mechanism.  Suppose the regular input set has already been quotiented by a visible symmetry and decomposed into labeled blocks whose target maps are injective on quotient orbits.

\begin{proposition}[Injective block-image reconstruction]\label{prop:block-incidence}
Let the regular inputs be partitioned into $h$-element orbits, assigned to blocks $1,\ldots,r$.  Assume that within each block no two distinct orbits have the same target.  Let $\mathcal T$ be the target universe after the exceptional targets have been removed, let $\mathcal A_i\subseteq\mathcal T$ be the target image of block $i$, and let $N_k$ be the number of targets in $\mathcal T$ belonging to exactly $k$ of the sets $\mathcal A_i$.  Put $I_0:=|\mathcal T|$ and, for $1\leqslant j\leqslant r$,
\[
 I_j:=\sum_{\substack{J\subseteq\{1,\ldots,r\}\\|J|=j}}
 \left|\bigcap_{i\in J}\mathcal A_i\right|.
\]
Then every regular fiber has size $hk$ for some $0\leqslant k\leqslant r$, and
\[
 I_j=\sum_{k=j}^r\binom{k}{j}N_k,
 \qquad
 N_k=\sum_{j=k}^r(-1)^{j-k}\binom{j}{k}I_j.
\]
Restoring the exceptional targets then gives the complete spectrum from the required image intersections.
\end{proposition}
\begin{proof}
A target lying in exactly $k$ block images has one quotient orbit from each of those blocks and hence $hk$ regular preimages.  It contributes $\binom{k}{j}$ to $I_j$.  The first identity follows by summing over targets, and the second is binomial inversion.
\end{proof}

This route can be strictly cheaper than constructing ordered self-fiber products of a single large cover.  In the YDQZ family, $h=2$ and $r=3$; blockwise injectivity reduces the complete spectrum to three pair intersections and one triple intersection of fixed rational maps.

\subsection{Support bounds, monodromy, and rank}

\begin{proposition}[Moments with restricted fiber support]\label{prop:few-support}
If every fiber size lies in a known set $S=\{s_1,\ldots,s_t\}$, then $M_0(m),\ldots,M_{t-1}(m)$ determine the spectrum.  The coefficient matrix $((s_i)_r)_{0\leqslant r<t,\,1\leqslant i\leqslant t}$ is invertible, with determinant, up to sign, $\prod_{i<j}(s_j-s_i)$.
\end{proposition}
\begin{proof}
The basis $1,z,(z)_2,\ldots,(z)_{t-1}$ is obtained from $1,z,\ldots,z^{t-1}$ by a unitriangular change of basis, so its evaluation determinant is Vandermonde.
\end{proof}

\begin{proposition}[Rank criterion for factorial moments]\label{prop:relation-aware-stopping}
Let the spectrum be supported on $S=\{s_1,\ldots,s_t\}$, and suppose exact linear relations $C\boldsymbol\omega=c$ are already known for $\boldsymbol\omega=(\omega_{s_1},\ldots,\omega_{s_t})^T$.  Put $v_r:=((s_i)_r)_{i=1}^t$.  The least depth for full-rank linear reconstruction over $\mathbb Q$ is the least $R\in\{-1,0,\ldots,t-1\}$ such that
\[
 \operatorname{rank}\!\begin{pmatrix}C\\v_0\\ \vdots\\v_R\end{pmatrix}=t,
\]
where the matrix is $C$ alone when $R=-1$.
At that depth the linear system uniquely determines the spectrum.
Below full rank, its equations have more than one rational solution;
additional nonnegativity or integrality constraints may nevertheless
uniquely determine a particular spectrum.
\end{proposition}
\begin{proof}
The known relations together with $M_r=v_r\boldsymbol\omega$ form a linear system in the $t$ supported coefficients.  Full column rank is equivalent to uniqueness; below full rank a nonzero homogeneous direction remains.
\end{proof}

\begin{remark}[Existence of a fiber of prescribed size]\label{rem:goal-adaptive}
To prove the occurrence of a fiber of size at least $r$, it suffices to produce one admissible rational point on the ordered $r$-configuration after boundary restoration.  If a fiber-size upper bound of $r$ is already known, this gives an $r$-fiber.  Exact spectrum recovery requires only the point counts selected by the reduced support matrix.  Thus a genus bound, together with a matching upper bound, may settle differential uniformity before exact Frobenius traces are known.
\end{remark}

Ordinary moments carry the same information through Stirling inversion.  The full double self-fiber product is stratified by the diagonal and the ordered distinct-pair locus $\operatorname{Conf}_2(X/Y)$.  For a finite-field function $F:\mathbb F_q\to\mathbb F_q$, define its four-variable differential energy by
\[
 E_4(F):=\sum_{a,b\in\mathbb F_q}\delta_F(a,b)^2.
\]
For a power function, the row $a=0$ contributes $q^2$ and every nonzero row is a target rescaling of the row $a=1$; hence
\begin{equation}\label{eq:global-energy}
 E_4(F)=q^2+(q-1)\sum_{b\in\mathbb F_q}\delta_F(1,b)^2.
\end{equation}
Thus the usual square moment is the point count of the full double self-fiber product, while its off-diagonal stratum is the second ordered configuration.

The \emph{branch locus} $B$ consists of the points of $Y$
over which $\pi$ is not \'etale.  Put $U:=Y\setminus B$ and
$V:=\pi^{-1}(U)$.  Let $L$ be the \emph{Galois closure} of
$k(X)/k(Y)$, its smallest normal overfield in a fixed separable
closure, and put $k_L:=L\cap\overline{k}$.
The \emph{arithmetic monodromy group} and \emph{geometric monodromy group}
are, respectively,
\[
 A:=\operatorname{Gal}(L/k(Y)),\qquad G:=\operatorname{Gal}(L/k_L(Y)).
\]
Then $G\triangleleft A\le S_\Omega$ on the geometric generic fiber $\Omega$, with $A/G\cong\operatorname{Gal}(k_L/k)$.  Let $A_m$ be the coset restricting to $q^m$-Frobenius on $k_L$.

\begin{proposition}[Frobenius-coset support]\label{prop:frobenius-support}
For $y\in U(k_m)$,
\[
 \nu_m^{\rm pr}(y)\in\{\#\operatorname{Fix}_\Omega(\gamma):\gamma\in A_m\}.
\]
\end{proposition}
\begin{proof}
The finite \`etale fiber is a Frobenius set; its $k_m$-points are the sheets fixed by its arithmetic Frobenius element in $A_m$.
\end{proof}
This is an exact support inclusion; it does not assert that every conjugacy class in $A_m$ occurs over a fixed field.

Let $\Omega^{(r)}$ be the ordered distinct $r$-tuples in $\Omega$.
\begin{proposition}[Orbit--component correspondence]\label{prop:orbit-component}
The geometrically irreducible components of $\operatorname{Conf}_r(V/U)\times_k\overline{k}$ are indexed by the $G$-orbits on $\Omega^{(r)}$.
\end{proposition}
\begin{proof}
The geometric generic fiber is $\Omega^{(r)}$, and connected components of a finite \`etale cover are the orbits of geometric monodromy.
\end{proof}
If $h_r(m)$ is the number of those components fixed by $q^m$-Frobenius, twisted Burnside gives
\[
 h_r(m)=\frac1{|G|}\sum_{\gamma\in A_m}(\#\operatorname{Fix}_\Omega(\gamma))_r,
\]
and therefore $\#\operatorname{Conf}_r(V/U)(k_m)=h_r(m)q^m+O(q^{m/2})$ for fixed $r$; exact spectra still require complete point counts and boundary corrections.

\begin{proposition}[Invariant-coordinate principle]\label{prop:quotient}
If a finite $H\subset\operatorname{PGL}_2(k)$ fixes $\Phi\in k(x)$, then $\Phi$ factors through $\mathbb P^1\to\mathbb P^1/H$.  Any generator of $k(x)^H$ is a quotient coordinate; generic fibers are unions of $H$-orbits, with deviations only at quotient or residual branch points.
\end{proposition}
\begin{proof}
The inclusion $k(\Phi)\subseteq k(x)^H$ gives the factorization; the remaining assertions are the standard quotient description.
\end{proof}

\subsection{Selecting the necessary configurations}\label{subsec:protocol}

Known fiber sizes and exact relations reduce the factorial moments
needed for reconstruction through
Proposition~\ref{prop:relation-aware-stopping}.  For derivatives with
injective input blocks, Proposition~\ref{prop:block-incidence} instead
uses intersections of target images.  Both methods distinguish
projective and affine fibers from the finitely many exceptional inputs.
The exponent relations below can simplify the required point counts.

\section{Exponent Equivalence and Kummer Covers}\label{sec:screen}

Changing the exponent can simplify the derivative of a power map.  Frobenius changes and inversion of power permutations preserve the complete spectrum.  A common Kummer presentation instead organizes the character sums used to count fibers and need not preserve the spectrum.  We separate these two relations before constructing the configuration covers.

\subsection{Equivalent power exponents}

Let $q=p^r$, $N=q-1$, and $F_d(x):=x^d$.  Exponent congruences modulo $N$ are interpreted first on $\mathbb F_q^\times$; the value at $x=0$ is retained as an explicit finite exceptional datum whenever a chosen representative has a pole or a zero residue.  Frobenius exponents $d\mapsto p^jd\pmod N$ preserve the spectrum.  Inverse exponents do so as well when $F_d$ is a permutation.

\begin{proposition}[Inverse exponents of power permutations]\label{prop:inverse-screening}
If $\gcd(d,N)=1$ and $d^*$ is the inverse of $d$ modulo $N$, then $x^d$ and $x^{d^*}$ have the same complete differential spectrum.  Thus every exponent in
\[
 \mathcal O_{\rm sp}(d):=\{p^jd,\ p^jd^*:0\leqslant j<r\}\pmod N
\]
gives the same spectrum and may be used in its place.
\end{proposition}
\begin{proof}
For every permutation $F$, $\delta_{F^{-1}}(a,b)=\delta_F(b,a)$.  For a power permutation, scaling identifies all nonzero input rows, and $\delta_F(1,0)=0$; as $b\ne0$ varies, $b^{-d}$ varies bijectively.  Thus DDT transposition preserves the row multiset.
\end{proof}

\subsection{Shared Kummer covers and character twists}

A different relation preserves geometry rather than spectrum.  Suppose $p\nmid h$, set $g=\gcd(h,N)$, and let $hd_0\equiv e\pmod N$.  The corresponding set of exponents is
\[
 \mathcal P_{h,e}:=\left\{d_a=d_0+a\frac{N}{g}:a\in\mathbb Z/g\mathbb Z\right\},
 \qquad \kappa_g(x):=x^{N/g}.
\]
Then, on $\mathbb F_q^\times$,
\begin{equation}\label{eq:kummer-packet-identity}
 x^{d_a}=\kappa_g(x)^a x^{d_0},\qquad (x^{d_a})^h=x^e.
\end{equation}
These exponents satisfy the same $h$-th-power equation, although their power functions can select different branches.

\begin{theorem}[Character-twist configuration transfer]\label{thm:kummer-transfer}
Fix a factorial-moment depth $j$.  All exponents in $\mathcal P_{h,e}$ admit the same ambient $h$-Kummer equations $u_i^h=x_i^e$ and $v_i^h=(x_i+1)^e$.  On the stratum $\kappa_g(x_i)=\alpha_i$, $\kappa_g(x_i+1)=\beta_i$, passing from $d_0$ to $d_a$ changes only the deck marking
\[
 (u_i,v_i)\longmapsto(\alpha_i^a u_i,\beta_i^a v_i).
\]
Consequently their factorial moments are finite character projections on the same Kummer cover and are linear combinations of its character sums.
\end{theorem}
\begin{proof}
For $x_i,x_i+1\ne0$, identity~\eqref{eq:kummer-packet-identity}
gives $\alpha_i,\beta_i\in\mu_g\subseteq\mu_h$.  Thus the
displayed scaling preserves the Kummer equations and changes each
derivative to $\beta_i^av_i-\alpha_i^au_i$.  Expanding the
indicators of the character strata by Fourier inversion on $\mu_g$
and summing gives the stated projections.  Inputs at $0$ or $-1$
are handled by the finite exceptional-input correction.
\end{proof}

The choice of character projector can change both the fiber sizes and the moments needed for reconstruction.

\subsection{Choosing the arithmetic model}

For a finite separable extension $K/k$, the \emph{norm-one torus}
is the kernel of the norm morphism
$\operatorname{Res}_{K/k}\mathbb G_m\to\mathbb G_m$.
Frobenius exponents and inverse exponents of power permutations can
yield simpler rational expressions for the derivative.  Kummer,
character--Frobenius, and norm-one-torus presentations may simplify
point counts, but a common Kummer cover need not preserve the spectrum.
After separating projective, affine, and exceptional finite fibers,
monodromy or uniformity bounds restrict the support; known relations
then determine the necessary moments by
Proposition~\ref{prop:relation-aware-stopping}.  The resulting spectrum
satisfies $\sum_j\omega_j=\sum_j j\omega_j=q$ and can be checked by
direct enumeration over a fixed field.

\paragraph{Comparison of the two ternary exponents.}
For odd $n$, put $q=3^n$, $d_{\rm YH}=(q-3)/4$ and $d_{\rm YDQZ}=(3q-5)/4$ \cite{YanHan2019,YuanEtAl2026}.  Then
\begin{align*}
 4d_{\rm YH}&\equiv4d_{\rm YDQZ}\equiv-2\pmod{q-1},\\
 d_{\rm YDQZ}-d_{\rm YH}&=\frac{q-1}{2}.
\end{align*}
Thus their fourth powers agree on $\mathbb F_q^\times$, while the exponents differ by a quadratic-character twist.  This relation does not imply equal spectra.  Independently, $d_{\rm YDQZ}$ is a permutation exponent with
\begin{align*}
 d_{\rm YDQZ}^{-1}&\equiv\frac{q-5}{2}\pmod{q-1},\\
 x^{d_{\rm YDQZ}^{-1}}&=\chi(x)x^{-2}.
\end{align*}
Section~\ref{sec:ydqz} uses the inverse exponent: its derivative splits into three injective blocks whose pair and triple image intersections determine the spectrum.  The Yan--Han argument instead uses the quartic Kummer cover and an Artin decomposition of a genus-$81$ curve.  Inversion preserves the spectrum, whereas sharing the quartic Kummer equation does not.

\begin{theorem}[Spectrum equivalence and moment reconstruction]\label{thm:screen-exactness}
Frobenius exponent changes and inversion of a power permutation
preserve the differential spectrum.  Kummer character projections
recover the corresponding factorial moments exactly.  If the required
moments and exceptional-fiber corrections are known and the
relation-augmented moment matrix has full column rank, the complete
spectrum is uniquely determined.  Over any fixed finite field it can
also be obtained by exhaustive enumeration.
\end{theorem}
\begin{proof}
The spectrum invariances follow from Frobenius and
Proposition~\ref{prop:inverse-screening}; the character projections
are Theorem~\ref{thm:kummer-transfer}.  Apply
Proposition~\ref{prop:relation-aware-stopping} and linear inversion.
Direct enumeration is finite.
\end{proof}

\section{A Ternary Permutation Family: Inverse Exponents and Image Intersections}\label{sec:ydqz}

The YDQZ exponent varies with the extension degree and belongs to
one quartic Kummer family with the Yan--Han exponent.  Inverting the
YDQZ permutation instead gives the fixed character-rational function
$\chi(x)x^{-2}$.  An involution and three injective blocks then reduce
its spectrum to intersections of fixed rational images.  The proof
establishes the exact parametrization degrees $4,8,8$ and the
multiplicities $32,32,64,256$, including all boundary corrections.

\subsection{Exact spectrum in terms of fixed covers}

Let $q=3^n$ with $n$ odd, let $k=\mathbb F_q$, and let $\chi:k\to\{0,\pm1\}$ be the quadratic character.  Put
\[
 d:=\frac{3q-5}{4},\qquad
 U:=k\setminus\{0,1,-1\}.
\]
Define rational functions over $\mathbb F_3$ by
\begin{equation}\label{ydqz:eq:fixed-maps}
\begin{aligned}
 f_A(t)&:=\frac{t^6(t^4+1)}{(t^4-1)^4},\\
 f_B(t)&:=\frac{t^{16}+t^{10}+t^6+1}{(t^4+1)^4},\\
 f_C(t)&:=\frac{t^{16}+t^{10}+t^6+1}{t^4(t^2-1)^4}.
\end{aligned}
\end{equation}
For the theorem, the three functions in \eqref{ydqz:eq:fixed-maps} are restricted to $U$.  Define the fixed open fiber products by their $k$-points:
\begin{align*}
 X_{AB}^{\circ}(k)&:=\{(r,s)\in U^2:f_A(r)=f_B(s)\},\\
 X_{AC}^{\circ}(k)&:=\{(r,s)\in U^2:f_A(r)=f_C(s)\},\\
 X_{BC}^{\circ}(k)&:=\{(r,s)\in U^2:f_B(r)=f_C(s)\},\\
 Z_U^{\circ}(k)&:=\{(r,s,t)\in U^3:f_A(r)=f_B(s)=f_C(t)\}.
\end{align*}
These are the rational points of fixed open curves over $\mathbb F_3$; only the extension field varies with $n$.

Set
\begin{align}\label{ydqz:eq:ThetaXi}
 \Theta_n^U
 &:=5q-2\#X_{AB}^{\circ}(k)-2\#X_{AC}^{\circ}(k)
   -\#X_{BC}^{\circ}(k),\notag\\
 \Xi_n^U&:=q-\#Z_U^{\circ}(k).
\end{align}
For $b\in k$, write
\begin{align*}
 \delta_F(1,b)&=\#\{x\in k:F(x+1)-F(x)=b\},\\
 \omega_j&=\#\{b\in k:\delta_F(1,b)=j\}.
\end{align*}

\begin{theorem}[Complete YDQZ fixed-cover spectrum]\label{ydqz:thm:main}
For every odd positive integer $n$, the complete differential spectrum of
\[
 F(x)=x^{(3^{n+1}-5)/4}\quad\text{over }\mathbb F_{3^n}
\]
is supported on $\{0,2,3,4,6\}$ and is given by
\begin{align*}
 \omega_0&=\frac{147q+128-4\Theta_n^U+\Xi_n^U}{256},\\
 \omega_2&=\frac{91q-384+8\Theta_n^U-3\Xi_n^U}{256},\\
 \omega_3&=1,\\
 \omega_4&=\frac{17q-4\Theta_n^U+3\Xi_n^U}{256},\\
 \omega_6&=\frac{q-\Xi_n^U}{256},
\end{align*}
with $\omega_j=0$ for every other $j$.
\end{theorem}

The permutation property and the upper bound are due to Yuan, Du,
Qiao, and Zhou \cite[Section~4 and Theorem~4.6]{YuanEtAl2026}.
Their Lemma~4.2 also gives the fibers over $0,-1,1$ as $0,0,3$;
the degree-one endpoint is already included in their theorem.
Our contribution is the complete multiplicity formula through fixed
image intersections.  The inverse representative below gives an
independent proof of the same upper bound.

\subsection{Inverse exponent and three image blocks}

\begin{lemma}[Inverse representative]\label{ydqz:lem:inverse}
The power map $F(x)=x^d$ is a permutation, and its complete differential spectrum equals that of
\[
 G(0):=0,\qquad G(x):=\chi(x)x^{-2}\quad(x\ne0).
\]
\end{lemma}

\begin{proof}
Since $n$ is odd, $q\equiv3\pmod8$; write $q=8m+3$.  Then
\begin{align*}
 d&=6m+1, & e&=\frac{q-5}{2}=4m-1,\\
 de-1&=(8m+2)(3m-1)=(q-1)(3m-1).
\end{align*}
Thus $e\equiv d^{-1}\pmod{q-1}$, and $F$ is a permutation whose inverse agrees with $x^e$ on $k^\times$.  By Proposition~\ref{prop:inverse-screening}, their complete differential spectra coincide.  Since $e=(q-1)/2-2$, this inverse is $\chi(x)x^{-2}$ for $x\ne0$ and is zero at $x=0$.
\end{proof}

Henceforth put
\[
 D(x):=G(x+1)-G(x).
\]
Since $q\equiv3\pmod4$, $-1$ is a nonsquare and $\chi(-1)=-1$.  In particular $G(-x)=-G(x)$ for every $x\in k$.

\begin{lemma}[Involution and block formulas]\label{ydqz:lem:blocks}
One has
\[
 D(-1-x)=D(x).
\]
For $x\notin\{0,-1\}$, put
\[
 s:=\chi(x),\qquad t:=\chi(x+1),\qquad u:=x(x+1).
\]
Then
\[
 D(x)=\frac{t x^2-s(x+1)^2}{u^2},
\]
and the regular inputs split into the following three blocks:
\begin{align*}
 A:&\quad s=t, & D(x)&=\frac{s(x-1)}{u^2};\\
 B:&\quad (s,t)=(1,-1), & D(x)&=\frac{u-1}{u^2};\\
 C:&\quad (s,t)=(-1,1), & D(x)&=\frac{1-u}{u^2}.
\end{align*}
The involution $\iota(x):=-1-x$ exchanges the two sign strata in $A$ and preserves $B$ and $C$.  Its unique fixed point is $x=1$.
\end{lemma}

\begin{proof}
Because $G$ is odd,
\begin{align*}
 D(-1-x)&=G(-x)-G(-1-x)\\
 &=-G(x)+G(x+1)=D(x).
\end{align*}
The displayed formula follows by putting both terms over the denominator $x^2(x+1)^2=u^2$ and using characteristic three.  The transformation of the signs is
\[
 (s,t)\longmapsto(-t,-s),
\]
because $\chi(-1)=-1$.  The remaining assertions are immediate.
\end{proof}

\subsection{The exceptional three-point fiber}

\begin{proposition}[Special fiber]\label{ydqz:prop:special}
For every odd $n$,
\[
 D(x)=1\iff x\in\{0,-1,1\}.
\]
Consequently $b=1$ is the unique odd-sized fiber and has size three.
\end{proposition}

\begin{proof}
Directly, $D(0)=D(-1)=D(1)=1$.  Suppose $x\notin\{0,-1\}$.  In block $B$, the equation $D(x)=1$ gives $(u+1)^2=0$, hence $u=-1$ and $(x-1)^2=0$, so $x=1$.  In block $C$, it gives $u^2+u-1=0$, whose discriminant $-1$ is nonsquare in $k$; thus there is no solution.

In block $A$, $D(x)=1$ gives
\[
 x^2(x+1)^2=s(x-1),\qquad s\in\{\pm1\}.
\]
Thus $x$ is a zero of one of
\begin{align*}
 P_+(X)&=X^4-X^3+X^2-X+1,\\
 P_-(X)&=X^4-X^3+X^2+X-1.
\end{align*}
The values of $P_+$ at $0,1,-1$ are $1,1,-1$, and those of $P_-$ are $-1,1,1$.  Thus neither has a linear factor over $\mathbb F_3$.  Every irreducible factor of a quartic without a linear factor has degree two or four, so its roots cannot lie in $\mathbb F_{3^n}$ when $n$ is odd.
\end{proof}

Removing the three inputs $0,-1,1$, the involution is fixed-point-free.  Thus the remaining $q-3$ inputs form exactly $(q-3)/2$ regular two-element orbits.

\subsection{Parametrizing the regular image blocks}

\subsubsection{Block A}

\begin{proposition}[Block $A$: degree four and target injectivity]\label{ydqz:prop:A}
Restricted to $U$, the map $f_A$ is exactly four-to-one onto the set of regular derivative targets contributed by block $A$.  Distinct $A$-orbits have distinct derivative targets.
\end{proposition}

\begin{proof}
Let $\tau\in U$ and put
\[
 a=\tau-\tau^{-1},\qquad c=\tau+\tau^{-1},\qquad x=a^2.
\]
Then
\[
 c^2-a^2=1,
\]
so $x$ and $x+1=c^2$ are nonzero squares.  Hence $x$ is the square-square representative of a regular $A$-orbit.  Moreover
\[
 u=x(x+1)=\frac{(\tau^4-1)^2}{\tau^4},
 \qquad
 x-1=\frac{\tau^4+1}{\tau^2},
\]
and therefore
\[
 D(x)=\frac{x-1}{u^2}=f_A(\tau).
\]
Conversely, each $A$-orbit has a unique square-square representative $x$.  Choosing square roots $a^2=x$, $c^2=x+1$ gives four sign pairs $(a,c)$.  For each pair
\[
 \tau=-(a+c)
\]
is the unique element satisfying $a=\tau-\tau^{-1}$ and $c=\tau+\tau^{-1}$.  Thus every $A$-orbit has exactly four parameters in $U$.

It remains to show that two different $A$-orbits cannot have the same target.  For $T,R\in U$, cross-multiplication gives, in $\mathbb F_3[T,R]$,
\begin{equation}\label{ydqz:eq:A-factor}
\begin{split}
&T^6(T^4+1)(R^4-1)^4-R^6(R^4+1)(T^4-1)^4\\
&\qquad=(R-T)(R+T)(RT-1)(RT+1)Q_A(T,R),
\end{split}
\end{equation}
where, writing
\[
 X=T^2,\quad Y=R^2,\quad S=X+Y,\quad \Pi=XY,
\]
one has the sum-of-two-squares identity
\begin{equation}\label{ydqz:eq:A-SOS}
 Q_A(T,R)
 =\bigl(S^2(\Pi-1)\bigr)^2
  +\bigl((X-Y)(\Pi+1)^2\bigr)^2.
\end{equation}
If $Q_A(T,R)=0$, then both squared terms in \eqref{ydqz:eq:A-SOS} vanish, because $-1$ is nonsquare in $k$.  The second term gives either $X=Y$ or $\Pi=-1$.  The latter is impossible, since $\Pi=(TR)^2$ is a square.  Hence $R^2=T^2$, so $R=\pm T$.

Therefore \eqref{ydqz:eq:A-factor} shows
\[
 f_A(T)=f_A(R)
 \quad\Longrightarrow\quad
 R\in\{\pm T,\pm T^{-1}\}.
\]
These four parameters are precisely the four sign choices attached to one $A$-orbit; they are distinct because $T\ne0$, $T^2\ne1$, and $T^2\ne-1$.  This proves both assertions.
\end{proof}

\subsubsection{Blocks B and C}

Put
\[
 h(u):=\frac{u-1}{u^2}.
\]
For every regular input in $B$ or $C$,
\begin{equation}\label{ydqz:eq:opposite-quotient}
 \chi(u)=-1,
 \qquad
 u+1=x^2+x+1=(x-1)^2\in k^{\times2}.
\end{equation}

\begin{lemma}[Opposite-sign quotient injectivity]\label{ydqz:lem:h-inj}
The function $h$ is injective on the set of $u\in k^\times$ satisfying \eqref{ydqz:eq:opposite-quotient}.
\end{lemma}

\begin{proof}
If $h(u)=h(v)$, then
\[
 (u-1)v^2-(v-1)u^2
 =-(u-v)(uv-u-v)=0.
\]
If $u\ne v$, then
\[
 v=\frac{u}{u-1}.
\]
Using \eqref{ydqz:eq:opposite-quotient} and $\chi(-1)=-1$ gives
\begin{align*}
 \chi(v)&=-\chi(u-1),\\
 \chi(v+1)&=\chi\!\left(-\frac{u+1}{u-1}\right)
 =-\chi(u-1).
\end{align*}
Thus $v$ cannot itself satisfy $\chi(v)=-1$ and $\chi(v+1)=1$.  Hence $u=v$.
\end{proof}

\begin{proposition}[Block $B$: degree eight]\label{ydqz:prop:B}
Restricted to $U$, the map $f_B$ is exactly eight-to-one onto the regular derivative targets contributed by block $B$, and distinct $B$-orbits have distinct targets.
\end{proposition}

\begin{proof}
The conic $\mathcal C_-:a^2+c^2=-1$ has the bijective parametrization
\[
 a=\frac{t^2+t-1}{t^2+1},\qquad
 c=\frac{1+t-t^2}{t^2+1},\qquad t\in\mathbb P^1(k),
\]
obtained from lines through $(1,1)$, with $t=\infty$ corresponding to $(1,-1)$.  The denominator never vanishes because $-1$ is nonsquare.
For a regular $B$-input choose $a^2=x$ and $c^2=-(x+1)$.  Each input has four sign choices, and each regular involution orbit contains two inputs.  The four parameters $t=0,1,-1,\infty$ lie over the special fixed input $x=1$; thus exactly eight parameters in $U$ represent each regular $B$-orbit.

For these parameters,
\[
 ac=-\frac{t^4+1}{(t^2+1)^2},\qquad
 u=x(x+1)=-(ac)^2,
\]
and simplification in characteristic three gives
\[
 h(u)=\frac{u-1}{u^2}
 =\frac{t^{16}+t^{10}+t^6+1}{(t^4+1)^4}=f_B(t).
\]
Again $t^4+1\ne0$.  The quotient $u$ determines the involution orbit, because the two roots of $X^2+X-u$ are $x$ and $-1-x$.  By Lemma~\ref{ydqz:lem:h-inj}, distinct $B$-orbits have distinct targets, completing both assertions.
\end{proof}

\begin{proposition}[Block $C$: degree eight]\label{ydqz:prop:C}
Restricted to $U$, the map $f_C$ is exactly eight-to-one onto the regular derivative targets contributed by block $C$, and distinct $C$-orbits have distinct targets.
\end{proposition}

\begin{proof}
The conic $\mathcal C_+:a^2+c^2=1$ has the bijective parametrization
\[
 a=\frac{t}{1+t^2},\qquad
 c=\frac{1-t^2}{1+t^2},\qquad t\in\mathbb P^1(k),
\]
from lines through $(0,1)$, with $t=\infty$ corresponding to $(0,-1)$.  For a regular $C$-input take $a^2=-x$ and $c^2=x+1$.  As for block $B$, four sign choices for each input and two inputs per involution orbit give eight parameters.  Here $t=0,\pm1,\infty$ are exactly the boundary parameters over $x=0,-1$, so every regular orbit has eight parameters in $U$.

The quotient and its target satisfy
\[
 u=x(x+1)=-(ac)^2=-\frac{t^2(1-t^2)^2}{(1+t^2)^4}
\]
and, using the polynomial identity
\[
 t^{16}+t^{10}+t^6+1
 =(t^2+1)^4(t^8-t^6+t^4-t^2+1)
 \quad\text{in }\mathbb F_3[t],
\]
we obtain
\[
 -h(u)=\frac{1-u}{u^2}
 =\frac{t^{16}+t^{10}+t^6+1}{t^4(t^2-1)^4}=f_C(t).
\]
Lemma~\ref{ydqz:lem:h-inj}, applied after negating the target, proves injectivity on $C$-orbits.
\end{proof}

\begin{remark}[Boundary parameters]\label{ydqz:rem:boundary}
The restriction to $U$ is essential and exactly accounts for the normalization in the fixed-cover formula.  The affine expression $f_A$ has an extraneous value at $t=0$, whereas the hyperbolic parametrization uses $t\in k^\times$.  In block $B$, the four projective parameter values over the fixed input $x=1$ are $0,1,-1,\infty$; the three finite ones are removed by $U$.  In block $C$, the four projective boundary parameters $0,1,-1,\infty$ lie over $x=0$ or $x=-1$.  No regular parameter is lost and no boundary parameter remains.
\end{remark}

\subsection{Fiber-product counts and the spectrum}

Let $\mathcal A,\mathcal B,\mathcal C\subset k$ be the three sets of regular derivative targets contributed by blocks $A,B,C$, respectively.  Propositions~\ref{ydqz:prop:A}--\ref{ydqz:prop:C} imply that every regular target contains at most one involution orbit from each block.  Hence every regular fiber has size $0,2,4$, or $6$.  Together with Proposition~\ref{ydqz:prop:special}, this already gives
\[
 \operatorname{supp}(\omega)\subseteq\{0,2,3,4,6\}
\]
and recovers $\Delta(F)\leqslant6$.

Put
\[
 P_{AB}:=|\mathcal A\cap\mathcal B|,
 \quad
 P_{AC}:=|\mathcal A\cap\mathcal C|,
 \quad
 P_{BC}:=|\mathcal B\cap\mathcal C|,
\]
\[
 P:=P_{AB}+P_{AC}+P_{BC},
 \qquad
 T:=|\mathcal A\cap\mathcal B\cap\mathcal C|.
\]
There are $(q-3)/2$ regular involution orbits in total.  Therefore the total incidence count is
\[
 |\mathcal A|+|\mathcal B|+|\mathcal C|=\frac{q-3}{2}.
\]
Elementary inclusion-exclusion gives the spectrum directly.

\begin{proposition}[Incidence reconstruction]\label{ydqz:prop:incidence}
One has
\begin{align*}
 \omega_6&=T,\\
 \omega_4&=P-3T,\\
 \omega_2&=\frac{q-3}{2}-2P+3T,\\
 \omega_0&=\frac{q+1}{2}+P-T,\\
 \omega_3&=1.
\end{align*}
\end{proposition}

\begin{proof}
On $\mathcal T=k\setminus\{1\}$, the three regular images have
$I_0=q-1$, $I_1=(q-3)/2$, $I_2=P$, and $I_3=T$.
Proposition~\ref{prop:block-incidence} with $h=2$ gives
$\omega_6=T$, $\omega_4=P-3T$,
$\omega_2=I_1-2P+3T$, and $\omega_0=I_0-I_1+P-T$.
Proposition~\ref{ydqz:prop:special} restores the unique three-fiber
at $b=1$, proving the formulas.
\end{proof}

The exact parametrization degrees determine the four normalization factors in the fixed fiber products.

\begin{corollary}[Fixed-cover multiplicities]\label{ydqz:cor:multiplicities}
With the open fiber products defined above,
\begin{align*}
 \#X_{AB}^{\circ}(k)&=32P_{AB},\\
 \#X_{AC}^{\circ}(k)&=32P_{AC},\\
 \#X_{BC}^{\circ}(k)&=64P_{BC},\\
 \#Z_U^{\circ}(k)&=256T.
\end{align*}
Consequently
\begin{align*}
 P&=\frac{\#X_{AB}^{\circ}(k)+\#X_{AC}^{\circ}(k)}{32}
   +\frac{\#X_{BC}^{\circ}(k)}{64},\\
 T&=\frac{\#Z_U^{\circ}(k)}{256}.
\end{align*}
\end{corollary}

\begin{proof}
Every common $A$--$B$ target has exactly $4\cdot8=32$ parameter pairs; the $A$--$C$ case is identical.  Every common $B$--$C$ target has $8\cdot8=64$ pairs, and every triple target has $4\cdot8\cdot8=256$ parameter triples.  Remark~\ref{ydqz:rem:boundary} guarantees that the open domain $U$ contains precisely the regular parameters and no exceptional ones.
\end{proof}

\begin{proof}[Proof of Theorem~\ref{ydqz:thm:main}]
By Corollary~\ref{ydqz:cor:multiplicities} and \eqref{ydqz:eq:ThetaXi},
\[
 P=\frac{5q-\Theta_n^U}{64},\qquad
 T=\frac{q-\Xi_n^U}{256}.
\]
Proposition~\ref{ydqz:prop:incidence} gives the stated five spectrum
coefficients.  Lemma~\ref{ydqz:lem:inverse} transfers the spectrum
from $G$ to the original exponent.
\end{proof}

\paragraph{Finite-field checks.}
The \texttt{--historical} option of
\nolinkurl{anc/NewDiffSpec_arXiv1_verify.py} independently enumerates
the original power derivative, its inverse-power representative, and
all four open-cover counts for $n=1,3,5,7,9$.  All three routes agree
with the spectrum above and with the historical rows embedded
in the same verification script.  These finite checks support
the formulas; the preceding proof applies to every odd degree.
The resulting uniformities also agree with
\cite[Table~3]{YuanEtAl2026} on these five fields.

\begin{remark}[Arithmetic form of the result]
The spectrum is expressed in terms of rational-point counts of fixed open covers, with normalization factors $32,32,64,256$.  Theorem~\ref{thm:spectral-recurrence} gives rational generating functions and finite Frobenius recurrences for these counts.  Proposition~\ref{ydqz:prop:BC-explicit} below evaluates the $B$--$C$ intersection by an elliptic trace and a genus-two trace.  The other two pair intersections and the triple intersection retain their fixed-cover expressions.
\end{remark}

\subsection{An explicit factorization of the \texorpdfstring{$B$--$C$}{B--C} intersection}

For a polynomial $P(T)=\prod_j(1-\lambda_jT)$ with constant term one,
write $\tau_n(P):=\sum_j\lambda_j^n$.  Define
\[
 E(T):=1-2T+3T^2,\qquad
 H(T):=1+2T+2T^2+6T^3+9T^4.
\]
These are the numerator polynomials of the smooth projective curves
\[
 E_0:y^2=-(z+1)(z^2+1),\qquad
 H_0:y^2=-z(z^2+1)(z^2+z-1),
\]
of genera one and two, respectively.

\begin{proposition}[Explicit $B$--$C$ intersection]\label{ydqz:prop:BC-explicit}
For every odd positive integer $n$, with $q=3^n$,
\begin{equation}\label{ydqz:eq:BC-trace}
 \#X_{BC}^{\circ}(\mathbb F_q)=q+1-8\tau_n(E)-6\tau_n(H).
\end{equation}
The smooth projective normalization $C_{BC}$ of this fixed fiber
product is geometrically connected of genus $105$.  Its zeta numerator
has the form $E(T)^8H(T)^6R(T^2)$, where $R\in\mathbb Z[T]$ has
degree $85$ and is given by the Artin factors below.
\end{proposition}

\begin{proof}
The two conic parametrizations identify the function field with
\[
 a^2+c^2=-1,\quad u=-a^2c^2=a^4+a^2,\qquad
 A^2+C^2=1,\quad v=-A^2C^2=A^4-A^2,
\]
subject to $h(u)=-h(v)$, where $h(u)=(u-1)/u^2$.
Putting $z=v/u$ gives the rational parametrization
\begin{equation}\label{ydqz:eq:BC-base}
 u=\frac{z^2+1}{z(z+1)},\qquad
 v=\frac{z^2+1}{z+1}.
\end{equation}
Each conic extension over $\overline{\mathbb F}_3(z)$ has group
$D_4$ of order eight, generated by sign changes and interchange
of its coordinates.  The simple pole of $u$ at $0$ gives inertia
of order four in $a^4+a^2=u$; its independent quadratic subfields
have square classes $1+u$ and $-u$.  The second extension has the
corresponding classes $1+v$ and $-v$ and order-four inertia at infinity.
The two extensions have no common quadratic subfield: classes
containing $1+u$ or $1+v$ are distinguished by the roots of
$z^2-z-1$ or $z^2+z-1$, while $-u$ and $-v$ are distinguished
at $0$ and infinity.  Their intersection is therefore trivial,
since a nontrivial Galois $2$-extension has a quadratic quotient.
Thus $C_{BC}\to\mathbb P^1_z$ is geometrically connected and
Galois with group $G=D_4\times D_4$.  All these automorphisms
are defined over $\mathbb F_3$, so the constant field is $\mathbb F_3$.

The six roots of
$(z^2+1)(z^2-z-1)(z^2+z-1)$ have inertia of order two.
The points $0,-1,\infty$ have inertia of order four.
At the common branch points the composite inertia is cyclic,
with order the least common multiple of the two tame indices.
Consequently
\[
 2g(C_{BC})-2
 =64\left(-2+6\cdot\frac12+3\cdot\frac34\right)=208.
\]

Write the irreducible representations of each $D_4$ as
$1,\beta,\alpha,\alpha\beta,\sigma$, where $\alpha$ and $\beta$
correspond to $1+u$ and $-u$, and $\sigma$ is the standard
two-dimensional signed-permutation representation; use $v$ on
the second factor.  Put $P^-(T):=P(-T)$ and define
\begin{align*}
 J&:=1+3T^2,& U&:=1+9T^4,\\
 Q&:=1-2T^2+9T^4,& B&:=1-4T^2+9T^4,\\
 D&:=1+2T+4T^2+6T^3+9T^4,\\
 V&:=1-T^2+4T^3-3T^4+27T^6,\\
 X&:=1-6T^4+81T^8.
\end{align*}
The complete Artin factors, including inertia invariants at the
branch points, are listed in Table~\ref{ydqz:tab:BC-Artin}.
The entry $(1,1)$ records the numerator of the rational base.

\begin{table}[H]
\caption{Artin numerator factors for $C_{BC}\to\mathbb P^1_z$.
Each entry has multiplicity $\dim(r)\dim(s)$.}
\label{ydqz:tab:BC-Artin}
\centering
\begin{tabular}{@{}c@{\quad}ccccc@{}}
\toprule
$r\backslash s$ & $1$ & $\beta$ & $\alpha$ & $\alpha\beta$ & $\sigma$\\
\midrule
$1$ & $1$ & $E$ & $E$ & $J$ & $U$\\
$\beta$ & $E$ & $1$ & $H$ & $J$ & $V^-$\\
$\alpha$ & $E$ & $H$ & $EE^-$ & $JQ$ & $EBD$\\
$\alpha\beta$ & $J$ & $J$ & $JQ$ & $J$ & $X$\\
$\sigma$ & $U$ & $V$ & $E^-BD^-$ & $X$ & $J^2EQH$\\
\bottomrule
\end{tabular}
\end{table}

Here is a complete finite calculation of the table.  Over
$\mathbb F_q$, put $\varepsilon:=\chi(-1)$ and
\[
 S_B(u):=\#\{a:a^4+a^2=u\}-1-\chi(1+u),\qquad
 S_C(v):=\#\{A:A^4-A^2=v\}-1-\chi(1+v).
\]
At an unramified $z$, put
$a_0:=\chi(1+u)$, $b_0:=\chi(-u)$,
$c_0:=\chi(1+v)$, and $d_0:=\chi(-v)$.
The five local traces are
$(1,b_0,a_0,a_0b_0,S_B)$ and
$(1,d_0,c_0,c_0d_0,S_C)$; their products give the tensor traces.
Indeed, the permutation representation on the four $a$-coordinates
is $1+\alpha+\sigma$.
For $r=\alpha^e\beta^f$, $s=\alpha^g\beta^h$ with
$e,f,g,h\in\{0,1\}$, the ramified traces are as follows;
$[\mathcal P]$ is one or zero according as $\mathcal P$ holds or fails.

\begin{center}
\footnotesize
\begin{tabular}{@{}ccccc@{}}
\toprule
Place & $r\otimes s$ & $\sigma\otimes s$ & $r\otimes\sigma$ & $\sigma\otimes\sigma$\\
\midrule
$u=v=0$ & $[f=h]$ & $1$ & $1$ & $2$\\
$u=-1$ & $[e=0]c_0^gd_0^h$ & $c_0^gd_0^h$ & $[e=0]S_C$ & $S_C$\\
$v=-1$ & $[g=0]a_0^eb_0^f$ & $[g=0]S_B$ & $\varepsilon a_0^eb_0^f$ & $\varepsilon S_B$\\
$z=0$ & $[e=f]\varepsilon^{e+g+h}$ & $0$ & $[e=f]\varepsilon^e S_C(1)$ & $0$\\
$z=\infty$ & $[g=h]\varepsilon^{e+f+g}$ & $[g=h]\varepsilon^g S_B(1)$ & $0$ & $0$\\
$z=-1$ & $[e+f+g+h\text{ even}]\varepsilon^{f+g}$ & $0$ & $0$ & $1+\varepsilon$\\
\bottomrule
\end{tabular}
\end{center}

For the simultaneous zero $z_0^2=-1$, write
$u=u_1\pi+\cdots$, $v=v_1\pi+\cdots$, where $\pi=z-z_0$.
Then $u_1v_1=z_0/(z_0+1)^2$.  The invariant directions of
$\sigma\otimes\sigma$ have leading squares $-1$ and $-u_1v_1$;
both are squares in $\mathbb F_9$, giving trace two whenever
the base point is rational.  The other common-zero entries follow
from the same leading units.  At $z=-1$, one has
$u\sim\pi^{-1}$ and $v\sim-\pi^{-1}$.
The branch residue ratios lie in $\mathbb F_9$, which contains
the fourth roots of $-1$.  For even $n$ Frobenius fixes both
rotation-invariant tensor directions; for odd $n$ it exchanges
them.  Their trace is therefore $1+\varepsilon$.
At the one-sided reflection points,
$(a^2-1)^2=0$ or $(A^2+1)^2=0$ gives traces $1$ or
$\varepsilon$.  At a one-sided pole the standard representation
has no inertia invariants.  These observations give every row.

The tame conductor formula
\[
 \deg L_\rho=-2\dim\rho+
       \sum_z\operatorname{codim}\rho^{I_z}
\]
gives degree at most six for the nontrivial rank-one factors,
degrees $4,6,10,8$ for each row or column involving exactly one
$\sigma$, and degree fourteen for $\sigma\otimes\sigma$.
For example, the last conductor sum is
$4+4+4+4+4+2=22$, hence its degree is $22-8=14$.
All representations here are orthogonal over $\mathbb Q$.
The alternating cup-product pairing therefore induces a
Frobenius pairing of multiplier three on each multiplicity space.
Every factor of degree $2d$ satisfies
$L(T)=3^dT^{2d}L(1/(3T))$.
Thus at most seven logarithmic trace sums determine each factor.
Summing the displayed local rules for $n=1,\ldots,7$ and using
$jc_j=\sum_{i=1}^j s_i c_{j-i}$, where
$\log L(T)=\sum_{i\geqslant1}s_iT^i/i$ and $c_0=1$,
gives Table~\ref{ydqz:tab:BC-Artin} by exact integer arithmetic.
The \texttt{--ydqz-artin} option of
\nolinkurl{anc/NewDiffSpec_arXiv1_verify.py} reproduces these sums
and the reconstruction, and compares the full point counts with a
separate conic-lift enumeration.  Extension degrees eight and nine
provide checks beyond the reconstruction data.

Multiplying the table entries with their indicated multiplicities gives
\begin{equation}\label{ydqz:eq:BC-factor}
 P_{C_{BC}}=E^{11}(E^-)^3J^{15}H^6U^4Q^6B^4X^4
                 (VV^-)^2(DD^-)^2.
\end{equation}
In \eqref{ydqz:eq:BC-factor}, removing $E^8H^6$ leaves an even polynomial of
degree $170$.  This proves the asserted factorization.
The quadratic characters $\beta$ on the second factor and
$\beta\otimes\alpha$ give the displayed models $E_0$ and $H_0$.

Finally let $n$ be odd.  The six quadratic branch points are not
$\mathbb F_q$-rational.  A point over $0,-1,$ or $\infty$ would
give a rational point at infinity on one of the conics, requiring
a square root of $-1$.  The excluded $B$ parameters lie over
$u=-1$, and the excluded $C$ parameters lie over $v=0$;
by \eqref{ydqz:eq:BC-base}, neither can occur at an odd-degree
rational base point.  Hence $X_{BC}^{\circ}$ and $C_{BC}$ have
the same $\mathbb F_q$-points.  The even factor contributes no
odd-degree trace, proving \eqref{ydqz:eq:BC-trace}.
\end{proof}

Consequently, for odd $n$, the first quantity in
\eqref{ydqz:eq:ThetaXi} may be written as
\[
 \Theta_n^U=4q-1-2\#X_{AB}^{\circ}(\mathbb F_q)
                 -2\#X_{AC}^{\circ}(\mathbb F_q)
                 +8\tau_n(E)+6\tau_n(H).
\]
The $A$--$B$, $A$--$C$, and triple intersections retain their
fixed-cover expressions.  Direct parameter counts for
$n=1,3,5,7,9$ give the $B$--$C$ values
$0,192,0,1344,19776$, respectively, and agree with
\eqref{ydqz:eq:BC-trace}.
Moreover, Hasse--Weil gives
$\#X_{BC}^{\circ}(\mathbb F_q)\geqslant q+1-40\sqrt q>0$
for odd $n\geqslant7$.  A common $B$--$C$ target contains two
distinct regular involution orbits, so the differential uniformity
is at least four in this range.  The direct checks at $n=3,5$
give the same lower bound.

\section{Inverse-Cubic Power Functions in Characteristic \texorpdfstring{$p\geqslant5$}{p>=5}}\label{sec:cover}

The exponent $q-4$ has rational core $x^{-3}$ in every characteristic,
but its branch structure depends on $p$.
Yuan, Du, Qiao, and Zhou proved the bounds
$5\leqslant\Delta\leqslant6$ for $p=5$ and $\Delta\leqslant8$ for
$p>5$ \cite[Theorems~3.2 and~3.4]{YuanEtAl2026}.
Their proof already identifies the five-fiber at $b=1$ in
characteristic five and shows that the bound eight is attained when
the extension degree is divisible by six in larger characteristics.
We determine all multiplicities and the exact characteristic-five
uniformity.  In characteristic five, two branch values coalesce and
the second factorial moment suffices; for $p>5$, the generic
four-value branch locus requires a third moment.  Both calculations
retain the ramification and boundary corrections.

\subsection{Characteristic five}\label{subsec:p5-inverse-cubic}

\begin{theorem}[Characteristic-five inverse-cubic spectrum]\label{thm:main-spectrum}
Let $q=5^m$ and $F_m(0):=0$, $F_m(x):=x^{q-4}=x^{-3}$ for $x\ne0$.  If $m$ is odd, then
\[
 (\omega_0,\omega_2,\omega_4,\omega_5)
 =\left(\frac{5q+7}{8},\frac{q-5}{4},\frac{q-5}{8},1\right),
\]
and if $m$ is even, then
\[
 (\omega_0,\omega_2,\omega_5,\omega_6)
 =\left(\frac{7(q-1)}{12},\frac{3q+5}{8},1,\frac{q-25}{24}\right),
\]
with all omitted coefficients zero.
\end{theorem}

\begin{corollary}[DDT consequences in characteristic five]\label{cor:ddt}
Every nonzero DDT row is a rescaling of the row $a=1$; its unique entry equal to five occurs at $b=a^{q-4}$.  Moreover
\[
 \Delta(F_m)=\begin{cases}
 5,&m\text{ odd or }m=2,\\
 6,&m\geqslant4\text{ even}.
 \end{cases}
\]
The maximum differential probability is $\operatorname{MDP}(F_m)=\Delta(F_m)/q$; $F_m$ is a permutation exactly for odd $m$, and
\begin{align*}
 \sum_b\delta_{F_m}(1,b)^2&=3q-10(-1)^m,\\
 E_4(F_m)&=q^2+(q-1)\bigl(3q-10(-1)^m\bigr).
\end{align*}
\end{corollary}

\begin{theorem}[Monodromy and the second moment]\label{thm:main-monodromy}
For $D(x):=(x+1)^{-3}-x^{-3}$, the geometric and arithmetic monodromy groups are respectively $S_4$ on $S_4/C_4$ and $C_2\wr S_3$, with Galois-closure constant field $\mathbb F_{25}$.  The regular fiber sizes lie in $\{0,2,4\}$ for odd $m$ and $\{0,2,6\}$ for even $m$, while the finite definition creates one five-fiber.  The only further moment needed for spectrum reconstruction is
\[
 M_2=2q-10(-1)^m.
\]
The double-configuration component not given by the involution normalizes to
\[
 S^2=2(W^2+1).
\]
\end{theorem}

Let $q=5^m$ and write
\[
 D(x)=(x+1)^{-3}-x^{-3}
\]
for the rational core, defined away from $x=0,-1$.  The involution
\[
 \iota(x)=-1-x
\]
fixes $D$.  Its quotient coordinate is $u:=x(x+1)$.  In characteristic five,
\[
 D(x)=\frac{2u+4}{u^3}.
\]
Writing $t=u^{-1}$ gives the factorization
\begin{equation}\label{eq:factorization}
 \mathbb P^1_x\xrightarrow{\ \kappa\ }\mathbb P^1_t
 \xrightarrow{\ P\ }\mathbb P^1_b,
\end{equation}
where
\[
 \kappa(x):=\frac1{x(x+1)},\qquad P(t):=t^2(2-t).
\]
The two maps have degrees two and three.

\subsubsection{Ramification}

The involution quotient $\kappa$ is ramified over $t=0,1$.  The cubic $P$ is ramified at $t=0,3,\infty$, with critical values $0,1,\infty$.  The complete ramification data are shown in Table~\ref{tab:ramification}.

\begin{table}[H]
\caption{Ramification of the factorization \eqref{eq:factorization}.}\label{tab:ramification}
\centering
\footnotesize
\begin{tabular}{@{}cccc@{}}
\toprule
$b$ & Points in the $t$-fiber & Ramification in $P$ & Composite partition \\
\midrule
$0$ & $t=0,2$ & $2,1$ & $(4,1,1)$ \\
$1$ & $t=1,3$ & $1,2$ & $(2,2,2)$ \\
$\infty$ & $t=\infty$ & $3$ & $(3,3)$ \\
\bottomrule
\end{tabular}
\end{table}

The total ramification contribution is
\[
 3+3+4=10=2\cdot6-2,
\]
as required by Riemann--Hurwitz.

\subsubsection{Geometric monodromy}

\begin{lemma}[The six-point $S_4$ action]\label{lem:geometric-group}
The geometric monodromy group of $D$ is the index-two subgroup
\[
 \begin{aligned}
 W^+=\{(v,\pi)\in(C_2)^3\rtimes S_3:\;&\\
 &|v|\equiv\operatorname{par}(\pi)\pmod2\},
 \end{aligned}
\]
with $W^+\subset C_2\wr S_3$, where $\operatorname{par}(\pi)=0$ for even $\pi$ and $1$ for odd $\pi$, and the three size-two blocks are the inverse images under $\kappa$ of the three roots of $P(t)-b$.  Hence
\[
 G=W^+\cong V_4\rtimes S_3\cong S_4,
\]
and its degree-six action is the coset action $S_4/C_4$.
\end{lemma}

\begin{proof}
The factorization \eqref{eq:factorization} gives a canonical system of three blocks of size two.  The induced action on the blocks is the geometric monodromy of the cubic $P(t)-b$.  Its two finite simple branch points give transpositions and the point at infinity gives a $3$-cycle, so the block action is $S_3$.

The three inertia cycle types on the six sheets are the composite partitions in Table~\ref{tab:ramification}.  The inertia elements above $b=0,1$ induce transpositions on the blocks, as shown by $P$; their respective types $(4,1,1)$ and $(2,2,2)$ therefore have odd within-block flip parity.  The inertia at infinity has type $(3,3)$, induces a $3$-cycle on the blocks, and has even flip parity.  Thus every inertia generator lies in $W^+$, so $G\subseteq W^+$ and $|G|\leqslant24$.

Let $g_0$ be inertia of type $(4,1,1)$.  Its square acts trivially on the blocks and flips exactly two of them.  Conjugating $g_0^2$ by the full block action $S_3$ produces the three nonzero even-weight vectors in $(C_2)^3$; hence the kernel of the block action contains
\[
 V=\{v\in(C_2)^3:|v|\equiv0\pmod2\}\cong V_4.
\]
Therefore $|G|\geqslant|V|\,|S_3|=24$, and $G=W^+$.
The map $\pi\mapsto(\operatorname{par}(\pi)(1,1,1),\pi)$ is
a section of the block action.  Since $S_3$ acts faithfully on the
three nonzero elements of $V$, this gives
$G\cong V_4\rtimes S_3\cong S_4$.

Finally, an inertia element of order four has two fixed sheets by the partition $(4,1,1)$.  The stabilizer of either fixed sheet has order $|G|/6=4$ and contains this order-four element, so it is cyclic.  Hence the transitive action is $S_4/C_4$.

For the explicit finite inventory below, choose a labeling of this action in which
\[
 \sigma_0=(1\ 2\ 3\ 4),\qquad
 \sigma_1=(1\ 3)(2\ 5)(4\ 6)
\]
generate $G$ and the block system is
$\{\{1,3\},\{2,4\},\{5,6\}\}$.
\end{proof}

\subsubsection{The normalizer}

\begin{lemma}[Normalizer of the six-point group]\label{lem:normalizer}
For the group $G$ of Lemma~\ref{lem:geometric-group},
\[
 N_{S_6}(G)=C_2\wr S_3,\qquad |N_{S_6}(G)|=48,
\]
\[
 N_{S_6}(G)/G\cong C_2.
\]
A representative of the nontrivial coset is
\[
 \tau=(5\ 6).
\]
\end{lemma}

\begin{proof}
In a faithful transitive action on $G/K$, the centralizer of $G$ is isomorphic to $N_G(K)/K$.  Here $K\cong C_4$, and its normalizer in $S_4$ is dihedral of order eight.  Hence the centralizer of $G$ in $S_6$ has order two.

The group $S_4$ is complete: its center is trivial and every automorphism is inner.  Conjugation by an element of $N_{S_6}(G)$ induces an automorphism of $G$, which can therefore be corrected by an element of $G$ to centralize $G$.  It follows that
\[
 |N_{S_6}(G)|\leqslant|G|\cdot2=48;
\]
see also the standard coset-action discussion in \cite{DixonMortimer1996}.

Directly,
\[
 \tau\sigma_0\tau^{-1}=\sigma_0,
 \qquad
 \tau\sigma_1\tau^{-1}
 =\sigma_0^2\sigma_1\sigma_0^2.
\]
Thus $\tau\in N_{S_6}(G)$ and $\tau\notin G$.  The group $\langle G,\tau\rangle$ has order $48$ and preserves the block system, hence is the full block stabilizer $C_2\wr S_3$.  This attains the upper bound and proves the result.
\end{proof}

\subsubsection{The constant field and arithmetic monodromy}

For a root $t$ of $P(t)-b$, the two points above it under $\kappa$ are obtained from
\[
 x^2+x=\frac1t.
\]
If $\eta=2x+1$, then
\[
 \eta^2=\frac{t-1}{t}.
\]

\begin{lemma}[The constant field]\label{lem:constant-field}
The Galois closure of $D(x)-b$ over $\mathbb F_5(b)$ contains $\sqrt2$.  Its full constant field is exactly $\mathbb F_{25}$.
\end{lemma}

\begin{proof}
Let $t_1,t_2,t_3$ be the roots of $P(t)-b$, and choose $\eta_i$ in the Galois closure with
\[
 \eta_i^2=\frac{t_i-1}{t_i}.
\]
The monic cubic is
\[
 t^3-2t^2+b,
\]
whose discriminant in characteristic five is
\[
 \mathfrak D=2b(1-b).
\]
Vieta's formulas give
\[
 \prod_{i=1}^3\frac{t_i-1}{t_i}
 =-\frac{1-b}{b}
 =2\left(\frac{\sqrt{\mathfrak D}}b\right)^2.
\]
Consequently,
\[
 \sqrt2=\frac{b\eta_1\eta_2\eta_3}{\sqrt{\mathfrak D}}
\]
belongs to the Galois closure.  Since $2$ is nonsquare in $\mathbb F_5$, the constant field is nontrivial.

The arithmetic monodromy group contains $G$ normally and acts on the same six sheets, hence is contained in $N_{S_6}(G)$.  By Lemma~\ref{lem:normalizer}, its constant-field quotient has order at most two.  The element $\sqrt2$ forces degree at least two.  Therefore the full constant field is $\mathbb F_{25}$.
\end{proof}

\begin{proposition}[Arithmetic group and Frobenius parity]\label{prop:arithmetic-group}
The arithmetic monodromy group is
\[
 A=C_2\wr S_3.
\]
The $5^m$-power Frobenius lies in $G$ for even $m$ and in $\tau G$ for odd $m$.
\end{proposition}

\begin{proof}
Lemma~\ref{lem:constant-field} shows that $A/G$ is nontrivial, while Lemma~\ref{lem:normalizer} bounds it by $C_2$.  Thus $A=N_{S_6}(G)$.  If $\alpha^2=2$, then
\[
 \alpha^5=4\alpha=-\alpha.
\]
Hence the $5^m$-power map acts trivially on $\mathbb F_{25}$ exactly when $m$ is even.
\end{proof}

\subsubsection{Fixed-point inventory}

A direct calculation with the displayed permutations gives Table~\ref{tab:cycle-inventory}.  The first five rows are the five conjugacy classes of $G$ in the six-point action; multiplying by $\tau$ gives the five $G$-conjugacy classes in the outer coset.

\begin{table}[t]
\caption{Cycle and fixed-point inventory in the two Frobenius cosets.}\label{tab:cycle-inventory}
\centering
\small
\begin{tabular}{@{}ccrr@{}}
\toprule
Set & Six-point cycle type & Number & Fixed points \\
\midrule
$G$ & $1^6$ & 1 & 6 \\
$G$ & $4\,1^2$ & 6 & 2 \\
$G$ & $2^2\,1^2$ & 3 & 2 \\
$G$ & $2^3$ & 6 & 0 \\
$G$ & $3^2$ & 8 & 0 \\
\midrule
$\tau G$ & $2\,1^4$ & 3 & 4 \\
$\tau G$ & $2^2\,1^2$ & 6 & 2 \\
$\tau G$ & $2^3$ & 1 & 0 \\
$\tau G$ & $4\,2$ & 6 & 0 \\
$\tau G$ & $6$ & 8 & 0 \\
\bottomrule
\end{tabular}
\end{table}

Thus the fixed-point distributions are
\begin{equation}\label{eq:fixed-distributions}
 \begin{aligned}
 G:&\quad(0,2,6)\longmapsto(14,9,1),\\
 \tau G:&\quad(0,2,4)\longmapsto(15,6,3).
 \end{aligned}
\end{equation}
By Proposition~\ref{prop:frobenius-support}, an unramified finite fiber has size in $\{0,2,4\}$ for odd $m$ and in $\{0,2,6\}$ for even $m$.

\subsubsection{Special finite fibers}

Let $\chi$ be the quadratic character of $\mathbb F_q$, extended by $\chi(0)=0$, and put
\[
 \varepsilon:=\chi(2)=(-1)^m.
\]
The equality follows because $2$ is nonsquare over $\mathbb F_5$ and becomes square precisely over even-degree extensions.

At $b=1$,
\[
 P(t)-1=-(t-1)(t-3)^2.
\]
The point $t=1$ has one lift under $\kappa$, while $t=3$ has two.  Thus the rational core has three affine points.  The excluded inputs $x=0,-1$ both give finite differential value $1$, so
\begin{equation*}
 \delta_{F_m}(1,1)=5.
\end{equation*}
At $b=0$, only the finite quotient point $t=2$ contributes to the affine core, and its lift discriminant is $3=-2$.  Since $\chi(-1)=1$ over $\mathbb F_{5^m}$,
\begin{equation}\label{eq:zero-fiber}
 \delta_{F_m}(1,0)=1+\chi(3)=1+\varepsilon.
\end{equation}
There are no other finite-function corrections.  Combining \eqref{eq:fixed-distributions}--\eqref{eq:zero-fiber} proves the support assertion in Theorem~\ref{thm:main-monodromy}:
\[
 \operatorname{supp}(\omega)\subseteq
 \begin{cases}
 \{0,2,4,5\},&m\text{ odd},\\
 \{0,2,5,6\},&m\text{ even},
 \end{cases}
\]
and the five-fiber is unique.

\subsubsection{The required double-configuration curve and exact moment}\label{sec:double}

Let
\[
 M_2=\sum_{b\in\mathbb F_q}
 \delta_{F_m}(1,b)\bigl(\delta_{F_m}(1,b)-1\bigr).
\]
We compute $M_2$ by counting ordered distinct pairs in equal fibers.

\paragraph{The two ordered-pair orbits.}
The group $G$ has two orbits on ordered pairs of distinct sheets.  By Proposition~\ref{prop:orbit-component}, these are the two geometric components of the ordered double-configuration cover.  One is induced by the nontrivial element of the centralizer of $G$ and geometrically gives the graph
\[
 y=-1-x.
\]
Removing the pole inputs $x=0,-1$ and the fixed point $x=2$ leaves
\begin{equation}\label{eq:graph-count}
 q-3
\end{equation}
rational ordered pairs.

The other orbit consists of pairs whose quotient coordinates are distinct.  Put
\[
 u=x(x+1),\qquad v=y(y+1).
\]
The equality
\[
 \frac{2u+4}{u^3}=\frac{2v+4}{v^3},
 \qquad u\ne v,
\]
can be parameterized by $z=v/u$.  Dividing by $z-1$ gives
\begin{equation}\label{eq:uv-parameterization}
 u=\frac{z^2+z+1}{2z(z+1)},
 \qquad
 v=\frac{z^2+z+1}{2(z+1)}.
\end{equation}
The admissible parameter set is
\begin{equation}\label{eq:Z-set}
 \mathcal Z:=
 \{z\in\mathbb F_q:
 z\ne0,-1,1,\ z^2+z+1\ne0\}.
\end{equation}
Since the discriminant of $z^2+z+1$ is $-3=2$ in characteristic five,
\begin{equation}\label{eq:Z-size}
 |\mathcal Z|=q-4-\varepsilon.
\end{equation}

The numbers of lifts of $u$ and $v$ to $x$ and $y$ are respectively
\[
 1+\chi(A(z)),\qquad 1+\chi(B(z)),
\]
where
\begin{equation}\label{eq:A-B}
 A(z):=\frac{(z-2)^2}{2z(z+1)},
 \qquad
 B(z):=\frac{2(z+2)^2}{z+1}.
\end{equation}
Thus the cross-orbit count is
\begin{equation}\label{eq:cross-count-start}
 \sum_{z\in\mathcal Z}
 \bigl(1+\chi(A(z))\bigr)
 \bigl(1+\chi(B(z))\bigr).
\end{equation}

\paragraph{The required character sums.}
The following sums evaluate the cross-orbit contribution directly.

\begin{lemma}[A basic quadratic sum]\label{lem:basic-quadratic-sum}
For every $c\in\mathbb F_q^\times$,
\[
 \sum_{z\in\mathbb F_q}\chi(z^2-c)=-1.
\]
In particular,
\[
 \sum_{z\in\mathbb F_q}\chi(z(z+1))=-1.
\]
\end{lemma}

\begin{proof}
The number of pairs $(z,w)\in\mathbb F_q^2$ satisfying $w^2=z^2-c$ equals
\[
 q+\sum_z\chi(z^2-c).
\]
On the other hand, $(z-w)(z+w)=c$.  Since the characteristic is odd, each choice of $z-w\in\mathbb F_q^\times$ determines a unique pair $(z,w)$, so the number is $q-1$.  The first identity follows.  Completing the square gives the second.
\end{proof}

\begin{lemma}[Restricted character sums]\label{lem:restricted-sums}
For the set $\mathcal Z$ in \eqref{eq:Z-set},
\begin{align*}
 \sum_{z\in\mathcal Z}\chi(z(z+1))&=-2-2\varepsilon,\\
 \sum_{z\in\mathcal Z}\chi(z+1)&=-2-2\varepsilon,\\
 \sum_{z\in\mathcal Z}\chi(z)&=-3-\varepsilon.
\end{align*}
\end{lemma}

\begin{proof}
The first unrestricted sum is $-1$ by Lemma~\ref{lem:basic-quadratic-sum}; the other two unrestricted sums are zero.  The removed values $0,-1,1$ contribute, where nonzero, $1$ and $\varepsilon=\chi(2)$.

If $\rho^2+\rho+1=0$, then $\rho(\rho+1)=-1$ and $\rho+1=-\rho^2$, so both corresponding character values equal one.  Such roots exist exactly when $\varepsilon=1$, and then there are two of them.  They are primitive cube roots of unity; since $q-1$ is then divisible by six, they are squares and $\chi(\rho)=1$.  Subtracting the removed contributions gives the three formulas.
\end{proof}

\begin{lemma}[The three lift sums]\label{lem:lift-sums}
The functions in \eqref{eq:A-B} satisfy
\begin{align*}
 \sum_{z\in\mathcal Z}\chi(A(z))
 =\sum_{z\in\mathcal Z}\chi(B(z))&=-2-3\varepsilon,\\
 \sum_{z\in\mathcal Z}\chi(A(z)B(z))&=-3-3\varepsilon.
\end{align*}
\end{lemma}

\begin{proof}
Away from $z=2$,
\[
 \chi(A(z))=\varepsilon\chi(z(z+1)).
\]
At $z=2$, the left-hand side is zero while the right-hand side is $\varepsilon$.  Lemma~\ref{lem:restricted-sums} therefore gives
\[
 \sum_{\mathcal Z}\chi(A)
 =\varepsilon(-2-2\varepsilon)-\varepsilon
 =-2-3\varepsilon.
\]
Similarly, away from $z=-2$,
\[
 \chi(B(z))=\varepsilon\chi(z+1),
\]
and the same correction gives the second equality.

Finally,
\[
 A(z)B(z)
 =\frac{(z^2-4)^2}{z(z+1)^2}.
\]
Thus $\chi(A(z)B(z))=\chi(z)$ away from $z=\pm2$.  Each of the two exceptional right-hand values is $\varepsilon$, and Lemma~\ref{lem:restricted-sums} yields
\[
 \sum_{\mathcal Z}\chi(AB)
 =(-3-\varepsilon)-2\varepsilon
 =-3-3\varepsilon.
\]
\end{proof}

Combining \eqref{eq:Z-size}, \eqref{eq:cross-count-start}, and Lemma~\ref{lem:lift-sums}, the cross component contributes
\begin{equation}\label{eq:cross-count}
 q-11-10\varepsilon.
\end{equation}

\paragraph{Explicit conic normalization.}
The character calculation has a direct geometric form.  Let
\[
 \xi=2x+1,\qquad \eta=2y+1.
\]
Then $\xi^2=A(z)$ and $\eta^2=B(z)$.  On the dense open set where the indicated factors are nonzero, put
\[
 R=\frac{\xi}{z-2},
 \qquad
 T=\frac{\eta}{z+2}.
\]
Equations \eqref{eq:A-B} become
\[
 R^2=\frac1{2z(z+1)},
 \qquad
 T^2=\frac2{z+1}.
\]
Now define
\[
 W=\frac{T}{2R},
 \qquad
 S=\frac2T.
\]
Then
\[
 W^2=z,
 \qquad
 S^2=2(W^2+1).
\]
Conversely,
\[
 z=W^2,
 \qquad
 T=\frac2S,
 \qquad
 R=\frac1{SW}
\]
on a dense open subset.  Hence the smooth projective normalization of the cross component is the conic
\[
 \mathcal Q:\quad S^2=2(W^2+1).
\]

The finite deletions in the parameterization are summarized in Table~\ref{tab:boundary}.

\begin{table*}[t]
\caption{Boundary and exceptional parameters in the cross component.}\label{tab:boundary}
\centering
\small
\begin{tabularx}{\textwidth}{@{}p{0.15\textwidth}X@{}}
\toprule
Parameter & Geometric meaning \\
\midrule
$z=0$ & $v=0$, so the $y$-coordinate belongs to the excluded pole locus. \\
$z=-1$ & The rational parameterization \eqref{eq:uv-parameterization} has a pole. \\
$z=1$ & $u=v$; the pair belongs to the diagonal or to the involution component rather than the cross orbit. \\
$z^2+z+1=0$ & Both quotient coordinates collapse to $0$, again giving excluded pole inputs. \\
$z=2$ & $A(z)=0$; there is one ramified $x$-lift, and the point is retained. \\
$z=-2$ & $B(z)=0$; there is one ramified $y$-lift, and the point is retained. \\
\bottomrule
\end{tabularx}
\end{table*}

Thus the term involving $\varepsilon$ in the exact point count is a boundary splitting contribution; the normalized component itself has genus zero and no nontrivial $H^1$ trace.

\paragraph{The exact second moment.}
\begin{proposition}[Second factorial moment]\label{prop:second-moment}
With $\varepsilon=(-1)^m$,
\[
 M_2=2q-10\varepsilon.
\]
\end{proposition}

\begin{proof}
The rational-core count is the sum of the graph contribution \eqref{eq:graph-count} and the cross contribution \eqref{eq:cross-count}:
\[
 (q-3)+(q-11-10\varepsilon)
 =2q-14-10\varepsilon.
\]
At $b=1$, the two exceptional inputs increase the fiber size from three to five.  The falling-factorial moment therefore increases by
\[
 (5)_2-(3)_2=20-6=14.
\]
Hence
\[
 M_2=2q-10\varepsilon.
\]
\end{proof}

This proves the moment assertion in Theorem~\ref{thm:main-monodromy}.

\paragraph{Reconstruction of the complete spectrum.}
\begin{proof}[Proof of Theorem~\ref{thm:main-spectrum}]
Suppose first that $m$ is odd.  Then $\varepsilon=-1$, the support is contained in $\{0,2,4,5\}$, and $\omega_5=1$.  The first three falling-factorial moments give
\begin{align*}
 \omega_0+\omega_2+\omega_4&=q-1,\\
 2\omega_2+4\omega_4&=q-5,\\
 2\omega_2+12\omega_4&=M_2-(5)_2=2q-10.
\end{align*}
Solving yields
\[
 \omega_4=\frac{q-5}{8},
 \qquad
 \omega_2=\frac{q-5}{4},
 \qquad
 \omega_0=\frac{5q+7}{8}.
\]

Now suppose that $m$ is even.  Then $\varepsilon=1$, the support is contained in $\{0,2,5,6\}$, and $\omega_5=1$.  Thus
\begin{align*}
 \omega_0+\omega_2+\omega_6&=q-1,\\
 2\omega_2+6\omega_6&=q-5,\\
 2\omega_2+30\omega_6&=M_2-(5)_2=2q-30.
\end{align*}
Hence
\[
 \omega_6=\frac{q-25}{24},
 \qquad
 \omega_2=\frac{3q+5}{8},
 \qquad
 \omega_0=\frac{7(q-1)}{12}.
\]
For odd $m$, $q\equiv5\pmod8$; for even $m$, $q\equiv1\pmod{24}$.  These congruences also show directly that every coefficient is integral.
\end{proof}

The leading proportions agree with the coset inventories.  For odd $m$, the proportions of $0$-, $2$-, and $4$-fibers tend to
\[
 \frac{15}{24},\qquad\frac6{24},\qquad\frac3{24},
\]
while for even $m$, the proportions of $0$-, $2$-, and $6$-fibers tend to
\[
 \frac{14}{24},\qquad\frac9{24},\qquad\frac1{24}.
\]
The exact constant terms come from branch fibers, exceptional inputs, and the boundary of the conic component.

\paragraph{Full DDT and cryptographic parameters.}
\begin{proof}[Proof of Corollary~\ref{cor:ddt}]
For $a\ne0$, substitute $x=ay$.  Homogeneity gives
\[
 F_m(x+a)-F_m(x)
 =a^{q-4}\bigl(F_m(y+1)-F_m(y)\bigr),
\]
which proves the row-rescaling identity and locates the unique five-entry.

The uniformity follows from Theorem~\ref{thm:main-spectrum}; for $m=2$, the coefficient $\omega_6$ is zero.  The maximum differential probability, denoted by $\operatorname{MDP}(F_m)$, is the largest DDT entry divided by $q$.

The power map is a permutation exactly when
\[
 \gcd(q-4,q-1)=\gcd(3,q-1)=1.
\]
Since $5^m\equiv(-1)^m\pmod3$, this occurs precisely for odd $m$.

Finally,
\[
 \sum_b\delta_{F_m}(1,b)^2
 =M_2+\sum_b\delta_{F_m}(1,b)
 =2q-10(-1)^m+q.
\]
Equation \eqref{eq:global-energy} gives the stated formula for $E_4(F_m)$.
\end{proof}

The zero DDT row is
\[
 \delta_{F_m}(0,0)=q,
 \qquad
 \delta_{F_m}(0,b)=0\quad(b\ne0).
\]
Together with Corollary~\ref{cor:ddt}, this describes the entire table.

\paragraph{Endpoint fields.}
At $q=5$, the theorem gives $(\omega_0,\omega_5)=(4,1)$,
consistent with $F_1(x)=x$.  At $q=25$, it gives
$(\omega_0,\omega_2,\omega_5)=(14,10,1)$; the first six-fiber
occurs at $q=5^4$.  These endpoint checks are not used in the proof.

\subsection{Characteristics greater than five}\label{subsec:generic-inverse-cubic}

Let $q=p^n$ with $p>5$ and $F(x):=x^{q-4}$, and let $\chi$ be
the quadratic character of $\mathbb F_q$, extended by $\chi(0)=0$.
Define
\begin{align*}
 \rho_q&:=1+\sum_{z^2-3z+6=0}\bigl(1+\chi(z)\bigr),\\
 \sigma_q&:=2+\sum_{u^3+3u+1=0}\bigl(1+\chi(1+4u)\bigr),
\end{align*}
so $\rho_q\in\{1,3,5\}$ and $\sigma_q\in\{2,4,6,8\}$.  Let $E_{\rm ic},G_{\rm ic},H_{\rm ic}$ be the smooth projective models
\begin{align*}
 E_{\rm ic}:&\quad Y^2=-3Z(Z+1)(Z^2+Z+4),\\
 G_{\rm ic}:&\quad Y^2=Z(Z^2+Z+4)(4Z^2+Z+1),\\
 H_{\rm ic}:&\quad Y^2=3(Z^2+Z+4)(4Z^2+Z+1)(4Z^2+7Z+4),
\end{align*}
with traces $a_E:=q+1-\#E_{\rm ic}(\mathbb F_q)$, $a_G:=q+1-\#G_{\rm ic}(\mathbb F_q)$, and $a_H:=q+1-\#H_{\rm ic}(\mathbb F_q)$.  Put
\begin{align*}
 M_2^{\rm ic}:={}&2q-17-3\chi(-1)-6\chi(-3)\\
 &{}-2a_E-a_G+4\sigma_q,\\
 M_3^{\rm ic}:={}&4q-13-18\chi(-1)-29\chi(-3)\\
 &{}-6\chi(5)-6\chi(-5)-3\chi(3)\\
 &{}-9a_E-6a_G-a_H-3\rho_q\\
 &{}+6\sigma_q^2-24\sigma_q,
\end{align*}
and
\begin{align*}
 T_q&:=q-\rho_q-\sigma_q,\\
 A_q&:=M_2^{\rm ic}-(\rho_q)_2-(\sigma_q)_2,\\
 B_q&:=M_3^{\rm ic}-(\rho_q)_3-(\sigma_q)_3.
\end{align*}
Define
\begin{align*}
 \nu_0&:=q-2-\frac{11T_q}{16}+\frac{3A_q}{16}-\frac{B_q}{48},\\
 \nu_2&:=\frac{15T_q-7A_q+B_q}{16},\\
 \nu_4&:=\frac{-5T_q+5A_q-B_q}{16},\\
 \nu_6&:=\frac{3T_q-3A_q+B_q}{48}.
\end{align*}

\begin{theorem}[Generic inverse-cubic spectrum]\label{thm:generic-inverse-cubic}
For $q=p^n$ with $p>5$, the complete differential spectrum of $F(x)=x^{q-4}$ is
\begin{align*}
 \omega_0&=\nu_0,&
 \omega_2&=\nu_2+\mathbf{1}_{\{\sigma_q=2\}},\\
 \omega_4&=\nu_4+\mathbf{1}_{\{\sigma_q=4\}},&
 \omega_6&=\nu_6+\mathbf{1}_{\{\sigma_q=6\}},\\
 \omega_8&=\mathbf{1}_{\{\sigma_q=8\}},&
 \omega_1&=\mathbf{1}_{\{\rho_q=1\}},\\
 \omega_3&=\mathbf{1}_{\{\rho_q=3\}},&
 \omega_5&=\mathbf{1}_{\{\rho_q=5\}},
\end{align*}
and every other coefficient is zero.
\end{theorem}

For $x\ne0,-1$, set
\[
 u=x(x+1),\qquad t=u^{-1}.
\]
Then
\begin{equation}\label{icg:eq:cubic-core}
 F(x+1)-F(x)=-\frac{3u+1}{u^3}=-t^2(t+3),
\end{equation}
and with $\xi=2x+1$,
\begin{equation*}
 \xi^2=1+4u=\frac{t+4}{t}.
\end{equation*}
Thus the same degree-two involution quotient is followed by the fixed cubic $P(t)=-t^2(t+3)$.

\begin{lemma}[Special targets and support]\label{icg:lem:support}
One has
\[
 \delta(16)=\rho_q,\qquad \delta(1)=\sigma_q.
\]
For every $b\notin\{1,16\}$,
\[
 \delta(b)\in\{0,2,4,6\}.
\]
Consequently odd fiber sizes can occur only at $b=16$.
\end{lemma}

\begin{proof}
For a root $t\ne0$ of $P(t)=b$, the number of lifts to $x$ is
\[
 1+\chi\!\left(\frac{t+4}{t}\right).
\]
It is $0$ or $2$, except at $t=-4$, where it is $1$; since $P(-4)=16$, all other ordinary fibers are even and have size at most six.  At $b=16$,
\[
 t^3+3t^2+16=(t+4)(t^2-t+4).
\]
The root $t=-4$ contributes one point.  For a root of $t^2-t+4=0$, set $z=(t+4)/t$; then $z^2-3z+6=0$, giving exactly the defining formula for $\rho_q$.  Finally, the finite polynomial differs from the rational core only at $x=0,-1$, and both inputs have derivative value $1$.  Setting $b=1$ in \eqref{icg:eq:cubic-core} gives $u^3+3u+1=0$, with lift number $1+\chi(1+4u)$, which gives $\sigma_q$.
\end{proof}

The four remaining regular multiplicities are determined by moments through order three, by Proposition~\ref{prop:relation-aware-stopping}.

\subsubsection{Second factorial moment}

Let $U=\mathbb F_q\setminus\{0,-1\}$ and let $\delta^\circ$ count solutions in $U$.  For a pair of ordinary solutions $x,y$ with equal derivative value, put $u=x(x+1)$ and $v=y(y+1)$.  Removing the component $u=v$, the equality of derivatives is
\begin{equation*}
 3u^2v+u^2+3uv^2+uv+v^2=0.
\end{equation*}
With $z=v/u$,
\begin{equation*}
 u=-\frac{z^2+z+1}{3z(z+1)},\qquad
 v=-\frac{z^2+z+1}{3(z+1)}.
\end{equation*}
Put
\begin{align*}
 A(z)&:=-\frac{z^2+z+4}{3z(z+1)},\\
 B(z)&:=-\frac{4z^2+z+1}{3(z+1)},\qquad
 Q:=z^2+z+1.
\end{align*}
The two lift numbers are $1+\chi(A(z))$ and $1+\chi(B(z))$.  For
\[
 S_2:=\mathbb F_q\setminus\bigl(\{0,1,-1\}\cup Z(Q)\bigr),
\]
define
\[
 C^\circ:=\sum_{z\in S_2}(1+\chi(A(z)))(1+\chi(B(z))).
\]

\begin{lemma}[Cross-pair count]\label{icg:lem:cross-pair}
One has
\begin{equation}\label{icg:eq:C}
 C^\circ=q-8-3\chi(-1)-6\chi(-3)-2a_E-a_G.
\end{equation}
\end{lemma}

\begin{proof}
Since $\#Z(Q)=1+\chi(-3)$, one has $|S_2|=q-4-\chi(-3)$.  The substitution $Y=3z(z+1)w$ on $w^2=A(z)$ and $w^2=A(z)B(z)$ gives the models $E_{\rm ic}$ and $G_{\rm ic}$, respectively.  Their projective point counts, with deleted parameters removed, give
\begin{align*}
 \sum_{S_2}\chi(A)&=-a_E-\chi(-1)-1-2\chi(-3),\\
 \sum_{S_2}\chi(AB)&=-a_G-\chi(-1)-2-\chi(-3).
\end{align*}
Indeed, $A(1)=-1$ and $A=1$ at the roots of $Q$; the defining polynomial $f_G$ of $G_{\rm ic}$ satisfies $f_G(-1)=-16$, $f_G(1)=36$, and is a nonzero square at the roots of $Q$.  The change $z\mapsto z^{-1}$ interchanges $A$ and $B$; expanding the product gives \eqref{icg:eq:C}.
\end{proof}

\begin{proposition}[Second moment for $p>5$]\label{icg:prop:M2}
The second factorial moment is $M_2^{\rm ic}$ defined above.
\end{proposition}

\begin{proof}
The same-$u$ ordered pairs are exactly $(x,-1-x)$ with $x\in U\setminus\{-1/2\}$, contributing $q-3$.  Hence $M_2^\circ=q-3+C^\circ$.  Restoring $x=0,-1$ raises the fiber at $b=1$ from $\sigma_q-2$ to $\sigma_q$, adding
\[
 (\sigma_q)_2-(\sigma_q-2)_2=4\sigma_q-6.
\]
Substitution of \eqref{icg:eq:C} yields the displayed formula for $M_2^{\rm ic}$.
\end{proof}

\subsubsection{Third factorial moment}

For an ordered triple of distinct ordinary solutions, the quotient coordinates have occupancy type $2+1$ or $1+1+1$.  To count the first type, remove from $C^\circ$ the $\rho_q-1$ cross-pairs whose first input is the fixed point $x=-1/2$; every remaining first input has a distinct involution partner.  The three positions of the singleton then give $3(C^\circ-\rho_q+1)$ ordered triples.  For three distinct quotient coordinates, if $t_1,t_2,t_3$ are the ordered roots of $t^3+3t^2+b=0$ and $z=t_1/t_2$, then with $Q=z^2+z+1$,
\begin{equation*}
 t_1=-\frac{3z(z+1)}Q,\qquad
 t_2=-\frac{3(z+1)}Q,\qquad
 t_3=\frac{3z}Q.
\end{equation*}
The first two lift discriminants are the previously defined $A(z),B(z)$; the third is
\[
 C(z):=\frac{4z^2+7z+4}{3z}.
\]
Distinctness is equivalent to
\[
 z\in S_3:=\mathbb F_q\setminus\bigl(\{0,-1,1,-2,-1/2\}\cup Z(Q)\bigr).
\]
Put
\[
 L^\circ:=\sum_{z\in S_3}(1+\chi(A(z)))(1+\chi(B(z)))(1+\chi(C(z))).
\]

\begin{lemma}[Three-block count]\label{icg:lem:L}
One has
\begin{align}
 L^\circ={}&q-16-9\chi(-1)-11\chi(-3)\notag\\
 &{}-6\chi(5)-6\chi(-5)-3\chi(3)\notag\\
 &{}-3a_E-3a_G-a_H.\label{icg:eq:L}
\end{align}
\end{lemma}

\begin{proof}
The transformations $z\mapsto z^{-1}$ and $z\mapsto-z-1$ preserve $S_3$ and exchange $A,B$ and $B,C$, respectively.  Thus the three single-character sums coincide, as do the three pair-character sums, and
\[
 L^\circ=|S_3|+3S_A+3S_{AB}+S_{ABC}.
\]
The first two sums reuse the cross-pair calculation: deleting the additional parameters $z=-2,-1/2$ removes $(A,AB)=(-1,-5)$ and $(5,-5)$, respectively.  Consequently,
\begin{align*}
 S_A&=\sum_{S_2}\chi(A)-\chi(-1)-\chi(5),\\
 S_{AB}&=\sum_{S_2}\chi(AB)-2\chi(-5).
\end{align*}
For $w^2=A(z)B(z)C(z)$, the substitution $Y=9z(z+1)w$ gives $H_{\rm ic}$.  Write $f_H$ for its defining polynomial; on $S_3$ one has $\chi(ABC)=\chi(f_H)$.  At $z=0,-1,1,-2,-1/2$ the values of $f_H$ are $48,48,1620,1620,405/16$, respectively, and $f_H=81$ at each root of $Q$.  Subtracting the deleted affine terms and accounting for the $1+\chi(3)$ points at infinity gives
\[
 S_{ABC}=-a_H-3\chi(3)-3\chi(5)-1-\chi(-3).
\]
Since $|S_3|=q-6-\chi(-3)$, substitution yields \eqref{icg:eq:L}.
\end{proof}

\begin{proposition}[Third moment for $p>5$]\label{icg:prop:M3}
The third factorial moment is $M_3^{\rm ic}$ defined above.
\end{proposition}

\begin{proof}
The ordinary contribution is
\[
 M_3^\circ=3(C^\circ-\rho_q+1)+L^\circ.
\]
Restoring $x=0,-1$ in the fiber $b=1$ adds
\[
 (\sigma_q)_3-(\sigma_q-2)_3=6\sigma_q^2-24\sigma_q+24.
\]
Insert \eqref{icg:eq:C} and \eqref{icg:eq:L} and simplify.
\end{proof}

\begin{proof}[Proof of Theorem~\ref{thm:generic-inverse-cubic}]
Remove the two targets $1$ and $16$, and let $\nu_j$ count the remaining targets of fiber size $j\in\{0,2,4,6\}$.  By Lemma~\ref{icg:lem:support} and Propositions~\ref{icg:prop:M2}--\ref{icg:prop:M3},
\begin{align*}
 \nu_0+\nu_2+\nu_4+\nu_6&=q-2,\\
 2\nu_2+4\nu_4+6\nu_6&=T_q,\\
 2\nu_2+12\nu_4+30\nu_6&=A_q,\\
 24\nu_4+120\nu_6&=B_q.
\end{align*}
The coefficient matrix is nonsingular, and solving gives precisely the four values $\nu_0,\nu_2,\nu_4,\nu_6$ defined above.  Finally reinsert the unique target $b=1$, of size $\sigma_q$, and the unique target $b=16$, of size $\rho_q$.
\end{proof}

\begin{remark}[Why characteristic five is exceptional]\label{icg:rem:p5}
The quadratic factors
\[
 Z^2+Z+4,\qquad 4Z^2+Z+1,\qquad 4Z^2+7Z+4
\]
all have discriminant $-15$, and their pairwise resultants are $216$.  Hence the three controlling curves are smooth for every $p>5$, while characteristic five forces a collision.  Equivalently, the generic branch values $-4$ and $16$ coincide when $p=5$.  Thus the conic calculation above is a genuine degeneration, not a specialization of the generic depth-three formulas.
\end{remark}

\section{Two Ternary APN Families with a Common Exact Spectrum}\label{sec:dmmp}

Let
\[
 q=3^n,\qquad n\geqslant1\ \text{odd},\qquad
 m=\frac{n+1}{2},\qquad \sigma:=3^m,
\]
and let $\chi:\mathbb{F}_{q}\to\{0,\pm1\}$ be the quadratic character, extended by $\chi(0)=0$.  Since $n$ is odd,
\[
 \chi(-1)=-1.
\]
The four branches considered here are the two cases of Theorems~2.1 and~2.2 of Dobbertin--Mills--M\"uller--Pott--Willems~\cite{DMMP2003}.  We use the following branch labels:
\begin{equation}
\begin{array}{c@{\qquad}c@{\qquad}c}
\toprule
\text{branch}& n\bmod4 & d\\
\midrule
A_3&3&\dfrac{\sigma-1}{2}\\[1.1ex]
A_1&1&\dfrac{\sigma-1}{2}+\dfrac{q-1}{2}\\[1.1ex]
B_3&3&\dfrac{3q-1}{8}\\[1.1ex]
B_1&1&\dfrac{3q-1}{8}+\dfrac{q-1}{2}.\\
\bottomrule
\end{array}
\label{dmmp6:eq:branches}
\end{equation}
For $F(x)=x^d$, put
\[
 N_d(b):=\#\{x\in\mathbb{F}_{q}:(x+1)^d-x^d=b\},
\]
\[
 \omega_j=\#\{b\in\mathbb{F}_{q}:N_d(b)=j\}.
\]
DMMP proved $N_d(b)\leqslant2$ for all four branches
\cite[Theorems~2.1--2.2]{DMMP2003}.  Leducq subsequently proved
$\Delta=2$ for every odd $n\geqslant5$ \cite{Leducq2012}.
The result below evaluates the frequencies and proves that the
$A$- and $B$-branches have the same spectrum; its exact-uniformity
consequence reproves this known result.  The upper bound gives
\begin{equation}
 \omega_0=\omega_2,\qquad \omega_1=q-2\omega_2,
 \label{dmmp6:eq:few-support}
\end{equation}
because
\[
 \omega_0+\omega_1+\omega_2=q,
 \qquad
 \omega_1+2\omega_2=q.
\]
Hence only $\omega_2$ has to be computed.

Let
\begin{equation}
 E/\mathbb{F}_{3}:\qquad Y^2=(X+1)(X^2+1),
 \qquad
 t_n:=q+1-\#E(\mathbb{F}_{q}).
 \label{dmmp6:eq:E}
\end{equation}

\begin{theorem}[Unified DMMP spectrum]
\label{dmmp6:thm:main}
For every odd $n$ and every exponent $d$ in \eqref{dmmp6:eq:branches},
\begin{equation}
 \begin{aligned}
 \omega_0=\omega_2&=\frac{q-7-2t_n}{4},\\
 \omega_1&=\frac{q+7+2t_n}{2},\qquad
 \omega_j=0\ (j\geqslant3).
 \end{aligned}
 \label{dmmp6:eq:main-spectrum}
\end{equation}
Moreover
\begin{equation}
 t_0=2,\qquad t_1=-2,\qquad
 t_n=-2t_{n-1}-3t_{n-2}\quad(n\geqslant2).
 \label{dmmp6:eq:trace-recurrence}
\end{equation}
Consequently the admissible branches are planar at the endpoint fields $n=1,3$ and have differential uniformity exactly two for every odd $n\geqslant5$.  For a fixed odd $n$, the corresponding $A$- and $B$-branches have identical complete differential spectra.
\end{theorem}

\begin{proposition}[Nontriviality of the DMMP differential isospectrality]
\label{dmmp8:prop:nontrivial-isospectrality}
Let $n\geqslant5$ be odd, let $N=3^n-1$, and let $d_A,d_B$ be the corresponding $A$- and $B$-branch exponents in \eqref{dmmp6:eq:branches}.  Then
\begin{equation}
 \gcd(d_A,N)=\gcd(d_B,N)=2
 \label{dmmp8:eq:gcd-exponents}
\end{equation}
and
\begin{equation}
 d_B\not\equiv \pm3^j d_A\pmod N
 \qquad(0\leqslant j<n).
 \label{dmmp8:eq:not-signed-frobenius}
\end{equation}
In particular, the two exponents lie in distinct Frobenius exponent orbits, so their common spectrum is not forced by Frobenius equivalence.  The signed exclusion in \eqref{dmmp8:eq:not-signed-frobenius} is purely arithmetic: negating a nonpermutation exponent need not preserve the differential spectrum.  The inverse-permutation relation does not apply, since neither exponent is invertible modulo $N$.  At the planar endpoints the situation is different: $d_{A_1}=d_{B_1}=2$ for $n=1$, whereas
\[
 d_{B_3}\equiv3^2d_{A_3}\pmod{26}
\]
for $n=3$.
\end{proposition}

\begin{proof}
Write $m=(n+1)/2$.  For the $A$-branch,
\[
 2d_A\equiv3^m-1\pmod N.
\]
Since $\gcd(m,n)=1$,
\[
 \gcd(3^m-1,N)=3^{\gcd(m,n)}-1=2.
\]
Thus $\gcd(d_A,N)\mid2$; because both $d_A$ and $N$ are even, the first equality in \eqref{dmmp8:eq:gcd-exponents} follows.  For the $B$-branch,
\[
 8d_B\equiv2\pmod N,
\]
so $\gcd(d_B,N)\mid2$, and again $d_B$ is even.  This proves \eqref{dmmp8:eq:gcd-exponents}.

For $0<d<N$, write
\[
 d=\sum_{i=0}^{n-1}a_i3^i,\qquad a_i\in\{0,1,2\},
\]
and define the cyclic transition number
\[
 \tau(d)=\#\{i\in\mathbb Z/n\mathbb Z:a_i\ne a_{i+1}\}.
\]
Multiplication by $3$ modulo $N$ cyclically shifts the length-$n$ ternary word, so $\tau(3^jd\bmod N)=\tau(d)$.  Moreover
\[
 N-d=\sum_{i=0}^{n-1}(2-a_i)3^i,
\]
hence digitwise complementation also preserves $\tau$.

Expanding the geometric sums in \eqref{dmmp6:eq:branches} gives the following ternary words, written from least to most significant digit; superscripts on digit strings indicate concatenated repetitions:
\[
\begin{array}{c|cc|cc}
 n & d_A & d_B & \tau(d_A)&\tau(d_B)\\ \hline
 4r+3 & 1^{\,2r+2}0^{\,2r+1} &(10)^{2r+1}1&2&n-1\\
 4r+1 & 2^{\,2r+1}1^{\,2r} &(21)^{2r}2&2&n-1
\end{array}
\qquad(r\geqslant1).
\]
For $n=4r+3$ the exponents are $d_A=\sum_{i=0}^{2r+1}3^i$ and $d_B=\sum_{i=0}^{2r+1}3^{2i}$; for $n=4r+1$, the additional term $(q-1)/2=\sum_{i=0}^{4r}3^i$ increases each base-three digit by one without carrying.  Since $2\ne n-1$ for odd $n\geqslant5$, cyclic shifts and digitwise complementation cannot relate $d_A$ and $d_B$.  This proves \eqref{dmmp8:eq:not-signed-frobenius}.

Finally, if $e\equiv3^jd\pmod N$, then
\[
 (x+1)^e-x^e
 =
 \bigl((x+1)^d-x^d\bigr)^{3^j},
\]
so such a relation automatically preserves the differential spectrum.  The proposition shows that this standard explanation is absent for the nonplanar DMMP pairs.
\end{proof}

\begin{remark}
The two maps are nevertheless related by a power permutation of the output.  Write $N=2h$, $d_A=2a$, and $d_B=2b$.  By \eqref{dmmp8:eq:gcd-exponents}, $a$ and $b$ are units modulo the odd integer $h$.  Choose $c\equiv ba^{-1}\pmod h$ and, among the two lifts modulo $N$, choose the odd one.  Then $\gcd(c,N)=1$ and $d_B\equiv cd_A\pmod N$, so $x^{d_B}$ is obtained from $x^{d_A}$ by a multiplicative power permutation on the output.  For $n\geqslant5$, however, Proposition~\ref{dmmp8:prop:nontrivial-isospectrality} shows that $c$ cannot be a Frobenius exponent $3^j$.  Nonlinear power postcomposition need not preserve additive derivatives and therefore does not explain the differential isospectrality.
\end{remark}

\begin{corollary}[DMMP DDT collision statistics]
\label{dmmp8:cor:ddt-energy}
For every nonzero input difference $a\in\mathbb{F}_{q}^\times$, the corresponding DDT row has the spectrum \eqref{dmmp6:eq:main-spectrum}.  In particular,
\begin{align*}
 \sum_{b\in\mathbb{F}_{q}}\delta_F(a,b)\bigl(\delta_F(a,b)-1\bigr)
 &=\frac{q-7-2t_n}{2},\\
 \sum_{b\in\mathbb{F}_{q}}\delta_F(a,b)^2
 &=\frac{3q-7-2t_n}{2}.
\end{align*}
For odd $n\geqslant5$, the maximum differential probability is $2/q$; at $n=1,3$ it is $1/q$.
\end{corollary}

\begin{proof}
By homogeneity, $\delta_F(a,b)=\delta_F(1,b/a^d)$ for $a\ne0$.  As the support is contained in $\{0,1,2\}$ and $\sum_b\delta_F(a,b)=q$,
\[
 \sum_b\delta_F(a,b)(\delta_F(a,b)-1)=2\omega_2,
 \qquad \sum_b\delta_F(a,b)^2=q+2\omega_2.
\]
Substituting \eqref{dmmp6:eq:main-spectrum} proves both identities.
\end{proof}

For all four branches the known support bound $\{0,1,2\}$ reduces spectrum recovery to a single double-collision count.  The same Kummer coordinate applies to both exponents; the collision count is then evaluated by two different arguments, one using a fixed quartic and the other a character sum involving the Frobenius automorphism $z\mapsto z^\sigma$.

\subsection{A common Kummer parameter}

Define
\begin{equation*}
 \lambda:=
 \begin{cases}
 \sigma=3^{(n+1)/2},&\text{for }A_3,A_1,\\
 3,&\text{for }B_3,B_1,
 \end{cases}
 \qquad r:=\lambda+1.
\end{equation*}

\begin{lemma}[Exponent congruence]
\label{dmmp6:lem:master-congruence}
For every one of the four branches, $d$ is even and
\begin{equation*}
 rd\equiv\frac{q+1}{2}\pmod{q-1}.
\end{equation*}
Moreover
\begin{equation}
 \gcd(r,q-1)=2.
 \label{dmmp6:eq:gcd-r}
\end{equation}
\end{lemma}

\begin{proof}
For the $B$-branches, $r=4$ and $4d\equiv(q+1)/2\pmod{q-1}$.  Since $q\equiv3\pmod4$, one has $\gcd(4,q-1)=2$.  The parity of $d$ follows from $q\equiv11\pmod{16}$ when $n\equiv3\pmod4$ and $q\equiv3\pmod{16}$ when $n\equiv1\pmod4$.

For the $A$-branches,
\[
 (\sigma+1)\frac{\sigma-1}{2}
 =\frac{\sigma^2-1}{2}
 =\frac{3q-1}{2}
 \equiv\frac{q+1}{2}\pmod{q-1},
\]
and adding $(q-1)/2$ to $d$ changes the left-hand side by $((\sigma+1)/2)(q-1)$.  If an integer divides both $\sigma+1$ and $q-1=3^{2m-1}-1$, then $3^m\equiv-1$ and $3^{2m-1}\equiv1$ modulo that integer; comparison with $3^{2m}\equiv1$ forces $3\equiv1$, hence the common divisor divides $2$.  Thus \eqref{dmmp6:eq:gcd-r} holds.  Finally, if $n\equiv3\pmod4$, then $m$ is even and $(\sigma-1)/2$ is even; if $n\equiv1\pmod4$, then both $(\sigma-1)/2$ and $(q-1)/2$ are odd, so their sum is even.
\end{proof}

\begin{lemma}[Exact Kummer coordinate]
\label{dmmp6:lem:kummer-coordinate}
For $x\in\mathbb{F}_{q}$, put $u=x^d$.  Then $u$ is a square and
\begin{equation}
 u^r=\chi(x)x.
 \label{dmmp6:eq:kummer-coordinate}
\end{equation}
Conversely, for fixed $x$, $x^d$ is the unique square $u$ satisfying \eqref{dmmp6:eq:kummer-coordinate}.
\end{lemma}

\begin{proof}
Since $d$ is even, $u$ is a square.  By Lemma~\ref{dmmp6:lem:master-congruence},
\[
 u^r=x^{rd}=x^{(q+1)/2}=\chi(x)x.
\]
Conversely, if $u$ is square and satisfies \eqref{dmmp6:eq:kummer-coordinate}, then, for $x\ne0$,
\[
 x^d=(\chi(x)u^r)^d=u^{rd}=u^{(q+1)/2}=u,
\]
because $d$ is even and $u^{(q-1)/2}=1$.  The case $x=0$ is immediate.
\end{proof}

\subsection{Derivative blocks and double fibers}

Let
\[
 D(x):=(x+1)^d-x^d.
\]
For an interior solution $D(x)=b$, $x\ne0,-1$, set
\[
 \begin{aligned}
 U&:=x^d, & V&:=(x+1)^d=U+b,\\
 s&:=\chi(x), & t&:=\chi(x+1).
 \end{aligned}
\]
By Lemma~\ref{dmmp6:lem:kummer-coordinate},
\begin{equation}
 tV^r-sU^r=1.
 \label{dmmp6:eq:master-derivative}
\end{equation}
The two cases $t=s$ and $t=-s$ will be called the same-character and opposite-character blocks.

The following result of DMMP applies to $r=\sigma+1$ and $r=4$.

\begin{lemma}[DMMP square separation]
\label{dmmp6:lem:dmmp-separation}
Let $Q$ be the set of squares in $\mathbb{F}_{q}$, and put
\[
 M_a^{(r)}:=\{u^r+(u+a)^r:u,u+a\in Q\}.
\]
For $a\ne0$ and $r\in\{\sigma+1,4\}$,
\begin{equation*}
 M_a^{(r)}\cap(-M_a^{(r)})=\varnothing.
\end{equation*}
\end{lemma}

\begin{proof}
This is exactly Lemma~3.1 of \cite{DMMP2003}, with $m=(n+1)/2$ in the first case.  The hypotheses are precisely odd $n$ and $r=3^{(n+1)/2}+1$ or $r=4$; in the first case the proof also uses $\sigma^2=3q$.
\end{proof}

\begin{lemma}[The two blocks are injective]
\label{dmmp6:lem:block-injective}
Fix $b\ne0$.  The union of the two signs in the same-character block contains at most one interior solution of $D(x)=b$, and the union of the two signs in the opposite-character block also contains at most one.
\end{lemma}

\begin{proof}
Normalize $U=bz$, $V=b(z+1)$.  Since $\lambda$ is a power of $3$,
\[
 (z+1)^r-z^r=z^\lambda+z+1.
\]
Thus the same-character equation is
\begin{equation}
 1=s b^r(z^\lambda+z+1).
 \label{dmmp6:eq:same-block}
\end{equation}
The linearized map $L_\lambda(z):=z^\lambda+z$ is a permutation of $\mathbb{F}_{q}$.  For $\lambda=3$ this follows from $z^2=-1$ having no solution.  For $\lambda=\sigma$, a kernel element satisfies $z^\sigma=-z$; applying $\sigma$ and using $\sigma^2=3q$ on $\mathbb{F}_{q}$ gives $z^3=z$, hence $z\in\mathbb{F}_{3}$, and then $z=0$.

For each sign $s$, \eqref{dmmp6:eq:same-block} has therefore at most one $z$.  If both signs occurred, their solutions $z_+,z_-$ would satisfy
\[
 L_\lambda(z_++z_-)=1=L_\lambda(-1),
\]
so $z_++z_-=-1$.  Hence $b(z_++1)=-bz_-$.  These two nonzero quantities would both be squares, impossible because $-1$ is a nonsquare.

For the opposite-character block write
\[
 U=b(w+1),\qquad V=b(w-1).
\]
Then \eqref{dmmp6:eq:master-derivative} becomes
\begin{equation}
 s=b^r(w^r+1).
 \label{dmmp6:eq:opposite-block}
\end{equation}
For fixed $s$, \eqref{dmmp6:eq:gcd-r} shows that there are at most two roots $w,-w$.  The corresponding pairs $(U,V)$ are $(U,V)$ and $(-V,-U)$, so at most one pair consists of nonzero squares.  Finally, simultaneous admissible solutions for $s=1$ and $s=-1$ would give two opposite elements of $M_b^{(r)}$, contradicting Lemma~\ref{dmmp6:lem:dmmp-separation}.
\end{proof}

The special values do not create double fibers.

\begin{lemma}[Boundary targets]
\label{dmmp6:lem:boundary}
For all four branches,
\[
 N_d(0)=N_d(1)=N_d(-1)=1.
\]
\end{lemma}

\begin{proof}
Because $d$ is even,
\[
 D(-1-x)=-D(x).
\]
The zero fiber contains the fixed point $x=1$.  Any additional zero would occur with its distinct mate $-1-x$, contradicting the DMMP bound $N_d(b)\leqslant2$.  Thus $N_d(0)=1$.

For $b=1$, the boundary point $x=0$ is a solution.  In \eqref{dmmp6:eq:same-block}, $s=1$ gives $L_\lambda(z)=0$, hence $z=0$, while $s=-1$ gives $L_\lambda(z)=1$, hence $z=-1$; the latter would make $U=-1$, a nonsquare.  In \eqref{dmmp6:eq:opposite-block}, $s=1$ gives $w=0$, for which $V=-1$ is nonsquare, while $s=-1$ gives $w^r=1$, hence $w=\pm1$; one coordinate is then zero or nonsquare.  Thus $N_d(1)=1$, and the identity $D(-1-x)=-D(x)$ gives $N_d(-1)=1$.
\end{proof}

It follows that every double fiber consists of one same-character and one opposite-character interior solution.

\begin{proposition}[Parametrization of double fibers]
\label{dmmp6:prop:master-collision}
Let
\begin{equation}
 \mathcal C_\lambda:\qquad z^\lambda+z=w^r,
 \qquad r=\lambda+1,
 \label{dmmp6:eq:master-curve}
\end{equation}
and let $M_\lambda$ be the number of points $(z,w)\in\mathcal C_\lambda(\mathbb{F}_{q})$ satisfying
\begin{equation}
 \chi(z)=\chi(z+1)=\chi(w-1)=\chi(w+1)\in\{\pm1\}.
 \label{dmmp6:eq:master-marking}
\end{equation}
Then
\begin{equation}
 \omega_2=2M_\lambda.
 \label{dmmp6:eq:omega2-master}
\end{equation}
\end{proposition}

\begin{proof}
For the same-character member of a double fiber, write
\[
 U_x=bz,
 \qquad V_x=b(z+1),
\]
and for the opposite-character member write
\[
 U_y=b(w+1),
 \qquad V_y=b(w-1).
\]
Because all four $U$-coordinates are squares,
\begin{equation}
 \chi(z)=\chi(z+1)=\chi(w-1)=\chi(w+1)=\chi(b).
 \label{dmmp6:eq:four-markings}
\end{equation}
If $s=\chi(x)=\chi(x+1)$ and $\rho=\chi(y)=-\chi(y+1)$, equations \eqref{dmmp6:eq:same-block} and \eqref{dmmp6:eq:opposite-block} are equivalent to
\[
 s=b^r(z^\lambda+z+1),
 \qquad
 \rho=b^r(w^r+1).
\]
Thus $s=\rho$ gives \eqref{dmmp6:eq:master-curve}, whereas $s=-\rho$ gives
\begin{equation}
 z^\lambda+z+w^r=1.
 \label{dmmp6:eq:second-sign-curve}
\end{equation}
The involution
\[
 (z,w)\longmapsto(-z-1,-w)
\]
sends \eqref{dmmp6:eq:second-sign-curve} to \eqref{dmmp6:eq:master-curve} and flips all four character values in \eqref{dmmp6:eq:four-markings}.  Hence the two relative-sign types have equal marked counts.

It remains to verify exact reconstruction.  Given a marked point of \eqref{dmmp6:eq:master-curve}, put
\[
 A=z^\lambda+z+1=w^r+1.
\]
Then $A\ne0$, because $A=0$ would give $w^r=-1$, impossible since the image of $u\mapsto u^r$ is the subgroup of squares.  Put $s=\chi(A)$.  The element $s/A$ is a square, and \eqref{dmmp6:eq:gcd-r} implies that
\[
 b^r=\frac{s}{A}
\]
has exactly two roots $b,-b$.  They have opposite quadratic character, so \eqref{dmmp6:eq:master-marking} selects the unique root with $\chi(b)=\chi(z)$.  The four resulting square coordinates reconstruct the two derivative solutions uniquely by Lemma~\ref{dmmp6:lem:kummer-coordinate}.  Hence each relative-sign type is in bijection with the marked set, proving \eqref{dmmp6:eq:omega2-master}.
\end{proof}

\subsection{The elliptic character sum}

Put
\begin{equation*}
 J_n:=\sum_{x\in\mathbb{F}_{q}}\chi\bigl((x+1)(x^2+1)\bigr).
\end{equation*}
By the elementary point-count formula for \eqref{dmmp6:eq:E},
\begin{equation}
 J_n=-t_n.
 \label{dmmp6:eq:J-trace}
\end{equation}
We shall prove the single marked-count identity
\begin{equation*}
 8M_\lambda=q-7+2J_n
\end{equation*}
for both $\lambda=3$ and $\lambda=\sigma$.  

\subsection{The fixed-quartic calculation: \texorpdfstring{$\lambda=3$}{lambda=3}}

Here
\[
 \mathcal C_3:\qquad z^3+z=w^4.
\]
Its smooth projective model
\[
 Z^3W+ZW^3=T^4
\]
is a nonsingular plane quartic, hence has genus $3$.  Its marked point count nevertheless reduces to the elliptic trace $J_n$.

\begin{proposition}[Fixed-quartic evaluation]
\label{dmmp6:prop:B-evaluation}
For odd $n$,
\[
 8M_3=q-7+2J_n.
\]
\end{proposition}

\begin{proof}
Delete $(0,0)$ and $(-1,\pm1)$ from $z^3+z=w^4$; these are exactly the points where a marking factor vanishes.  On the remaining curve put
\[
 \varepsilon_1=\chi(z),\quad \varepsilon_2=\chi(z+1),\quad
 \varepsilon_3=\chi(w-1),\quad \varepsilon_4=\chi(w+1).
\]
The map $z\mapsto z^3+z$ is a permutation of $\mathbb F_q$, so the affine curve has $q$ points.  Expanding the four character projectors yields
\begin{equation}
 8M_3=q-3+\sum_{i<j}S_{ij}+S_{1234},
 \label{dmmp6:eq:B-expansion}
\end{equation}
where the sums are over the punctured curve.  Direct summation over $z$ and $w$, respectively, gives
\[
 S_{12}=J_n-2,\qquad S_{34}=0;
\]
the first uses the fiber size $1+\chi(z^3+z)$ and the deleted point $z=0$, while the second uses $\sum_w\chi(w^2-1)=-1$ and the same deleted point.

To evaluate the other pair terms, use the chart $Z=1$ and set
\[
 r_1=1+W,\qquad r_2=W-T,\qquad r_3=W+T,
 \qquad U=r_1r_2+r_2r_3+r_3r_1.
\]
The quartic equation gives $r_1+r_2+r_3=1$ and $r_1r_2r_3=U(U+1)$.  The three omitted chart points correspond to the permutations of $(1,0,0)$; thus the punctured curve is invariant under the permutations of $r_1,r_2,r_3$.  In these coordinates the six pair products sum to
\[
 \sum_{i<j}\varepsilon_i\varepsilon_j
 =\chi(r_1)-\chi(r_2)+\chi(r_3)
  -\chi(r_1r_2)+\chi(r_1r_3)-\chi(r_2r_3).
\]
By permutation symmetry, the sum of the three single-character terms is $S_{12}$, and the sum of the three pair-character terms is $S_{34}$.  Hence
\begin{equation}
 \sum_{i<j}S_{ij}=J_n-2.
 \label{dmmp6:eq:B-pair-total}
\end{equation}

For the fourfold product put $\eta(U)=\chi(-U(U+1))$.  Since
\[
 \varepsilon_1\varepsilon_2\varepsilon_3\varepsilon_4
 =\chi(-r_1r_2r_3)=\eta(U),
\]
the omitted points, whose product is zero, may be restored.  The ordered triple $(r_1,r_2,r_3)$ consists of the roots of
\[
 P_U(R)=R^3-R^2+UR-U(U+1),
\]
with discriminant $-U(U^2+1)$.  Writing $\rho(U)$ for its number of $\mathbb F_q$-roots and $N(U)$ for the number of ordered triples of roots, the three separable splitting types and the ramified case $U=0$ give
\[
 N(U)=2\rho(U)-1+\chi(-U(U^2+1)).
\]
Set $R_n=\sum_U\eta(U)\rho(U)$.  The standard quadratic sum and removal of the square factor $U^2$ give
\[
 \sum_U\eta(U)=1,
 \qquad \sum_U\eta(U)\chi(-U(U^2+1))=J_n-1.
\]
Consequently
\begin{equation}
 S_{1234}=2R_n+J_n-2.
 \label{dmmp6:eq:B-fourfold}
\end{equation}

We show that $R_n=0$ by an involution on the elliptic curve
\[
 E_0:\qquad V^2=X^3+X+1.
\]
Indeed, $P_U(X)=0$ is equivalent to
\[
 U^2+(1-X)U+X^2-X^3=0;
\]
with $V=1-X-U$ this gives the stated curve and $-U(U+1)=f(X,V)$, where $f=1-(X+V)^2$.  Thus
\[
 R_n=\sum_{(X,V)\in E_0(\mathbb F_q)_{\rm aff}}\chi(f(X,V)).
\]
Let $P_0=(0,1)$ and $\tau(P)=P_0-P$.  For $X\ne0$ the group law gives
\[
 X_\tau=\frac{X-1-V}{X^2},\qquad
 V_\tau=\frac{X^3+1+(X+1)V}{X^3}.
\]
For $B=(X+1)^2+V$ and $\bar B=(X+1)^2-V$, the curve equation implies
\[
 B\bar B=X^4,\qquad
 X_\tau+V_\tau=1+\frac{B}{X^3},\qquad
 X^3-B=-\frac{Bf}{X^2}.
\]
It follows that
\begin{equation}
 f\circ\tau=-\left(\frac{(X+1)^2+V}{X^4}\right)^2f.
 \label{dmmp6:eq:self-twist}
\end{equation}
The four points $O,(0,1),(0,-1),(1,0)$ are paired by $\tau$ as $O\leftrightarrow(0,1)$ and $(0,-1)\leftrightarrow(1,0)$; at the three affine points $f=0$.  On their complement the square factor in \eqref{dmmp6:eq:self-twist} is nonzero.  Since $\chi(-1)=-1$, the remaining character summands cancel in $\tau$-pairs.  Therefore $R_n=0$.

Equations~\eqref{dmmp6:eq:B-expansion}, \eqref{dmmp6:eq:B-pair-total}, and~\eqref{dmmp6:eq:B-fourfold} now yield
\[
 8M_3=q-3+(J_n-2)+(J_n-2)=q-7+2J_n.
\]
\end{proof}

\begin{remark}
Permutation symmetry reduces the marked count on this genus-three curve to the elliptic trace $J_n$ and the character sum $R_n$.  The self-twist \eqref{dmmp6:eq:self-twist} makes $R_n$ vanish.
\end{remark}

\subsection{Frobenius character sums for \texorpdfstring{$\lambda=\sigma$}{lambda=sigma}}

For $\lambda=\sigma$, define
\begin{equation}\label{dmmp6:eq:S1S2}
 \begin{aligned}
 S_1(n)&:=\sum_{z\in\mathbb{F}_{q}}\chi\bigl((z^\sigma+z)z(z+1)\bigr),\\
 S_2(n)&:=\sum_{z\in\mathbb{F}_{q}}\chi\bigl((z^\sigma+z)(z-1)\bigr).
 \end{aligned}
\end{equation}

\begin{theorem}[Frobenius--Dickson elliptic trace theorem]
\label{dmmp6:thm:semilinear-dickson}
Let $n\geqslant1$ be odd, $q=3^n$, $\sigma=3^{(n+1)/2}$, and let $\chi$ be the quadratic character of $\mathbb{F}_{q}$, extended by $\chi(0)=0$.  Let
\[
 E/\mathbb{F}_{3}:\qquad Y^2=(X+1)(X^2+1),
 \qquad t_n=q+1-\#E(\mathbb{F}_{q}).
\]
Then the two character correlations in \eqref{dmmp6:eq:S1S2} satisfy
\begin{equation}
 S_1(n)=S_2(n)=
 \sum_{x\in\mathbb{F}_{q}}\chi\bigl((x+1)(x^2+1)\bigr)-1
 =-t_n-1.
 \label{dmmp6:eq:semilinear-dickson-trace}
\end{equation}
In particular, their common value is governed by the fixed elliptic curve $E/\mathbb{F}_{3}$, although the defining exponent $\sigma=3^{(n+1)/2}$ grows with $n$.
\end{theorem}

Let $D_j(X)=D_j(X,1)$ denote the \emph{Dickson polynomial}
characterized by $D_j(\eta+\eta^{-1})=\eta^j+\eta^{-j}$.
The \emph{Coulter--Matthews exponents} are $(3^\ell+1)/2$ with
$\ell$ a positive odd integer and $\gcd(\ell,n)=1$ over $\mathbb F_{3^n}$
\cite[Theorem~4.1]{CoulterMatthews1997}.
We prove Theorem~\ref{dmmp6:thm:semilinear-dickson} below by a M\"obius symmetry, reduction to a Coulter--Matthews exponent, and a Dickson square-class conjugacy.  We first record how these correlations arise from the marked collision problem.  The curve
\[
 \mathcal C_\sigma:\qquad z^\sigma+z=w^{\sigma+1}
\]
need not be normalized.

\begin{lemma}[Elimination of the $w$-coordinate]
\label{dmmp6:lem:A-elimination}
Let
\[
 \mu_n:=
 \#\Bigl\{z\in\mathbb{F}_{q}:
 \chi(z^\sigma+z)=
 \chi(z(z+1))=
 \chi(z-1)=1
 \Bigr\}.
\]
Then
\begin{equation}
 M_\sigma=\mu_n.
 \label{dmmp6:eq:M-mu}
\end{equation}
\end{lemma}

\begin{proof}
Set $L(z)=z^\sigma+z$.  If $w^{\sigma+1}=L(z)$, apply $\sigma$ and use $\sigma^2=3q$ on $\mathbb{F}_{q}$ to get
\[
 w^2=\frac{z^3+z^\sigma}{z^\sigma+z}.
\]
Hence
\begin{equation}
 w^2-1=\frac{z(z-1)(z+1)}{L(z)}.
 \label{dmmp6:eq:w2minus1}
\end{equation}
Because $\gcd(\sigma+1,q-1)=2$, a nonzero $L(z)$ has exactly two roots $w,-w$ precisely when $\chi(L(z))=1$.  The condition $\chi(z(z+1))=1$ is exactly $\chi(z)=\chi(z+1)$.  Under these two conditions, \eqref{dmmp6:eq:w2minus1} gives
\[
 \chi(w^2-1)=\chi(z-1).
\]
Thus $\chi(z-1)=1$ is equivalent to $\chi(w-1)=\chi(w+1)$.  Passing from $w$ to $-w$ flips both latter character values, so exactly one of the two roots has the same common character as $z,z+1$.  This proves \eqref{dmmp6:eq:M-mu}.
\end{proof}

\begin{lemma}[Minimal character expansion]
\label{dmmp6:lem:A-expansion}
One has
\begin{equation}
 8M_\sigma=q-5+S_1(n)+S_2(n).
 \label{dmmp6:eq:A-expansion}
\end{equation}
\end{lemma}

\begin{proof}
For $A(z)=z^\sigma+z$, $B(z)=z(z+1)$, and $C(z)=z-1$, the indicator of the three conditions in Lemma~\ref{dmmp6:lem:A-elimination} is
\[
 \frac18(1+\chi(A))(1+\chi(B))(1+\chi(C))
\]
away from the zeros.  At $z=-1$ this expression contributes the spurious value $4/8$; $z=0,1$ contribute zero.  Hence
\[
 8M_\sigma=
 \sum_z(1+\chi(A))(1+\chi(B))(1+\chi(C))-4.
\]
Now
\[
 \sum_z\chi(A)=0,\qquad
 \sum_z\chi(B)=-1,\qquad
 \sum_z\chi(C)=0,
\]
while
\[
 \sum_z\chi(BC)=\sum_z\chi(z^3-z)=0
\]
by $z\mapsto-z$.  Finally
\[
 \sum_z\chi(ABC)=0,
\]
because for $z\ne0$ the substitution $z\mapsto z^{-1}$ multiplies $ABC$ by $-z^{-(\sigma+5)}$, whose nonconstant factor is a square.  This proves \eqref{dmmp6:eq:A-expansion}.
\end{proof}

\begin{lemma}[M\"obius symmetry of the character correlations]
\label{dmmp6:lem:S1S2-symmetry}
For every odd $n$,
\begin{equation*}
 S_1(n)=S_2(n).
\end{equation*}
\end{lemma}

\begin{proof}
Substitute $z=-u^{-1}$ in $S_2(n)$.  For $u\ne0$,
\[
 ( (-u^{-1})^\sigma-u^{-1})(-u^{-1}-1)
 =\frac{(u^\sigma+u)(u+1)}{u^{\sigma+2}}.
\]
Compared with $(u^\sigma+u)u(u+1)$, the ratio is $u^{-(\sigma+3)}$, a square because $\sigma+3$ is even.  The terms at $u=0$ and the corresponding omitted point contribute zero, so there is no boundary correction.  Hence $S_1(n)=S_2(n)$.
\end{proof}

It remains to evaluate one Frobenius character sum.  A Dickson polynomial substitution reduces it to the fixed elliptic trace $J_n$.

\paragraph{Dickson reduction.}

Let
\[
 H:=(\mathbb{F}_{q}^\times)^2,
 \qquad h:=|H|=\frac{q-1}{2},
 \qquad e:=\sigma-1.
\]
For $s$ coprime to $h$, define
\begin{equation*}
 C_s:=\sum_{r\in H}\chi\bigl((r+1)(r^s+1)\bigr).
\end{equation*}

\begin{lemma}[Reduction to a Coulter--Matthews type exponent]
\label{dmmp6:lem:CM-reduction}
There is an odd integer $\ell$ with $\gcd(\ell,n)=1$ such that, for
\[
 d_0:=\frac{3^\ell+1}{2},
\]
one has
\begin{equation}
 S_2(n)=C_{d_0}.
 \label{dmmp6:eq:S2Cd}
\end{equation}
More precisely, $\ell=m$ if $m$ is odd and $\ell=m-1$ if $m$ is even.
\end{lemma}

\begin{proof}
Since $z^\sigma+z=z(z^e+1)$, pairing the two elements $r,-r$ with $r\in H$ gives
\[
 S_2(n)=
 \sum_{r\in H}\chi((r^e+1)(r-1))
 +\sum_{r\in H}\chi((r^e+1)(r+1)).
\]
The first sum is zero under $r\mapsto r^{-1}$: its integrand changes by the nonsquare factor $-1$, and the only fixed point of inversion in the odd-order group $H$ is $r=1$, where the summand is zero.  Hence $S_2(n)=C_e$.

Put $\rho=(\sigma+1)/2$.  Since
\[
 e\rho=\frac{\sigma^2-1}{2}=3h+1,
\]
we have $e^{-1}\equiv\rho\pmod h$.  The substitutions $r\mapsto r^s$ and Frobenius invariance of $\chi$ give
\begin{equation}
 C_s=C_{s^{-1}},\qquad C_{3^js}=C_s.
 \label{dmmp6:eq:C-invariances}
\end{equation}
Thus $C_e=C_\rho$.  If $m$ is odd, take $\ell=m$.  If $m$ is even, take $\ell=m-1$; then
\[
 3^{m-1}\rho\equiv\frac{3^{m-1}+1}{2}\pmod h,
\]
so \eqref{dmmp6:eq:C-invariances} again gives \eqref{dmmp6:eq:S2Cd}.  In either case $\ell$ is odd and $\gcd(\ell,n)=1$.
\end{proof}

\begin{lemma}[M\"obius transfer]
\label{dmmp6:lem:mobius-transfer}
For $d_0=(3^\ell+1)/2$ as above, put
\[
 A_{d_0}:=\sum_{x\in\mathbb{F}_{q}}\chi\bigl((x+1)D_{d_0}(x)\bigr).
\]
Then
\begin{equation}
 C_{d_0}=A_{d_0}-1.
 \label{dmmp6:eq:C-vs-A}
\end{equation}
\end{lemma}

\begin{proof}
The M\"obius transformation
\[
 r=-\frac{x-1}{x+1}
\]
gives a bijection
\[
 H\longleftrightarrow
 \{x\in\mathbb{F}_{q}:\chi(x^2-1)=-1\}.
\]
Indeed $\chi(r)=1$ is equivalent to $\chi(x^2-1)=-1$.  Since $d_0=(3^\ell+1)/2$, characteristic $3$ gives the polynomial identity
\begin{equation*}
 (x-1)^{d_0}+(x+1)^{d_0}=-D_{d_0}(x).
\end{equation*}
To verify it, write $x=\eta+\eta^{-1}$; then
\[
 x-1=\frac{(\eta+1)^2}{\eta},
 \qquad
 x+1=\frac{(\eta-1)^2}{\eta},
\]
and use Frobenius on $(\eta\pm1)^{3^\ell+1}$.  It follows that the summand in $C_{d_0}$ becomes
\[
 \chi\bigl((x+1)D_{d_0}(x)\bigr).
\]

Write this last character as $A(x)$.  Since $A(-1)=0$ and $A(1)=1$,
\[
 2C_{d_0}
 =A_{d_0}-\sum_x\chi(x^2-1)A(x)-1.
\]
If
\[
 B_{d_0}=\sum_x\chi((x-1)D_{d_0}(x)),
\]
then the point $x=-1$ gives the sole boundary correction
\[
 \sum_x\chi(x^2-1)A(x)=B_{d_0}+1.
\]
Because $d_0$ is even, $D_{d_0}(-x)=D_{d_0}(x)$, and $x\mapsto-x$ gives $B_{d_0}=-A_{d_0}$.  Therefore $2C_{d_0}=2A_{d_0}-2$, proving \eqref{dmmp6:eq:C-vs-A}.
\end{proof}

\begin{lemma}[Dickson square-class conjugacy]
\label{dmmp6:lem:Dickson-conjugacy}
Put
\[
 a=\frac{d_0}{2}=\frac{3^\ell+1}{4}.
\]
Then $D_a$ permutes $\mathbb{F}_{q}$, and
\begin{equation}
 \chi(D_a(x)+1)=\chi(x+1)
 \qquad(x\in\mathbb{F}_{q}).
 \label{dmmp6:eq:square-class-conjugacy}
\end{equation}
Consequently
\begin{equation*}
 A_{d_0}=J_n.
\end{equation*}
\end{lemma}

\begin{proof}
Since $\ell$ is odd, $3^\ell\equiv3\pmod8$, so $a$ is odd.  Moreover
\begin{equation}
 \gcd(a,q^2-1)=1.
 \label{dmmp6:eq:gcd-a}
\end{equation}
Indeed, if an odd prime $p$ divided both numbers, then $3^\ell\equiv-1\pmod p$, while the order of $3$ modulo $p$ would divide $2n$.  It also divides $2\ell$ but not $\ell$.  Since $\gcd(\ell,n)=1$ and both are odd, that order would be $2$, forcing $p\mid4$, a contradiction.  The factor $2$ is excluded because $a$ is odd.

For every $x\in\mathbb{F}_{q}$, choose $\eta\in\mathbb F_{q^2}^\times$ with $x=\eta+\eta^{-1}$.  If $D_a(x)=D_a(y)$, then the unordered pairs $\{\eta^a,\eta^{-a}\}$ associated with $x$ and $y$ agree.  By \eqref{dmmp6:eq:gcd-a}, the power map on $\mathbb F_{q^2}^\times$ is bijective, so $x=y$.  Thus $D_a$ permutes $\mathbb{F}_{q}$.

For odd $a$ one has the exact polynomial identity
\begin{equation}\label{dmmp6:eq:Dickson-square}
 D_a(X)+1=(X+1)H_a(X)^2,
\end{equation}
where
\[
 H_a(X)=1+\sum_{j=1}^{(a-1)/2}D_j(X).
\]
Indeed, with $X=\eta+\eta^{-1}$,
\[
 H_a(X)=\sum_{j=-(a-1)/2}^{(a-1)/2}\eta^j
 =\eta^{(1-a)/2}\frac{\eta^a-1}{\eta-1},
\]
and $X+1=(\eta-1)^2/\eta$ in characteristic $3$.  Since $D_a$ is a permutation and $D_a(-1)=-1$, the only zero of $D_a(X)+1$ in $\mathbb{F}_{q}$ is $X=-1$; hence \eqref{dmmp6:eq:Dickson-square} implies \eqref{dmmp6:eq:square-class-conjugacy} for every $x$.

Finally $d_0=2a$ and Dickson composition gives
\[
 D_{d_0}=D_2\circ D_a,
 \qquad
 D_2(X)=X^2+1
\]
in characteristic $3$.  Therefore
\[
 \begin{aligned}
 A_{d_0}
 &=\sum_x\chi\bigl((D_a(x)+1)D_2(D_a(x))\bigr)\\
 &=\sum_y\chi\bigl((y+1)(y^2+1)\bigr)=J_n,
 \end{aligned}
\]
because $D_a$ permutes $\mathbb{F}_{q}$.
\end{proof}

\begin{proof}[Proof of Theorem~\ref{dmmp6:thm:semilinear-dickson}]
By Lemmas~\ref{dmmp6:lem:CM-reduction}, \ref{dmmp6:lem:mobius-transfer}, and~\ref{dmmp6:lem:Dickson-conjugacy},
\[
 S_2(n)=C_{d_0}=A_{d_0}-1=J_n-1.
\]
By Lemma~\ref{dmmp6:lem:S1S2-symmetry}, $S_1(n)=S_2(n)$.  Finally, the elementary point-count identity \eqref{dmmp6:eq:J-trace} gives $J_n=-t_n$.  Therefore
\[
 S_1(n)=S_2(n)=J_n-1=-t_n-1,
\]
which is \eqref{dmmp6:eq:semilinear-dickson-trace}.
\end{proof}

\begin{proposition}[Frobenius collision evaluation]
\label{dmmp6:prop:A-evaluation}
For odd $n$,
\[
 8M_\sigma=q-7+2J_n.
\]
\end{proposition}

\begin{proof}
By Lemma~\ref{dmmp6:lem:A-expansion} and Theorem~\ref{dmmp6:thm:semilinear-dickson},
\[
 8M_\sigma=q-5+2(J_n-1)=q-7+2J_n.
\]
\end{proof}

\subsection{Completion of the proof}

\begin{proof}[Proof of Theorem~\ref{dmmp6:thm:main}]
Apply Proposition~\ref{dmmp6:prop:B-evaluation} to the $B$-branches and Proposition~\ref{dmmp6:prop:A-evaluation} to the $A$-branches.  In all four cases
\[
 8M_\lambda=q-7+2J_n=q-7-2t_n
\]
by \eqref{dmmp6:eq:J-trace}.  The double-fiber formula \eqref{dmmp6:eq:omega2-master} gives
\[
 \omega_2=2M_\lambda=\frac{q-7-2t_n}{4}.
\]
Equation \eqref{dmmp6:eq:few-support} then yields \eqref{dmmp6:eq:main-spectrum}.

Directly, $\#E(\mathbb{F}_{3})=6$, so $t_1=-2$.  If $\alpha,\beta$ are the Frobenius eigenvalues of $E/\mathbb{F}_{3}$, then
\[
 \alpha+\beta=-2,
 \qquad
 \alpha\beta=3,
 \qquad
 t_n=\alpha^n+\beta^n,
\]
and therefore \eqref{dmmp6:eq:trace-recurrence}; this is the standard elliptic Frobenius recurrence, cf.~\cite[Chap.~V]{Silverman2009}.

The recurrence gives $\omega_2=0$ for $n=1,3$.  For odd $n\geqslant5$, Hasse's bound $|t_n|\leqslant2\sqrt q$ implies
\[
 q-7-2t_n\geqslant q-7-4\sqrt q>0,
\]
so $\omega_2>0$.  Since DMMP already gives $N_d(b)\leqslant2$, the differential uniformity is exactly two.
\end{proof}

\paragraph{Finite-field examples.}
The supplementary material gives sample spectra and direct finite-field computations; these are not used in the proof for arbitrary odd extension degrees.

\section{The Yan--Han Family: Complete Spectrum from Artin Factors}\label{sec:yanhan-screen}

Let
\[
 q=3^m,\qquad m\geqslant3\text{ odd},\qquad d:=\frac{q-3}{4},
\]
and put
\[
 D_m(x):=(x+1)^d-x^d.
\]
Let $\chi$ denote the quadratic character of $\mathbb F_q$, extended by $\chi(0)=0$.  Since $q\equiv3\pmod8$, one has $\chi(-1)=-1$ and $d$ is even.  Yan and Han proved that every equation $D_m(x)=b$ has at most three solutions \cite{YanHan2019}.  Their proof separates the regular inputs $x\ne0,-1$ according to the four pairs $(\chi(x+1),\chi(x))\in\{1,-1\}^2$.  The exponent congruence
\[
 4d=q-3\equiv-2\pmod{q-1}
\]
shows that this is a fixed quartic-Kummer problem rather than a cover whose effective degree grows with $q$.

\begin{lemma}[Quartic Kummer reduction]\label{lem:yanhan-kummer}
For every $z\in\mathbb F_q^\times$,
\[
 (z^d)^2=\frac{\chi(z)}z.
\]
For a regular input $x$, put $s=\chi(x+1)$ and $t=\chi(x)$.  There are unique \emph{square elements} $\alpha,\beta\in\mathbb F_q^\times$ such that
\[
 x+1=s\alpha^2,\qquad x=t\beta^2.
\]
They satisfy
\[
 (x+1)^d=\alpha^{-1},\qquad x^d=\beta^{-1},
\]
and hence
\[
 s\alpha^2-t\beta^2=1,\qquad D_m(x)=\alpha^{-1}-\beta^{-1}.
\]
\end{lemma}

\begin{proof}
The exponent identity $2d=(q-1)/2-1$ gives the first formula.  For $z\ne0$, the element $\chi(z)z$ is a square because $\chi(-1)=-1$.  Its two square roots differ by $-1$, so exactly one of them is itself a square; call it $\gamma$.  Then $z=\chi(z)\gamma^2$.  Since $d$ is even and $\gamma$ is a square,
\[
 z^d=\gamma^{2d}=\gamma^{(q-3)/2}=\gamma^{-1}.
\]
Apply this construction to $z=x+1$ and $z=x$.
\end{proof}

Call the four regular input strata I--IV in the order
\[
 (s,t)=(1,1),\ (-1,-1),\ (1,-1),\ (-1,1),
\]
and let $A,B,C$ be the corresponding target images
\[
 A:=\operatorname{Im}(\mathrm I)\cup\operatorname{Im}(\mathrm{II}),
 \qquad B:=\operatorname{Im}(\mathrm{III}),
 \qquad C:=\operatorname{Im}(\mathrm{IV}).
\]
The argument in the proof of \cite[Th.~1]{YanHan2019}, specifically its four quadratic-character cases, shows that I and II together contribute at most one solution to any regular differential fiber, and III and IV contribute at most one each.  Thus the differential map is injective on each of the three indicated input blocks.

\begin{lemma}[Sizes of the three input blocks]\label{lem:yanhan-block-sizes}
The four regular strata have sizes
\[
 |\mathrm I|=|\mathrm{II}|=|\mathrm{III}|=\frac{q-3}{4},
 \qquad
 |\mathrm{IV}|=\frac{q+1}{4}.
\]
Consequently
\[
 |A|=\frac{q-3}{2},\qquad |B|=\frac{q-3}{4},\qquad |C|=\frac{q+1}{4}.
\]
\end{lemma}

\begin{proof}
For $s,t\in\{\pm1\}$, the number of regular $x$ with $(\chi(x+1),\chi(x))=(s,t)$ equals
\[
 \frac14\sum_{x\ne0,-1}(1+s\chi(x+1))(1+t\chi(x)).
\]
Using
\[
 \sum_{x\ne0,-1}\chi(x+1)=-1,\qquad
 \sum_{x\ne0,-1}\chi(x)=1,
\]
and $\sum_x\chi(x(x+1))=-1$ gives the four displayed values.  Injectivity on the three input blocks then gives the target-image cardinalities.
\end{proof}

The conics in Lemma~\ref{lem:yanhan-kummer} admit the following rational parameters.  In Cases I and II put $u=\alpha-\beta$.  Then respectively
\begin{align*}
 (\alpha,\beta)&=\left(\frac{u+u^{-1}}2,\frac{u^{-1}-u}2\right),\\
 (\alpha,\beta)&=\left(\frac{u-u^{-1}}2,-\frac{u+u^{-1}}2\right).
\end{align*}
Both cases yield the same unmarked target map.  For Cases III and IV one may take
\[
 \alpha_B=\frac{v^2-1}{v^2+1},\qquad
 \beta_B=\frac{v}{v^2+1},
\]
\[
 \alpha_C=\frac{w^2+w-1}{w^2+1},\qquad
 \beta_C=\frac{-w^2+w+1}{w^2+1}.
\]
Substitution gives
\begin{equation}\label{eq:yanhan-targetmaps}
 \begin{aligned}
 b_A(u)&:=\frac{u^3}{u^4-1},\qquad
 b_C(w):=\frac{1-w^4}{1+w^4},\\
 b_B(v)&:=\frac{-v^4+v^3+v+1}{v^3-v}.
 \end{aligned}
\end{equation}
The actual strata retain the square markings from Lemma~\ref{lem:yanhan-kummer}; the next lemmas prove that the unmarked rational images needed later introduce no spurious finite targets.

\begin{lemma}[Target-only description of $C$]\label{lem:yanhan-C}
For $b\in\mathbb F_q$,
\[
 b\in C\quad\Longleftrightarrow\quad
 \chi(1-b)=\chi(1+b)=1.
\]
In particular, $0\in C$ and $\pm1\notin C$.
\end{lemma}

\begin{proof}
From $b=b_C(w)$,
\[
 1+b=-\frac1{1+w^4},\qquad 1-b=-\frac{w^4}{1+w^4}.
\]
The two characters are therefore equal.  If the point is actually in the marked Case~IV stratum, then $\alpha_C$ and $\beta_C$ are both squares, so their product
\[
 \alpha_C\beta_C=-\frac{1+w^4}{(1+w^2)^2}
\]
is a square.  Because $-1$ is nonsquare, this implies $\chi(1+w^4)=-1$, and hence $\chi(1\pm b)=1$.

Conversely, suppose both target characters are one.  Then $(1-b)/(1+b)$ is a square.  Since $\gcd(4,q-1)=2$, fourth powers and squares coincide in $\mathbb F_q^\times$, so choose $w$ with
\[
 w^4=\frac{1-b}{1+b}.
\]
The denominator $1+w^2$ is nonzero because $-1$ is nonsquare.  The condition $\chi(1+w^4)=-1$ now shows only that $\alpha_C(w)$ and $\beta_C(w)$ have the same quadratic character.  The two choices $w$ and $-w$ satisfy
\[
 (\alpha_C(-w),\beta_C(-w))=(-\beta_C(w),-\alpha_C(w)).
\]
Multiplication by the nonsquare $-1$ flips both characters simultaneously.  Hence exactly one of the two choices has both marked coordinates square.  The endpoint statements follow directly.
\end{proof}

\begin{lemma}[The cyclic quartic block]\label{lem:yanhan-B}
The map $b_B:\mathbb P^1_v\to\mathbb P^1_b$ is a cyclic cover of degree four over $\mathbb F_3$.  The M\"obius transformation
\[
 \sigma(v)=\frac1{v+1}
\]
has order four and satisfies $b_B\circ\sigma=b_B$.  The quadratic intermediate coordinate
\[
 t=\frac{v^2+1}{v-1}
\]
satisfies
\begin{equation}\label{eq:yanhan-B-intermediate}
 b=\frac{t(1-t)}{t+1},\qquad v^2-tv+t+1=0.
\end{equation}
The two geometric branch values are the roots of $b^2+1$, and each is totally ramified of index four; $b=\infty$ is unramified.  In particular there is no rational branch value over an odd extension.  Furthermore
\[
 b_B(\mathbb P^1(\mathbb F_q))\setminus\{\infty\}=B,
\]
and each $b\in B$ has exactly four rational $v$-preimages.
\end{lemma}

\begin{proof}
Substitution gives $b_B\circ\sigma=b_B$ and
$\sigma^2(v)=(v+1)/(v-1)$, so $\sigma$ has order four and generates
the deck group.  Every nonidentity element fixes precisely the two
roots of $v^2+v-1$; at either root, $b_B(v)=1-v$ and
$b_B(v)^2+1=0$.  Thus these are the only ramification points,
each of index four.  Eliminating $v$ gives
\eqref{eq:yanhan-B-intermediate}.

The fiber over $\infty$ consists of $0,1,-1,\infty$.
Over an odd extension the branch values are not rational, and a
rational point in a finite fiber has four distinct rational deck
translates.  The finite rational image therefore has size
$(q+1-4)/4=(q-3)/4$.  It contains $B$ by the marked parametrization
and equals $B$ by Lemma~\ref{lem:yanhan-block-sizes}.
\end{proof}

\begin{lemma}[Two-to-one rational fibers of $b_A$]\label{lem:yanhan-A}
Let $u\in\mathbb F_q^\times$ with $u^4\ne1$.  Then $b_A(X)=b_A(u)$ has exactly two solutions in $\mathbb F_q$.  Consequently
\[
 b_A(\mathbb P^1(\mathbb F_q))=A\cup\{0,\infty\};
\]
every $b\in A$ has exactly two unmarked rational $u$-preimages, $b=0$ has the two preimages $0,\infty$, and $b=\infty$ has the two rational poles $u=\pm1$.
\end{lemma}

\begin{proof}
After cross multiplication and removal of the root $X=u$, the remaining factor is
\[
 u^3X^3+X^2+uX+u^2.
\]
Writing $X=u+y$ and $s=u/y$ transforms the remaining equation into
\[
 s^3+u^{-4}s=-1.
\]
The $\mathbb F_3$-linear map $L_u(s)=s^3+u^{-4}s$ has trivial kernel: a nonzero kernel element would satisfy $s^2=-u^{-4}$, impossible because $u^{-4}$ is a square and $-1$ is nonsquare.  Thus $L_u$ is bijective, giving exactly one companion root.

The two special fibers $b=0$ and $b=\infty$ are as stated.  Removing their four rational source points leaves $q-3$ points, and every remaining rational image fiber has size two; hence the finite nonzero image has size $(q-3)/2$.  The Case I/II parametrization gives $A$ as a subset of this image, while Lemma~\ref{lem:yanhan-block-sizes} gives $|A|=(q-3)/2$.  Therefore the two sets coincide.
\end{proof}

Yan and Han also show that the targets $b=\pm1$ have only the two exceptional solutions $x=0,-1$, respectively, and no regular solution \cite[Lem.~2]{YanHan2019}.  Table~\ref{tab:yanhan-block-dictionary} summarizes the resulting blocks and their lift multiplicities.

\begin{table*}[t]
\caption{Yan--Han input strata and unmarked target models.}\label{tab:yanhan-block-dictionary}
\centering\footnotesize
\begin{tabularx}{\textwidth}{@{}p{0.13\textwidth}p{0.20\textwidth}p{0.26\textwidth}p{0.14\textwidth}X@{}}
\toprule
Target block & Input stratum & Unmarked model & Target size & Lift multiplicity used below \\
\midrule
$A$ & I $\cup$ II & $b_A(u)=u^3/(u^4-1)$ & $(q-3)/2$ & $2$ for $b\in A$ \\
$B$ & III & $b_B(v)$ in \eqref{eq:yanhan-targetmaps} & $(q-3)/4$ & $4$ for $b\in B$ \\
$C$ & IV & $\chi(1-b)=\chi(1+b)=1$ & $(q+1)/4$ & $2\times2$ in the target Kummer lift $z_\pm^2=1\pm b$ \\
\bottomrule
\end{tabularx}
\end{table*}

\subsection{Recovering the spectrum from image intersections}

Put
\[
 P_{AB}:=|A\cap B|,\qquad P_{AC}:=|A\cap C|,\qquad P_{BC}:=|B\cap C|,
\]
\[
 P:=P_{AB}+P_{AC}+P_{BC},\qquad T:=|A\cap B\cap C|.
\]

\begin{proposition}[Intersection reconstruction]\label{prop:yanhan-reconstruction}
For every odd $m\geqslant3$,
\begin{equation}\label{eq:yanhan-reconstruction}
 \begin{aligned}
 \omega_0&=P-T, & \omega_1&=q-2P+3T,\\
 \omega_2&=P-3T, & \omega_3&=T.
 \end{aligned}
\end{equation}
and $\omega_i=0$ for $i\geqslant4$.  Equivalently,
\[
 M_2=2P,\qquad M_3=6T.
\]
\end{proposition}

\begin{proof}
On $\mathcal T=\mathbb F_q\setminus\{\pm1\}$, the block images have
$I_0=I_1=q-2$ by Lemma~\ref{lem:yanhan-block-sizes},
$I_2=P$, and $I_3=T$.  Proposition~\ref{prop:block-incidence}
with $h=1$ gives the regular frequencies
$N_0=P-T$, $N_1=q-2-2P+3T$, $N_2=P-3T$, and $N_3=T$.
The two exceptional targets each have one preimage, yielding
\eqref{eq:yanhan-reconstruction}; ordered pairs and triples
then give $M_2=2P$ and $M_3=6T$.
\end{proof}

The support bound $\{0,1,2,3\}$ requires no moment beyond order three.  To prove $\Delta=3$, it suffices to show $T>0$; the complete spectrum requires the exact values of both $P$ and $T$.

\subsection{Point counts for pair and triple intersections}\label{sec:yanhan-exact}

\subsubsection{The three pair intersections}

Let $Y$ be the smooth projective normalization of
\begin{equation}\label{eq:yanhan-Y}
 b_A(u)=b_B(v).
\end{equation}
Write
\[
 A_0:=u^4-1,\qquad B_+:=u^4+u^3-1,\qquad B_-:=u^4-u^3-1,
\]
\[
 H:=u^8+u^6+u^4+1.
\]
Since
\[
 b_A(u)^2+1=\frac{H}{A_0^2}
\]
and $H$ has simple zeros, its square class remains nontrivial in $\overline{\mathbb F}_3(u)^\times/\overline{\mathbb F}_3(u)^{\times2}$.  Hence the unique quadratic subfield of the cyclic $C_4$ cover $b_B$ is not contained even after base change to $\overline{\mathbb F}_3(u)$.  Any nontrivial intersection with a cyclic quartic extension contains that quadratic subfield, so the pullback \eqref{eq:yanhan-Y} is geometrically connected of degree four over the $u$-line.

A direct differentiation gives
\[
 b_A'(u)=-\frac{u^6}{(u^4-1)^2}.
\]
Locally $b_A(u)=-u^3+O(u^7)$ at $u=0$, so this is a wild point of index three above $b=0$.  All four roots of $A_0$ are simple poles, and at infinity, with $r=1/u$, one has $b_A=r/(1-r^4)$.  Thus the only geometric branch value of $b_A$ is $0$.  The two branch values $b^2+1=0$ of $b_B$ are disjoint from it, and each has four unramified $u$-preimages.  At every point of the projective fiber product, one of the two maps to the $b$-line is \'etale; base change therefore shows that this fiber product is smooth.  Since it is geometrically connected, it is already $Y$.  In particular, for every $n\geqslant1$,
\[
 \#Y(\mathbb F_{3^n})
 =\sum_{b\in\mathbb P^1(\mathbb F_{3^n})}N_A(b)N_B(b),
\]
where $N_A,N_B$ count projective rational preimages, including those at infinity.  The cover $Y\to\mathbb P^1_u$ has eight tame ramification points of index four.  Riemann--Hurwitz gives
\begin{equation*}
 2g(Y)-2=4(-2)+8(4-1)=16,
 \qquad g(Y)=9.
\end{equation*}

For odd $m$, a target in $A\cap B$ has $2\cdot4=8$ rational lifts.  The only rational projective boundary contribution is $b=\infty$: $b_A$ has the two poles $u=\pm1$ and $b_B$ has the four poles $v=0,\pm1,\infty$.  At $b=0$ the $B$-fiber has no rational point because
\[
 -v^4+v^3+v+1=-(v^2+1)(v^2-v-1)
\]
and both quadratic factors have nonsquare discriminant $-1$ in odd extensions.  The values $\pm1$ are not in the finite $B$-image, and the two branch values $b^2+1=0$ are not rational over an odd extension.  Hence
\begin{equation}\label{eq:PAB}
 P_{AB}=\frac{\#Y(\mathbb F_q)-8}{8}
 =\frac{q-7-\tau_Y(m)}8,
\end{equation}
where $\tau_Y(m):=\operatorname{Tr}(\operatorname{Frob}_3^m|H^1(Y_{\overline{\mathbb F}_3},\mathbb Q_\ell))$.

\begin{lemma}[Smooth low-genus models]\label{lem:yanhan-smooth-models}
The polynomials defining the curves $C_1,C_2,E,G_1,G_2$ below are squarefree over $\mathbb F_3$.  Their smooth projective genera are respectively $3,3,1,2,2$.  The even-degree models $C_1,C_2,G_2$ are monic and therefore have two rational points at infinity over every extension.  The auxiliary quartic
\[
 E_0:y^2=-(t^2+1)(t^2+t-1)
\]
is squarefree of genus one; its leading coefficient is $-1$, so it has no rational point at infinity over an odd extension.
\end{lemma}

\begin{proof}
Taking gcds with the formal derivatives in $\mathbb F_3[u]$ or $\mathbb F_3[t]$ gives one in every case.  For a squarefree model $y^2=f(x)$ in odd characteristic, the smooth projective genus is $\lfloor(\deg f-1)/2\rfloor$.  An odd degree gives one point at infinity; an even degree gives two rational points or none according as the leading coefficient is a square or a nonsquare.  Apply these formulas to degrees $8,8,3,5,6,4$ and the displayed leading coefficients.
\end{proof}

For $A\cap C$, define the genus-three curves
\begin{equation*}
 C_1:y^2=B_+A_0,
 \qquad
 C_2:y^2=B_+B_-,
\end{equation*}
and let $\tau_1(m),\tau_2(m)$ denote their Frobenius traces.

\begin{lemma}[The $A$--$C$ count]\label{lem:PAC}
For every odd $m$,
\begin{equation}\label{eq:PAC}
 P_{AC}=\frac{q-7-2\tau_1(m)-\tau_2(m)}8.
\end{equation}
\end{lemma}

\begin{proof}
For $u\ne\pm1$, Lemma~\ref{lem:yanhan-C} gives the indicator
\[
 \frac14\bigl(1+\chi(1+b_A(u))\bigr)\bigl(1+\chi(1-b_A(u))\bigr).
\]
Since
\[
 \chi(1+b_A(u))=\chi(B_+A_0),\qquad
 \chi(1-b_A(u))=\chi(B_-A_0),
\]
the two single-character sums agree under $u\mapsto-u$ and are governed by $C_1$, while the product is governed by $C_2$.  Both even-degree monic models have two rational points at infinity, so their affine character sums are $-\tau_i(m)-1$.  At $u=\pm1$ the single characters vanish and the product character equals $-1$.  The value $u=0$ represents one of the two unmarked preimages of the boundary target $b=0$; the other is $u=\infty$.  Consequently
\[
 2P_{AC}+1
 =\frac14\bigl(q+2(-\tau_1-1)+(-\tau_2-1)\bigr),
\]
which is \eqref{eq:PAC}.
\end{proof}

For $B\cap C$, use the intermediate coordinate $t$ from \eqref{eq:yanhan-B-intermediate}.  Put
\begin{equation*}
 \begin{aligned}
 E&:y^2=t^3-t^2-1,\\
 G_1&:y^2=t^5+t^4+t+1,\\
 G_2&:y^2=t^6+t^4+t^2+1.
 \end{aligned}
\end{equation*}
and write $\tau_E(m),\tau_{G_1}(m),\tau_{G_2}(m)$ for their Frobenius traces.

\begin{lemma}[The $B$--$C$ count]\label{lem:PBC}
For every odd $m$,
\begin{equation}\label{eq:PBC}
 P_{BC}=\frac{q+1+2\tau_E(m)+2\tau_{G_1}(m)+\tau_{G_2}(m)}{16}.
\end{equation}
\end{lemma}

\begin{proof}
Put
\[
 D=t^2-t-1,\qquad B_t=(t+1)(t^2+1),
\]
\[
 A_t=-(t+1)(t^2+t-1).
\]
For $t\ne-1$, these square classes are respectively the discriminant of $v^2-tv+t+1$ and the target classes $1+b$, $1-b$.  Every target in $B\cap C$ has two rational $t$-preimages, each with two rational $v$-preimages.  The affine value $t=-1$ contributes a spurious amount two to the expanded indicator, hence
\begin{equation}\label{eq:PBC-indicator}
 16P_{BC}
 =\sum_{t\in\mathbb F_q}(1+\chi(D))(1+\chi(A_t))(1+\chi(B_t))-2.
\end{equation}
The substitution $t\mapsto1-t$ identifies the $A_t,B_t$ sums and the $DA_t,DB_t$ sums.  Moreover
\[
 A_t=-(t^3-t^2-1),\qquad
 DA_t=-(t^5+t^4+t+1),
\]
\[
 DA_tB_t=-(t+1)^2(t^6+t^4+t^2+1),
\]
and $\sum_t\chi(D)=-1$.

For the remaining product, remove the square factor $(t+1)^2$ and use the genus-one curve $E_0$ of Lemma~\ref{lem:yanhan-smooth-models}.  It has four $\mathbb F_3$-points, so its numerator is $L_{E_0}(T)=1+3T^2$ and its Frobenius trace vanishes in every odd extension.  Its leading coefficient is nonsquare, hence it has no rational point at infinity over such an extension.  Therefore the affine character sum of the squarefree quartic equals one.  At $t=-1$ that quartic has character value $-1$, whereas the original product $A_tB_t$ vanishes because of the removed factor $(t+1)^2$.  Hence the correction is $+1$ and
\[
 \sum_t\chi(A_tB_t)=1+1=2.
\]
Expanding \eqref{eq:PBC-indicator} gives
\[
 16P_{BC}=q-2S_E-2S_{G_1}-S_{G_2}.
\]
For the odd-degree models $E,G_1$, $S_H=-\tau_H(m)$; for the even-degree monic model $G_2$, $S_{G_2}=-\tau_{G_2}(m)-1$.  This is \eqref{eq:PBC}.
\end{proof}

\subsubsection{The triple configuration cover}

Put $b=b_A(u)$ and let $Z$ be the smooth projective normalization of
\begin{equation}\label{eq:Zfield}
 \begin{gathered}
 \mathbb F_3(u,v,z_+,z_-),\qquad b_B(v)=b_A(u),\\
 z_+^2=1+b,\qquad z_-^2=1-b.
 \end{gathered}
\end{equation}
The cyclic quartic pullback has quadratic subfield with square class
$a:=b^2+1=H/A_0^2$.  Put $c_+:=1+b=B_+/A_0$ and
$c_-:=1-b=B_-/A_0$.  The four polynomials $A_0,H,B_+,B_-$ are
squarefree and pairwise coprime over $\mathbb F_3$.  At a geometric
root of $H$, $B_+$, or $B_-$, the valuation parities of
$(a,c_+,c_-)$ are $(1,0,0)$, $(0,1,0)$, or $(0,0,1)$,
respectively.  Thus these square classes are independent over
$\overline{\mathbb F}_3(u)$.  In particular, the quartic pullback
still has degree four, and its quadratic subfield cannot lie in
the biquadratic extension defined by $c_+,c_-$.  As every nontrivial
subfield of a cyclic quartic extension contains the quadratic one,
the two extensions are linearly disjoint.  Their degree-$16$ compositum
remains a field after adjoining $\overline{\mathbb F}_3$, so its full
constant field is $\mathbb F_3$.  Let $\rho$ act by
$v\mapsto1/(v+1)$ while fixing $u,z_+,z_-$, and let
$\epsilon_\pm$ change the sign of $z_\pm$ while fixing the other
generators.  These three deck transformations are defined over
$\mathbb F_3$, so the arithmetic and geometric groups coincide.

\begin{proposition}[The abelian triple cover]\label{prop:Z}
The curve $Z$ is geometrically connected,
\[
 \operatorname{Gal}(Z/\mathbb P^1_u)\cong C_4\times C_2^2,
\]
and $g(Z)=81$.  For every odd $m$,
\begin{equation}\label{eq:ZT}
 \#Z(\mathbb F_{3^m})=32T.
\end{equation}
\end{proposition}

\begin{proof}
Linear disjointness gives the asserted group of order $16$.
Because the degree is prime to three, all ramification is tame.
The inertia subgroups are $\langle\rho\rangle$ at the eight roots
of $H$, $\langle\epsilon_+\rangle$ and
$\langle\epsilon_-\rangle$ at the four roots of $B_+$ and
$B_-$, respectively, and $\langle\epsilon_+\epsilon_-\rangle$
at the four roots of $A_0$; their orders are $4,2,2,2$.
Indeed, $H=0$ pulls back the two totally ramified values of
$b_B$ without multiplicity; at $B_\pm=0$ exactly one target
square root ramifies.  At $A_0=0$, both $c_\pm$ have simple
poles, while $c_+/c_-\to-1$ is a square over the geometric
residue field, giving diagonal inertia.  Infinity is unramified:
$b_A\to0$, $a,c_+,c_-\to1$, and $b=0$ is unramified in
the quartic cover.  Riemann--Hurwitz therefore gives
\[
 2g(Z)-2=-32+8(12)+8(8)+4(8)=160,
\]
so $g(Z)=81$.

For $b\in A\cap B\cap C$, Lemmas~\ref{lem:yanhan-A}, \ref{lem:yanhan-B}, and \ref{lem:yanhan-C} give $2,4,2,2$ rational choices of $u,v,z_+,z_-$, respectively, hence $32$ rational points.  Conversely, every rational nonboundary point of $Z$ gives such a triple target by the block dictionary.

The remaining targets contribute no rational points over an odd extension.  At $b=0$, the numerator of $b_B$ factors over $\mathbb F_3$ as
\[
 -v^4+v^3+v+1=-(v^2+1)(v^2-v-1),
\]
with both quadratics irreducible.  At $b=\pm1$, the $u$-fiber is defined by $u^4\mp u^3-1$, two irreducible quartics over $\mathbb F_3$.  Thus none of these three finite targets has the required rational $u$- or $v$-preimages.

At the poles $u^4=1$ of $b_A$, the ratio
$z_+^2/z_-^2=(1+b)/(1-b)$ has residue $-1$, a nonsquare in every odd extension; consequently the normalization has no rational point above $b=\infty$.  Finally, the branch values $b^2+1=0$ lie in $\mathbb F_9$.  This exhausts the boundary and proves \eqref{eq:ZT}.
\end{proof}

The point count already determines when $\Delta=3$, without computing the zeta function.  Hasse--Weil gives
\begin{equation}\label{eq:yanhan-Weil}
 32T=\#Z(\mathbb F_q)\geqslant q+1-162\sqrt q.
\end{equation}
The right-hand side is positive for odd $m\geqslant11$.  For $m=5,7,9$,
explicit three-root witnesses are encoded and verified by
\nolinkurl{anc/NewDiffSpec_arXiv1_verify.py}.
Thus $T>0$ for every odd $m\geqslant5$.

\subsubsection{Artin characters and explicit quotient/Prym models}

Use complex characters only as labels for the character table; equivalently, fix $\ell\ne2,3$ and a coefficient field in $\overline{\mathbb Q}_\ell$ containing $\mu_4$.  With the generators $\rho,\epsilon_+,\epsilon_-$ fixed above, write
\[
 \chi_{j,e_+,e_-}(\rho)=i^j,\qquad
 \chi_{j,e_+,e_-}(\epsilon_\pm)=(-1)^{e_\pm},
\]
with $j\in\mathbb Z/4\mathbb Z$ and $e_\pm\in\{0,1\}$.  Put
$\psi:=\chi_{1,0,0}$, $\xi_+:=\chi_{0,1,0}$, and
$\xi_-:=\chi_{0,0,1}$; thus $\xi_\pm$ are characters,
whereas $\epsilon_\pm$ are group elements.  The involution
\[
 \iota:(u,v,b,z_+,z_-)\longmapsto(-u,-1/v,-b,z_-,z_+)
\]
is defined over $\mathbb F_3$ and preserves \eqref{eq:Zfield}.
Indeed, $b_A(-u)=-b_A(u)$ and $b_B(-1/v)=-b_B(v)$.
Conjugation sends $\rho$ to $\rho^{-1}$ and exchanges
$\epsilon_+,\epsilon_-$.  Hence it takes $(j,e_+,e_-)$ to
$(-j,e_-,e_+)$ and, commuting with Frobenius, identifies the
corresponding rational single-twist factors.

A \emph{hyperelliptic curve} is a smooth projective curve of genus at least two admitting a degree-two map to $\mathbb P^1$.  The quadratic characters below have explicit such realizations.  Here
\[
 A_0'=u^3,\qquad B_+'=B_-'=u^3,\qquad H'=u^3(1-u^4)
\]
in characteristic three, and direct gcds with the defining polynomials show that $A_0,B_+,B_-,H$ are squarefree and pairwise coprime in every combination used below.  Define
\[
 X:s^2=H.
\]
Then $g(X)=3$, and the relevant quadratic character factors are represented by
\[
 H_1:\ y^2=A_0B_+,
 \qquad
 H_2:\ y^2=B_+B_-,
\]
\[
 Q_+:\ y^2=HA_0B_+,
 \qquad
 Q_0:\ y^2=HB_+B_-.
\]
Thus $H_1,H_2$ have genus $3$ and $Q_+,Q_0$ have genus $7$.

For the order-four pair, set
\[
 B_0=1+u^3-u^4,\qquad t=\frac{B_0-s}{A_0},\qquad D=t^2-t-1.
\]
The relation $b=t(1-t)/(t+1)$ from \eqref{eq:yanhan-B-intermediate}, after $b=b_A(u)$ is substituted, is quadratic in $t$ with discriminant $b_A(u)^2+1=H/A_0^2$.  Thus $\mathbb F_3(X)=\mathbb F_3(u,t)$ is the quadratic intermediate field of $Y/\mathbb P^1_u$, and $Y/X$ is the double cover defined by
\[
 v^2-tv+t+1=0,
\]
whose discriminant is $D$.  On $X$ one has the useful identity
\[
 D=\frac{s(u^3-s)}{A_0^2}.
\]
Twisting this last quadratic extension by the target classes $c_+$ and $c_+c_-$ gives double covers $W_+$ and $W_0$ of $X$.  Let $W_+$ and $W_0$ denote their smooth projective normalizations; their square classes are
\[
 D c_+\sim s(u^3-s)A_0B_+,
 \qquad
 D c_+c_-\sim s(u^3-s)B_+B_-.
\]
The function $s$ has eight simple zeros.  Moreover
\[
 (u^3-s)(u^3+s)=u^6-H=-\,(u^4-1)^2.
\]
At each root of $A_0$, exactly one of $u^3-s$ and $u^3+s$ has a double zero; each has poles of order four at the two points at infinity.  Thus $\operatorname{div}(u^3-s)$ is even.  In particular, $(u^3-s)/A_0^2=-1/(u^3+s)$ gives the normalized $W_0$ square class $B_+B_-$ at $A_0=0$, $s=u^3$.
The branch locus of $W_+\to X$ consists of the eight zeros of $s$ and the eight points over each of $A_0=0$ and $B_+=0$.  For $W_0\to X$, it consists of the eight zeros of $s$ and the eight points over each of $B_+=0$ and $B_-=0$.  These points are distinct by squarefreeness and coprimality.  At either point at infinity both defining square classes have leading term $-u^{16}$ and normalized residue $-1$, hence are unramified.  Thus each cover has exactly $24$ ramification points.  Riemann--Hurwitz gives
\[
 2g(W_\bullet)-2=2(2g(X)-2)+24=32,
 \qquad g(W_\bullet)=17.
\]
The \emph{Prym variety} of each double cover $\pi:W_\bullet\to X$
is $P_\bullet:=(\ker\mathrm{Nm}_\pi)^0$, the identity component
of the norm kernel on Jacobians.  For the covering involution $\tau$,
the identities $\mathrm{Nm}_\pi\pi^*=[2]$ and
$\pi^*\mathrm{Nm}_\pi=1+\tau$ give
$P_\bullet=\operatorname{im}(1-\tau)$ and an
$\mathbb F_3$-isogeny $J(X)\times P_\bullet\longrightarrow J(W_\bullet)$.
These identities apply in characteristic three, since $2$ is
invertible; see \cite[\S\S1--3, pp.~328--332]{Mumford1974}.
Thus $\dim P_\bullet=g(W_\bullet)-g(X)=17-3=14$.
For $\ell\ne2,3$, the idempotents $(1\pm\tau)/2$ split
cohomology, compatibly with Frobenius, as
\[
 H^1(W_\bullet)=H^1(X)\oplus H^1(W_\bullet)^-,
\]
and identify the anti-invariant summand with $H^1(P_\bullet)$.
We denote its numerator by $R_\bullet$.  Equivalently,
\[
 L_{W_\bullet}=L_XR_\bullet,
 \qquad
 s_{R_\bullet}(n)=\#X(\mathbb F_{3^n})-\#W_\bullet(\mathbb F_{3^n}).
\]

\begin{lemma}[Character-to-Prym identification]\label{lem:yanhan-prym-characters}
Let $G=C_4\times C_2^2=\langle\rho\rangle\times\langle\epsilon_+\rangle\times\langle\epsilon_-\rangle$ and put $K=\mathbb F_3(u)=\mathbb F_3(Z)^G$.  The fixed field of $\langle\epsilon_+,\epsilon_-\rangle$ is the cyclic quartic field $\mathbb F_3(Y)$, while
\[
 \mathbb F_3(X)=\mathbb F_3(Z)^{\langle\rho^2,\epsilon_+,\epsilon_-\rangle}
\]
is its unique quadratic intermediate over $K$; equivalently $X/K$ is cut out by the character $\psi^2$.  The residual double cover $Y/X$ has deck involution induced by $\rho^2$, and its anti-invariant cohomology realizes the rational conjugate pair $\{\psi,\psi^3\}$.

Twisting $Y/X$ by $c_+$, $c_-$, and $c_+c_-$ gives, respectively, the character pairs $\{\psi,\psi^3\}\xi_+$, $\{\psi,\psi^3\}\xi_-$, and $\{\psi,\psi^3\}\xi_+\xi_-$.  The Prym numerators $R_+$ and $R_0$ are the corresponding rational Artin factors; the $+$ and $-$ factors agree under $u\mapsto-u$.
\end{lemma}

\begin{proof}
Put $J=\langle\rho^2,\epsilon_+,\epsilon_-\rangle$.
The defining square classes $Dc_+$ and $Dc_+c_-$ give characters
$\eta_+,\eta_0$ of $J$ with values $(-1,-1,1)$ and
$(-1,-1,-1)$ on these generators.  Their extensions to $G$ give
\[
 \operatorname{Ind}_{J}^{G}\eta_+
  =\psi\xi_+\oplus\psi^3\xi_+,
 \qquad
 \operatorname{Ind}_{J}^{G}\eta_0
  =\psi\xi_+\xi_-\oplus\psi^3\xi_+\xi_-.
\]
For inertia generators $\rho,\epsilon_+,\epsilon_-$, and
$\epsilon_+\epsilon_-$, the invariant dimensions of these two
induced representations are respectively $(0,0,2,0)$ and
$(0,0,0,2)$.  The induction identity holds also at the ramified places:
for a closed place $v$ of $K$ with $I_v\subseteq J$ and
$\eta_\bullet|_{I_v}=1$, its local determinant is
\[
 \begin{cases}
 (1-\eta_\bullet(F_v)T^{\deg v})^2,&F_v\in J,\\
 1-\eta_\bullet(F_v^2)T^{2\deg v},&F_v\notin J,
 \end{cases}
\]
where $F_v$ is a Frobenius lift.  In the other cases the determinant
is one, since $\eta_\bullet(\rho^2)=-1$ when
$I_v=\langle\rho\rangle$.  These are precisely the products of
the local factors of the quadratic character over the places of $X$
above $v$.  Consequently,
\[
 R_\bullet(T)=L_{W_\bullet}(T)/L_X(T)
            =L_K(T,\operatorname{Ind}_{J}^{G}\eta_\bullet).
\]
At $B_-=0$, $A_0=u^3$ and $s=\pm\zeta u^3$
($\zeta\in\mathbb F_9$, $\zeta^2=-1$),
so the $W_+$ radicands lie in $u^{12}\mathbb F_9^\times$
and are squares over $\mathbb F_{81}$.  The $\operatorname{Ind}_J^G\eta_+$ local
determinant is $(1-T^4)^2$.  At $A_0=0$, the normalized
$W_0$ classes $-1,-u^6$ give $(1+T)^2$ at $u=\pm1$
and $(1-T^2)^2$ at $u^2+1=0$.
For the twist exchange, put
$r=(u^3+s)/A_0$.  On $X$, the identity
$(u^3-s)(u^3+s)=-A_0^2$ gives
$D(-u,s)=r^2D(u,s)$ and $r(-u,s)r(u,s)=1$.
Thus $(u,s,w)\mapsto(-u,s,rw)$ identifies
$w^2=Dc_+$ with $w^2=Dc_-$ over $\mathbb F_3$.
This function-field isomorphism extends to the smooth projective
models and commutes with Frobenius, proving equality of their
Prym factors.
\end{proof}

\begin{table}[H]
\caption{Rational Artin pieces for the Yan--Han triple cover.}\label{tab:yanhan-artin-dictionary}
\centering\footnotesize
\begin{tabularx}{\textwidth}{@{}>{\raggedright\arraybackslash}p{0.19\textwidth}p{0.28\textwidth}p{0.19\textwidth}X@{}}
\toprule
Characters & Explicit realization & Degree of numerator & Source of the degree \\
\midrule
$j=1,2,3$, no target twist & cyclic quartic pullback $Y$ & $18$ & $g(Y)=9$ \\
$\xi_+$ or $\xi_-$ & $H_1:y^2=A_0B_+$, and $u\mapsto-u$ & $6$ each & genus $3$ \\
$\xi_+\xi_-$ & $H_2:y^2=B_+B_-$ & $6$ & genus $3$ \\
$\psi^2\xi_+$ or $\psi^2\xi_-$ & $Q_+:y^2=HA_0B_+$, and $u\mapsto-u$ & $14$ each & genus $7$ \\
$\psi^2\xi_+\xi_-$ & $Q_0:y^2=HB_+B_-$ & $14$ & genus $7$ \\
$\{\psi,\psi^3\}\xi_+$ or $\{\psi,\psi^3\}\xi_-$ & Prym$(W_+/X)$ & $28$ each & Prym dimension $14$ \\
$\{\psi,\psi^3\}\xi_+\xi_-$ & Prym$(W_0/X)$ & $28$ & Prym dimension $14$ \\
\bottomrule
\end{tabularx}
\end{table}

\begin{samepage}
\begin{proposition}[Rational character decomposition]\label{prop:yanhan-character-decomp}
Let $L_Y,L_{H_1},L_{H_2},L_{Q_+},L_{Q_0},R_+,R_0$ denote the numerator polynomials attached to the indicated curve or Prym pieces.  Then
\begin{equation}\label{eq:Artin-natural}
 L_Z=L_YL_{H_1}^2L_{H_2}L_{Q_+}^2L_{Q_0}R_+^2R_0.
\end{equation}
The degrees in Table~\ref{tab:yanhan-artin-dictionary}, counted with the multiplicities in \eqref{eq:Artin-natural}, sum to $162=2g(Z)$.
\end{proposition}

\begin{proof}
At the loci $H=0$, $B_+=0$, $B_-=0$, and $A_0=0$, the tame
inertia groups are generated by $\rho$, $\epsilon_+$, $\epsilon_-$,
and $\epsilon_+\epsilon_-$, respectively.  Thus the conductor degree
of any nontrivial $\chi_{j,e_+,e_-}$ is
\[
 8\mathbf1_{j\ne0}+4e_++4e_-+4(e_+\oplus e_-),
\]
where $\oplus$ denotes addition modulo two.  The rank-one
Grothendieck--Ogg--Shafarevich formula on $\mathbb P^1$ gives
$\deg L_{\chi}=\deg\operatorname{cond}(\chi)-2$:
six characters have degree $6$, and nine have degree $14$.
These are exactly the $15$ characters grouped in
Table~\ref{tab:yanhan-artin-dictionary}, with the order-four pairs
identified by Lemma~\ref{lem:yanhan-prym-characters}.  Artin
factorization gives \eqref{eq:Artin-natural}, of total degree
$6\cdot6+9\cdot14=162$.
\end{proof}
\end{samepage}

\subsubsection{Fixed Weil polynomials and exact reconstruction}

For
\[
 L(T)=1+a_1T+\cdots+a_rT^r=\prod_{j=1}^r(1-\alpha_jT),
\]
write $s_L(n):=\sum_j\alpha_j^n$.  Newton identities recover the first half of a weight-one Weil numerator from its first $g$ power sums, and the functional equation recovers the second half; for $n>r$,
\[
 s_L(n)+a_1s_L(n-1)+\cdots+a_rs_L(n-r)=0.
\]

\begin{definition}[Fixed Yan--Han trace polynomials]\label{def:yanhan-weil}
Let
\[
 L_0:=1+3T^2,\quad L_1:=1-2T+3T^2,\quad L_2:=1+2T+3T^2,
\]
\[
 L_3:=1-2T^2+9T^4,
\quad L_4:=1+2T+2T^2+6T^3+9T^4,
\]
\[
 L_5:=1-2T+4T^2-6T^3+9T^4,
\]
\[
 L_6:=1-2T+2T^2-6T^3+9T^4.
\]
\begin{align*}
 L_7:={}&1+2T+3T^2+4T^3\\
      &+9T^4+18T^5+27T^6,
\end{align*}
\begin{align*}
 L_8:={}&1+4T^3-T^4+4T^5+16T^6\\
      &+12T^7-9T^8+108T^9+729T^{12},
\end{align*}
\begin{align*}
 L_9:={}&1+2T+2T^2-2T^3-11T^4-16T^5\\
 &+8T^6+32T^7+50T^8-60T^9-148T^{10}\\
 &-180T^{11}+450T^{12}+864T^{13}+648T^{14}\\
 &-3888T^{15}-8019T^{16}-4374T^{17}\\
 &+13122T^{18}+39366T^{19}+59049T^{20}.
\end{align*}
Define
\begin{equation*}
 \Theta_m:=6s_{L_2}(m)+4s_{L_5}(m)+4s_{L_7}(m)+2s_{L_6}(m)-2s_{L_1}(m)
\end{equation*}
and
\begin{align*}
 \Xi_m:={}&3s_{L_0}(m)+6s_{L_1}(m)+6s_{L_2}(m)\\
 &+2s_{L_3}(m)+3s_{L_4}(m)+4s_{L_5}(m)\\
 &+5s_{L_6}(m)+2s_{L_7}(m)+2s_{L_8}(m)+2s_{L_9}(m).
\end{align*}
\end{definition}

\begin{proposition}[Frobenius and Artin factorization]\label{prop:yanhan-Artin}
The raw point counts of the explicit pair-intersection curves and of the Artin/Prym pieces above reconstruct the numerator factors as
\[
 L_Y=L_2^3L_4L_5^2,\qquad
 L_{H_1}=L_7,\qquad
 L_{H_2}=L_2L_6,
\]
\[
 L_{Q_+}=L_2L_8,\qquad
 L_{Q_0}=L_0^3L_3^2,
\]
\[
 L_E=L_1,\qquad L_{G_1}=L_4,\qquad L_{G_2}=L_2^2,
\]
\[
 R_+=L_1^2L_6L_9,\qquad
 R_0=L_1^2L_6^2L_5^2L_4^2.
\]
Consequently
\begin{equation}\label{eq:LZ}
 L_Z=L_0^3L_1^6L_2^6L_3^2L_4^3L_5^4L_6^5L_7^2L_8^2L_9^2.
\end{equation}
In particular, the displayed $L_i$ are genuine $3$-Weil factors, and
\begin{equation}\label{eq:TXi}
 T=\frac{q+1-\Xi_m}{32}.
\end{equation}
\end{proposition}

\begin{proof}
For a genus-$g$ curve, its first $g$ exact projective counts determine
the Weil numerator by Newton identities and the functional equation.
The first fourteen relative traces
$\#X(\mathbb F_{3^n})-\#W_\bullet(\mathbb F_{3^n})$
likewise determine the degree-$28$ Prym numerators.  The single
standard-library script \nolinkurl{anc/NewDiffSpec_arXiv1_verify.py}
recomputes these counts from the displayed equations: the
\texttt{--raw-curves} option counts the nine curve models through
their respective genera, and \texttt{--raw-prym --max-degree 14}
evaluates all fourteen relative traces.  The latter uses the
normalized character sums for $W_\bullet\to X$, including the
even-order zero at $u^4=1$ and both points at infinity.
Newton reconstruction and exact polynomial multiplication then give
the displayed factors and \eqref{eq:LZ}.  The stored counts also
allow a shorter default verification.  The \texttt{--historical}
option independently enumerates the five archived derivative rows
at $m=3,5,7,9,11$ and the six $Z$-counts at $m=1,3,5,7,9,11$,
including the projective boundary corrections.  These provide
consistency checks separate from the Weil-factor reconstruction.
The geometric identifications follow from
Proposition~\ref{prop:yanhan-character-decomp} and
Lemma~\ref{lem:yanhan-prym-characters}; purity gives reciprocal-root
modulus $\sqrt3$.  Finally, \eqref{eq:ZT} implies \eqref{eq:TXi}.
\end{proof}

Combining \eqref{eq:PAB}, \eqref{eq:PAC}, \eqref{eq:PBC} with the reconstructed pair factors gives
\begin{equation}\label{eq:PTheta}
 P=\frac{5q-27-\Theta_m}{16}.
\end{equation}

\subsubsection{Complete spectrum and exact uniformity}

The bound in \cite[Theorem~1]{YanHan2019} remains valid at every
odd $m\geqslant3$.  The equality conjectured immediately after that
theorem \cite[Section~3, Conjecture]{YanHan2019} has a planar
exception at $m=3$.  The corollary below proves equality for all
remaining odd $m\geqslant5$ and determines the asymptotic distribution.

\begin{theorem}[Complete Yan--Han spectrum]\label{thm:yanhan-spectrum}
Let $m\geqslant3$ be odd, $q=3^m$, and $d=(q-3)/4$.  With $\Theta_m,\Xi_m$ as above, the differential spectrum of $x^d$ is
\[
 \omega_0=\frac{9q-55-2\Theta_m+\Xi_m}{32},
\]
\[
 \omega_1=\frac{15q+111+4\Theta_m-3\Xi_m}{32},
\]
\[
 \omega_2=\frac{7q-57-2\Theta_m+3\Xi_m}{32},
 \qquad
 \omega_3=\frac{q+1-\Xi_m}{32},
\]
and $\omega_i=0$ for $i\geqslant4$.
\end{theorem}

\begin{proof}
Substitute \eqref{eq:PTheta} and \eqref{eq:TXi} into Proposition~\ref{prop:yanhan-reconstruction}.
\end{proof}

\begin{corollary}[Uniformity and asymptotic distribution]\label{cor:yanhan-asymptotic}
For $m=3$, the derivative is the permutation $1-x^3$ and
\[
 (\omega_0,\omega_1,\omega_2,\omega_3)=(0,27,0,0).
\]
For every odd $m\geqslant5$, the differential uniformity is exactly three.  Moreover
\[
 (\omega_0,\omega_1,\omega_2,\omega_3)
 =\frac q{32}(9,15,7,1)+O(\sqrt q).
\]
\end{corollary}

\begin{proof}
The case $m=3$ follows from $(x+1)^6-x^6=1-x^3$.  For odd $m\geqslant5$, Theorem~\ref{thm:yanhan-spectrum} already shows that no fiber size exceeds three.  Positivity of $T=\omega_3$ follows from \eqref{eq:yanhan-Weil} for $m\geqslant11$ and from the explicit small-field witnesses for $m=5,7,9$, so a three-fiber actually occurs and hence $\Delta=3$.  Proposition~\ref{prop:yanhan-Artin} shows that every nontrivial trace in $\Theta_m$ or $\Xi_m$ is a sum of weight-one Frobenius eigenvalues of absolute value $\sqrt3$, giving the asymptotic formula.
\end{proof}

\section{The Felke Family: Spectra from Character Projections}\label{sec:felke}
Let
\[
 p\equiv1\pmod4,\qquad n=2m-1\geqslant1,\qquad q=p^n,
\]
and put
\[
 \sigma:=p^m,\qquad h:=\frac12,\qquad
 d:=\frac{q-1}{2}+\sigma+1.
\]
Let $\chi$ be the quadratic character of $\mathbb{F}_q$, extended by $\chi(0)=0$. Since $n$ is odd, its restriction to $\mathbb{F}_p$ is the quadratic character $\chi_p$. On $\mathbb{F}_q^\times$,
\[
 F(x):=x^d=\chi(x)x^{\sigma+1}.
\]
Center the derivative:
\[
 \Phi(t):=F(t+h)-F(t-h),\qquad t^*:=t^\sigma.
\]
For $t\ne\pm h$, the four sign pairs $\bigl(\chi(t-h),\chi(t+h)\bigr)$ give
\begin{equation}\label{felke:eq:block-values}
 \Phi(t)=
 \begin{cases}
 t+t^*, &(+,+),\\
 2tt^*+h,&(-,+),\\
 -2tt^*-h,&(+,-),\\
 -t-t^*,&(-,-).
 \end{cases}
\end{equation}
These are Felke's four fundamental cases
\cite[Section~4.1]{Felke2020}.  The injectivity and root restrictions
used below are already present in his Sections~4.2--4.3; we express
them as explicit character projectors to count all target
multiplicities.  We retain his exact uniformity result, including
$(p,n)=(17,1)$, and do not treat the branches $p\equiv3\pmod4$.

\subsection{The three marked blocks outside the prime field}
Set
\[
 L(u):=u+u^*,\qquad Q(v):=2vv^*+h.
\]

\begin{lemma}[The linearized block]\label{felke:lem:L}
The map $L:\mathbb{F}_q\to\mathbb{F}_q$ is a permutation and
\begin{equation}\label{felke:eq:Linv}
 L^{-1}(b)=\frac12\sum_{j=0}^{n-1}(-1)^j b^{\sigma^j}.
\end{equation}
If $b\notin\mathbb{F}_p$ and $u=L^{-1}(b)$, then the two same-sign cases in \eqref{felke:eq:block-values} contribute exactly one solution to the fiber over $b$ if and only if
\[
 \chi(u^2-h^2)=1.
\]
\end{lemma}

\begin{proof}
Since $*$ has odd order $n$,
\[
 (1+*)\sum_{j=0}^{n-1}(-*)^j=2,
\]
giving \eqref{felke:eq:Linv}.  For $b\notin\mathbb F_p$,
$u\ne\pm h$.  The only same-sign candidates are $t=u$
in $(+,+)$ and $t=-u$ in $(-,-)$.  Their endpoint
characters agree exactly when $\chi(u^2-h^2)=1$;
as $\chi(-1)=1$, precisely one candidate then has
the required sign.
\end{proof}

\begin{lemma}[The mixed block quotient]\label{felke:lem:mixed}
One has
\[
 \gcd(\sigma+1,q-1)=2.
\]
Consequently $v\mapsto vv^*=v^{\sigma+1}$ maps $\mathbb{F}_q^\times$ two-to-one onto the nonzero squares. Define
\begin{equation*}
 \kappa:=\frac{2\sigma-q-1}{p-1}\in\mathbb Z.
\end{equation*}
Then
\begin{equation}\label{felke:eq:kappa-cong}
 (\sigma+1)\kappa\equiv2\pmod{q-1}.
\end{equation}
Hence, for $v\ne0$, if $s=vv^*$, then $v^2=s^\kappa$.

For $b\notin\mathbb{F}_p$ put $s=(b-h)/2$. The $(-,+)$ block contributes exactly one solution to the fiber over $b$ if and only if
\begin{equation*}
 \chi(s)=1,\qquad \chi(s^\kappa-h^2)=-1.
\end{equation*}
The $(+,-)$ block at $b$ is the same condition applied to $-b$.
\end{lemma}

\begin{proof}
Let $g=\gcd(\sigma+1,q-1)$.  Since $p^m\equiv-1$ and
$p^{2m-1}\equiv1\pmod g$, we have $p\equiv1\pmod g$,
so $g\mid2$; both arguments are even, hence $g=2$.
As $\sigma\equiv q\equiv1\pmod{p-1}$, $\kappa$ is integral, and
\[
 (\sigma+1)\kappa-2=\frac{q-1}{p-1}(2p-\sigma-1),
\]
proving \eqref{felke:eq:kappa-cong}.  For $s\ne0$, the equation
$v^{\sigma+1}=s$ has exactly two roots $\pm v$ when
$\chi(s)=1$, with $v^2=s^\kappa$.  Their endpoint characters
are opposite iff $\chi(v^2-h^2)=-1$.  Since $\chi(-1)=1$,
negation exchanges the two mixed cases; exactly one root
belongs to $(-,+)$.  The $(+,-)$ case follows by $b\mapsto-b$.
\end{proof}

For $s\in\mathbb{F}_q\setminus\mathbb{F}_p$ set
\[
 b_s:=2s+h,\qquad U_s:=L^{-1}(b_s),\qquad s^\vee:=-s-h.
\]
The arguments of every character below are nonzero. For example, if $s^\kappa=h^2$, then raising to the $(\sigma+1)$-st power and using \eqref{felke:eq:kappa-cong} gives $s^2=h^4$, hence $s=\pm h^2\in\mathbb{F}_p$, a contradiction. Define
\begin{align*}
 A(s)&:=\frac{1+\chi(U_s^2-h^2)}2,\\
 B(s)&:=\frac{(1+\chi(s))(1-\chi(s^\kappa-h^2))}{4},\\
 C(s)&:=B(s^\vee).
\end{align*}
Thus $A(s),B(s),C(s)\in\{0,1\}$.

\begin{samepage}
\begin{proposition}[Targetwise fiber formula outside $\mathbb{F}_p$]\label{felke:prop:targetwise}
For every $s\in\mathbb{F}_q\setminus\mathbb{F}_p$,
\begin{equation*}
 \delta_F(1,b_s)=A(s)+B(s)+C(s).
\end{equation*}
In particular every fiber over a target outside $\mathbb F_p$ has size at most three.
\end{proposition}

\begin{proof}
The four cases in \eqref{felke:eq:block-values} partition the nonboundary inputs. By Lemma~\ref{felke:lem:L}, the two same-sign cases together contribute $A(s)$. By Lemma~\ref{felke:lem:mixed}, the $(-,+)$ case contributes $B(s)$. Since $b_{s^\vee}=-b_s$, the $(+,-)$ case contributes $B(s^\vee)=C(s)$.
\end{proof}
\end{samepage}

\subsection{Prime-field targets}
The following closure property is proved in
\cite[Section~4.4]{Felke2020}; we include its short argument in the
centered coordinates used here.
\begin{lemma}[Prime-field closure]\label{felke:lem:primeclosure}
If $\Phi(t)\in\mathbb{F}_p$, then $t\in\mathbb{F}_p$.
\end{lemma}

\begin{proof}
In a same-sign block $L(t)=\pm\Phi(t)\in\mathbb{F}_p$. Applying $*$ gives $t^{\sigma^2}=t$, hence $t^p=t$ since $\sigma^2=pq$. In a mixed block $tt^*\in\mathbb{F}_p$; applying $*$ gives $t^*(t^p-t)=0$, so again $t\in\mathbb{F}_p$. The endpoints $\pm h$ lie in $\mathbb{F}_p$.
\end{proof}

Thus the prime-field part is exactly the derivative of
\[
 g_p(x):=\chi_p(x)x^2=x^{(p+3)/2}
 \quad\text{on }\mathbb{F}_p.
\]
Put $\varepsilon:=\chi_p(2)$ and define
\begin{align}
 P_0&:=\#\left\{b\in\mathbb{F}_p:
 \chi_p(b-1)=\chi_p(b+1)=-\varepsilon,
 \ \chi_p(2b-1)=1\right\},\label{felke:eq:P0}\\
 Q_0&:=\#\left\{b\in\mathbb{F}_p:
 \chi_p(b-1)=\chi_p(b+1)=-\varepsilon,
 \ \chi_p(2b-1)=\chi_p(2b+1)=1\right\}.
 \label{felke:eq:Q0}
\end{align}

The spectrum of $x^{(Q+3)/2}$ over $\mathbb F_Q$, for
$Q\equiv1\pmod4$ in characteristic different from three, was
determined in \cite[Theorem~2]{JiangLiLiQu2022}.  The following
$Q=p$ specialization is written in the intersection notation needed
for the extension-field calculation.

\begin{proposition}[Prime-field spectrum contribution]\label{felke:prop:base}
The contribution of the $p$ prime-field targets to the spectrum is
\begin{equation}\label{felke:eq:base-spectrum}
 (\beta_0,\beta_1,\beta_2,\beta_3)
 =\bigl(2P_0,\ p-4P_0+Q_0,\ 2P_0-2Q_0,\ Q_0\bigr).
\end{equation}
\end{proposition}

\begin{proof}
For $t\in\mathbb F_p\setminus\{\pm h\}$ the four block values are
$\pm2t$ and $\pm(2t^2+h)$.  At target $b$, the same-sign
block contributes when $\chi_p(b^2-1)=1$, and the two mixed
blocks contribute under, respectively,
\[
\begin{aligned}
 &\chi_p(2b-1)=1,\quad\chi_p(b-1)=-\varepsilon,\\
 &\chi_p(2b+1)=1,\quad\chi_p(b+1)=-\varepsilon.
\end{aligned}
\]
Each mixed image meets the same-sign image in $P_0$ targets,
by \eqref{felke:eq:P0} and $b\mapsto-b$.  The mixed images
intersect in exactly $Q_0$ targets, all in the same-sign image.
At $b=\pm1$ the three regular indicators vanish, while the
boundary inputs $t=\pm h$ contribute one point each.  Hence
there are $2P_0+Q_0$ unordered pairs and $Q_0$ unordered
triples in prime-field fibers.  These $p$ targets have $p$
preimages in total, so factorial-moment inversion gives
\eqref{felke:eq:base-spectrum}.
\end{proof}

The same prime-field contribution can be evaluated through two
elliptic character sums.  Put
\[
 \Gamma_1:=\sum_{x\in\mathbb F_p}\chi_p\bigl(x(x-1)(x+3)\bigr),
 \qquad
 \Gamma_2:=\sum_{x\in\mathbb F_p}\chi_p\bigl(x(x-1)(x+2)\bigr).
\]
Specializing \cite[Theorem~2]{JiangLiLiQu2022} to $Q=p$ and
comparing its zero- and three-fiber frequencies with
\eqref{felke:eq:base-spectrum} gives
\begin{equation}\label{felke:eq:base-elliptic}
 P_0=\frac{p+\Gamma_1+1}{8},\qquad
 Q_0=\frac{p+3\Gamma_1-2\Gamma_2+1}{16}.
\end{equation}
Here $\Gamma_i$ is the negative Frobenius trace of the corresponding
smooth cubic model.  Formulas \eqref{felke:eq:base-elliptic} evaluate the known prime-field
part; the remaining extension-field terms are the sums below.

\subsection{The complete spectrum in odd extensions}
Define the three intersection sums
\begin{equation}\label{felke:eq:PQT}
 \begin{aligned}
 P&:=\sum_{s\in\mathbb{F}_q\setminus\mathbb{F}_p}A(s)B(s),\\
 Q&:=\sum_{s\in\mathbb{F}_q\setminus\mathbb{F}_p}B(s)B(s^\vee),\\
 T&:=\sum_{s\in\mathbb{F}_q\setminus\mathbb{F}_p}A(s)B(s)B(s^\vee).
 \end{aligned}
\end{equation}
Each summand is explicit through the quadratic-character formulas for $A$ and $B$ above.

\begin{theorem}[Complete differential spectrum]\label{felke:thm:main}
Let $p\equiv1\pmod4$ be prime and $n=2m-1$ be odd. For
\[
 F(x)=x^{(p^n-1)/2+p^m+1}\quad\text{on }\mathbb{F}_{p^n},
\]
let $P_0,Q_0$ be \eqref{felke:eq:P0}--\eqref{felke:eq:Q0} and $P,Q,T$ be \eqref{felke:eq:PQT}. Then the complete differential spectrum is
\begin{equation}\label{felke:eq:main-spectrum}
 \begin{aligned}
 \omega_0&=2P_0+2P+Q-T,\\
 \omega_1&=q-4P_0+Q_0-4P-2Q+3T,\\
 \omega_2&=2P_0-2Q_0+2P+Q-3T,\\
 \omega_3&=Q_0+T,
 \end{aligned}
\end{equation}
and $\omega_j=0$ for $j\geqslant4$. The formula includes $(p,n)=(17,1)$, where $\omega_3=0$.
\end{theorem}

\begin{proof}
For $s\notin\mathbb F_p$, $\delta=A+B+C$.  Under
$s\mapsto s^\vee$, $A$ is invariant and $B,C$ are exchanged;
hence $\sum(AB+AC+BC)=2P+Q$ and $\sum ABC=T$.
Lemma~\ref{felke:lem:primeclosure} gives $q-p$ inputs and
targets in this sector.  Binomial inversion for the three
binary images (Proposition~\ref{prop:block-incidence}, $h=1$),
followed by addition of \eqref{felke:eq:base-spectrum}, gives
\eqref{felke:eq:main-spectrum}.
\end{proof}

\begin{corollary}[Moment form]\label{felke:cor:moments}
The two nontrivial factorial moments are
\[
 M_2=4P_0+2Q_0+4P+2Q,\qquad
 M_3=6(Q_0+T),
\]
and equivalently
\[
 \omega_3=\frac{M_3}{6},\quad
 \omega_2=\frac{M_2-M_3}{2},\quad
 \omega_1=q-M_2+\frac{M_3}{2},\quad
 \omega_0=\frac{M_2}{2}-\frac{M_3}{6}.
\]
\end{corollary}

\section{Geometric Comparisons and Limitations}\label{sec:elliptic}\label{sec:scope}

\subsection{A known elliptic configuration curve}
Let $p\geqslant5$, $q=p^n$, and let $F(0)=0$, $F(x)=x^{q-3}$ for $x\ne0$.  After the centered change $y=2x+1$, the rational differential core is
\[
 \phi(y)=-\frac{16y}{(y^2-1)^2}.
\]
The complete differential spectrum was determined in~\cite{YanEtAl2022}; the following proposition identifies its second-moment configuration curve.

\begin{proposition}[Intrinsic elliptic double configuration]\label{prop:quartic-monodromy}\label{prop:elliptic-birational}
The geometric and arithmetic monodromy groups of $\phi$ are $S_4$.  The non-diagonal ordered double-configuration curve is geometrically irreducible of genus one and is birational to
\[
 \mathcal C:\qquad Y^2=(X^2-4)(-3X^2-4).
\]
Consequently, the elliptic Frobenius trace occurring in the known differential-spectrum formulas is the intrinsic weight-one contribution of the configuration curve required by the second factorial moment.
\end{proposition}

\begin{proof}
The branch cycles have types $(3,1)$, $(2,2)$, $(2,1,1)$, $(2,1,1)$, which force $S_4$; their action on ordered pairs of distinct letters is transitive, and Riemann--Hurwitz gives genus one.  Cancelling the diagonal in $\phi(y)-\phi(z)$ and using $a=y+z$, $c=yz$, $w=y-z$ yields $ca^2=(c+1)^2$.  On the dense open set $ar\ne0$, put $r=(c+1)/a$; then $r^2=c$ and $a=r+r^{-1}$.  The substitutions
\[
 X=2r,\qquad Y=4rw
\]
give
\[
 Y^2=(X^2-4)(-3X^2-4).
\]
Conversely, on the dense open set $X\ne0$, recover $r=X/2$, $c=r^2$, $a=r+r^{-1}$, $w=Y/(4r)$, and then $y=(a+w)/2$, $z=(a-w)/2$.  Thus the two curves are birational; the omitted points belong only to the projective boundary of these dense opens and do not affect the genus calculation.
\end{proof}

\subsection{Scope of the geometric reductions}

After cancellation, projective completion, and boundary corrections, the moment identities apply to any nonconstant separable row $x\mapsto f(x+a)-cf(x)$.  For a rational function $f$ the two-parameter map is
\[
 \begin{aligned}
 \Psi_{f,c}:\;&\mathbb A^1_x\times\mathbb G_{m,a}
 \dashrightarrow \mathbb G_{m,a}\times\mathbb A^1_b,\\
 &(x,a)\longmapsto(a,f(x+a)-cf(x)).
 \end{aligned}
\]
The monodromy over $k(a,b)$ need not specialize unchanged: degree, branch type, separability, or the constant field may vary with $a$.  Scaling identifies all nonzero directions for power maps, not for general rational functions.  Proposition~\ref{prop:block-incidence} requires an injective quotient on each input block, and Kummer transfer alone does not evaluate point counts.  Even for the same rational core, arithmetic depth can depend on the characteristic: inverse-cubic powers require no third moment when $p=5$.

\subsection{Semilinear and toric reductions}

A Frobenius-dependent reduction can give a \emph{projective semilinear map} $J_b=M_b\circ\operatorname{Frob}^r$, where $M_b$ is a M\"obius map and $\operatorname{Frob}^r(x)=x^{p^r}$, instead of a fixed-base cover.  For DMMP, Theorem~\ref{dmmp6:thm:semilinear-dickson} uses M\"obius and Dickson transformations to evaluate the two character correlations of $z^\sigma+z=w^{\sigma+1}$ by the same elliptic trace as in the fixed-quartic case; no general semilinear configuration theorem follows.

For $K=\mathbb F_{q^2}$ over $k=\mathbb F_q$, the group $U:=\{u\in K^\times:u^{q+1}=1\}$ is the group of $k$-points of the norm-one torus $\operatorname{Res}^1_{K/k}\mathbb G_m$.  The homomorphism $K^\times\to U$, $x\mapsto x^{q-1}$, is surjective.  Such substitutions occur for exponents $r(q-1)$ with integer $r$ \cite{HuEtAl2023} and for \emph{normalized Niho exponents}, of the form $d=s(q-1)+1$ over $\mathbb F_{q^2}$ with $q=2^m$ \cite{XieEtAl2023}; the associated incidence questions for translates of $U$ still require case-specific normalization, boundary corrections, and descent.

\section{Conclusion}

Once fiber support and exceptional targets are controlled, factorial
moments or intersections of injective block images recover exact
spectra.  The five applications give explicit conic and elliptic
evaluations, a genus-$81$ Artin decomposition for Yan--Han, and
character-sum formulas for the Felke extension-field terms.  For YDQZ,
the genus-$105$ $B$--$C$ cover has an explicit $D_4\times D_4$
decomposition whose odd-degree traces reduce to genera one and two;
the other two pair covers and the triple cover remain to be evaluated.
The proofs keep the exceptional fibers and the arithmetic constant
fields explicit.  A single standard-library script reproduces the
finite arithmetic used in the factorizations and independently checks
the archived derivative and configuration counts.

\newcommand{\etalchar}[1]{$^{#1}$}
\providecommand{\bysame}{\leavevmode\hbox to3em{\hrulefill}\thinspace}
\providecommand{\MR}{\relax\ifhmode\unskip\space\fi MR }
\providecommand{\MRhref}[2]{%
  \href{http://www.ams.org/mathscinet-getitem?mr=#1}{#2}
}
\providecommand{\href}[2]{#2}

\end{document}